\documentclass[review]{elsarticle}
\usepackage{amsmath,amssymb}
\usepackage{lineno,hyperref}
\usepackage{multirow}
\usepackage{array}
\usepackage{threeparttable}
\usepackage{siunitx}
\usepackage{empheq}
\usepackage{subcaption}
\usepackage{rotating}
\usepackage{apptools}
\usepackage{longtable}
\usepackage{tabularx}
\usepackage{booktabs}

\usepackage[heightadjust=all]{floatrow}
\counterwithin*{equation}{section}

\newtheorem{Remark}{Remark}[section]
\newtheorem{Proposition}[Remark]{Proposition}
\newtheorem{Theorem}[Remark]{Theorem}

\newtheorem{Proof}{Proof}
\newtheorem{Propositionn}{Proposition}

\AtAppendix{
}
	
\journal{Journal of \LaTeX\ Templates}

\begin{document}

\begin{frontmatter}

\title{Phase field models for thermal fracturing and their variational structures \tnoteref{mytitlenote}}

\author[mymainaddress,mysecondaryaddress]{S. Alfat\corref{mycorrespondingauthor}}
\cortext[mycorrespondingauthor]{Corresponding author}
\ead{sayahdin.alfat@yahoo.com}

\author[mythirdaddress]{M. Kimura}
\ead{mkimura@se.kanazawa-u.ac.jp}

\author[mysecondaryaddress]{Alifian. M. M.}
\ead{fiansmamda@gmail.com}

\address[mymainaddress]{Physics Education Department, Halu Oleo University, Kendari, Southeast Sulawesi, Indonesia}
\address[mysecondaryaddress]{Division of Mathematical and Physical Sciences, Graduate School on Natural Science and Technology, Kanazawa University, Kakuma, Kanazawa, 920-1192, Japan}
\address[mythirdaddress]{Faculty of Mathematics and Physics, Kanazawa University, Kakuma, Kanazawa, 920-1192, Japan}

\begin{abstract}
It is often observed that thermal stress enhances crack propagation in materials, and conversely, crack propagation can contribute to temperature shifts in materials. In this study, we first consider the thermoelasticity model proposed by M. A. Biot (1956) and study its energy dissipation property. The Biot thermoelasticity model takes into account the following effects. Thermal expansion and contraction are caused by temperature changes, and conversely, temperatures decrease  in expanding areas but increase in contracting areas. In addition, we examine its thermomechanical properties through several numerical examples and observe that the stress near a singular point is enhanced by the thermoelastic effect. In the second part, we propose two crack propagation models under thermal stress by coupling a phase field model for crack propagation and the Biot thermoelasticity model and show their variational structures. In our numerical experiments, we investigate how thermal coupling affects the crack speed and shape. In particular, we observe that the lowest temperature appears near the crack tip, and the crack propagation is accelerated by the enhanced thermal stress.
\end{abstract}

\begin{keyword}
Thermoelasticity \sep Crack Propagation \sep Crack Path \sep Phase Field Model \sep Variational Structure \sep Energy Equality \sep Adaptive Finite Element Method
\end{keyword}

\end{frontmatter}


\section{Introduction}

Cracking is a phenomenon that occurs everywhere  in our lives, but if it is allowed to continue, it can cause fatal damage. A crack in a material occurs when the material experiences a continuous overload. However, several other factors, such as thermal expansion and contraction due to temperature changes \cite{Mackin2002,Nara2011,Vivekanandan2020}, fluid pressure (e.g., in hydraulic fracturing) \cite{Kou2019}, the diffusion of hydrogen (or hydrogen embrittlement) \cite{Dwivedi2018,Louthan1972}, chemical reactions \cite{Freiman1984}, and humidity \cite{Nara2011}, cause cracks in materials. In particular, among these phenomena, cracks due to thermal expansion are interesting to study from the viewpoint of the energy balance between elastic, thermal, and surface energies.

M. A. Biot proposed a theoretical framework for coupled thermoelasticity based on the principle of minimum entropy production \cite{Biot1956}. Biot's model is now widely known as the traditional coupled thermoelasticity model, and it has been extended to dynamical theory \cite{LordShulman1967} and to various other situations \cite{Entezari2018,GreenLindsay1972,GreenNaghdi1991,Kouchakzadeh2015,Zheng2015,Zhou2011}. As shown in Section~2.2, it satisfies an energy balance equality between the elastic and thermal energies.

In fracture mechanics, especially in the modeling and simulation of crack propagation,
a phase field approach has been recently recognized as a powerful tool. The phase field model (PFM) for fractures was first proposed by Bourdin et al. \cite{Bourdin2000} and Karma et al. \cite{Karma2001}. Then, based on the framework of variational fracture theory \cite{Bourdin2008,Francfort1998}, the techniques and applications of PFM have been extensively developed, for example \cite{Alfat2018,Amor2009,Bourdin2007,Miehe2015,Kimura2009}. We refer to \cite{Kimura2021} for further information on the development of PFM for fracture mechanics. PFM for fracture mechanics is derived as a gradient flow of the total energy, which consists of the elastic energy and the surface energy and is known to be consistent with the classical Griffith theory \cite{Bourdin2000,Kimura2021}. It allows us to handle the complex geometry of multiple, kinked, or branching cracks in both 2D and 3D without a crack path search. Comparisons with the experimental results are investigated in \cite{Nguyen2016}.

\begin{figure}[!h]
	\begin{tabular}{c}
		{\includegraphics[trim=4.5cm 12cm 4.5cm 12cm,width=0.8\textwidth]{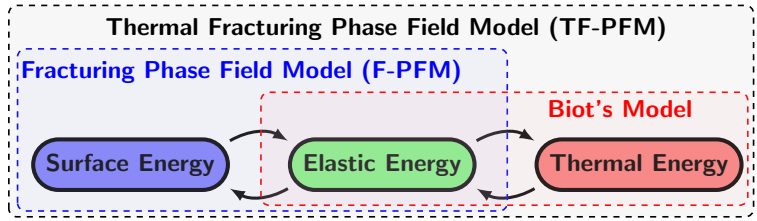}}
	\end{tabular}
	\caption{A conceptual diagram of energy balance for Biot's model, F-PFM, and TF-PFM} \label{fig:untitled-1}
\end{figure}

In this study, we deal with the modeling of thermal fracturing in an isotropic and homogeneous body by coupling the Biot thermoelasticity model and PFM. Naturally, three kinds of energy, i.e., elastic, thermal, and surface energies, appear in our stage, and the exchange and dissipation of those energies are the main interests of our research. An illustration is shown in Figure \ref{fig:untitled-1}. There are several previous works that address thermal fracturing using PFM \cite{Weilong2019,Bourdin2014,Duflot2008,Kimura2021,Miehe2015,Nguyen2017}, but they neglect the strain's influence on the heat transfer. To the best of our knowledge, a peridynamics model that employs the coupled thermoelastic equation was proposed by Gao and Oterkus \cite{Gao2019}.

The organization of this paper is as follows. In Section~\ref{sec:2}, we introduce the linear thermoelasticity model by M.A. Biot and derive its variational principle and energy dissipation property. In addition, we numerically investigate the effect of the thermal coupling term on the elastic and thermoelastic energies in an expanding region.

Section~\ref{sec:3} is devoted to PFMs for crack propagation under thermal stress. In Section~\ref{subsec:3.1}, we give a brief review of the irreversible fracturing phase field model (F-PFM) and its energy equality, which guarantees the energy dissipation property (Theorem~\ref{Theorem2:a}) and follows the works \cite{Kimura2021,Kimura2009}. In Sections~\ref{subsec:3.2} and \ref{subsec:3.3}, we propose two types of thermal fracturing phase field models (TF-PFMs). The first model, TF-PFM1, is a straightforward coupling of F-PFM and the Biot thermoelasticity model. Based on the variational principle of the Biot model (Proposition~\ref{Proposition:2}), we show a partial energy equality for a fixed temperature (Theorem~\ref{Theorem2:b}). However, it does not satisfy the energy equality for the total energy, which consists of the elastic, thermal, and surface energies.

The second model, TF-PFM2, presented in Section~\ref{subsec:3.3} is another natural coupling of F-PFM and the Biot thermoelasticity model based on the energy equality of the Biot model (Theorem~\ref{Theorem1}). We prove an energy equality for TF-PFM2 in Theorem \ref{Theorem3}. Since we consider several models (Biot's model, F-PFM, and TF-PFMs) and their energy qualities, for the readers' convenience, we list the energies and energy equalities for each model in Tables \ref{tab:1} and \ref{tab:2}.

In Section \ref{sec-4}, we show some numerical comparisons between two TF-PFMs using nondimensionalized equations. We investigate the effects of the thermal coupling in TF-PFM1 and TF-PFM2 on the crack speed and the crack path by changing a dimensionless coupling parameter $\delta$. The last section shows some conclusions and comments on further topics. 

\linespread{1.}{
\begin{table}[!h]
	\caption{List of energies}
	\label{tab:1}
	\begin{tabular}{p{2.5cm}m{7.2cm}p{1.2cm}}
		\hline\noalign{\smallskip}
		Type of Energy & Definition & Equation  \\
		\noalign{\smallskip}\hline\noalign{\smallskip}
		Elastic & $E_{el}(u) := \dfrac{1}{2}\displaystyle\int_{\Omega} \sigma[u]:e[u]~dx$ & \eqref{original-elastic}\\[17pt]
		Thermoelastic  &  $E_{el}^{*}(u,\Theta) := \dfrac{1}{2}\displaystyle\int_{\Omega} \sigma^{*}[u,\Theta]:e^{*}[u,\Theta]~dx$ & \eqref{Original-Thermo-Elastic} \\[17pt]
		Thermal & $E_{th}(\Theta) := \dfrac{\chi}{2\Theta_{0}}\displaystyle\int_{\Omega} \big| \Theta(x) - \Theta_{0}\big|^{2}~dx$ & \eqref{ThermalEn}\\[17pt]
		Modified elastic & $\mathcal{E}_{el}(u,z) := \dfrac{1}{2}\displaystyle\int_{\Omega} (1-z)^{2}\sigma[u]:e[u]~dx$ & \eqref{Modify-Elastic}\\[17pt]
		Modified thermoelastic & $\mathcal{E}_{el}^{*}(u,\Theta,z) := \dfrac{1}{2}\displaystyle\int_{\Omega} (1-z)^{2}\sigma^{*}[u,\Theta]:e^{*}[u,\Theta]~dx$ & \eqref{Thermo-Elastic}\\[17pt]
		Surface & $E_{s}(z) := \dfrac{1}{2}\displaystyle\int_{\Omega} \gamma_{*}\Big( \epsilon |\nabla z|^{2} + \frac{|z|^{2}}{\epsilon}\Big)~dx$ &  \eqref{Surface-Energy}\\
		\noalign{\smallskip}\hline
	\end{tabular}
\end{table}}

\begin{table}[!h]
	\begin{threeparttable}[b]
		\caption{Different forms of energy equalities}
		\label{tab:2}
		\begin{tabular}{p{2.6cm}p{2.5cm}m{3.8cm}m{1.5cm}}
			\hline\noalign{\smallskip}
			Model & Strong Form & Energy & Energy Equality \\
			\noalign{\smallskip}\hline\noalign{\smallskip}
			Linear elasticity & \eqref{Eq:LinElastic} & $E_{el}(u)$ & -\\[6pt]
			Biot's model & \eqref{Biot-Model1} - \eqref{Biot-Model2} & $E_{el}(u) + E_{th}(\Theta)$ & \eqref{EnergyEquality:1}\\[6pt]
			F-PFM & \eqref{PML1} - \eqref{PML2}  &  $\mathcal{E}_{el}(u,z) + E_{s}(z)$ & \eqref{Eq:Theorem2:a}\\[6pt]
			TF-PFM1 & \eqref{TMF1a} - \eqref{TMF1c} & $\mathcal{E}_{el}^{*}(u,\Theta,z) + E_{s}(z) $ & \eqref{Eq:Theorem2:b}\tnote{a} \\[6pt]
			TF-PFM2 & \eqref{TMF2a} - \eqref{TMF2c} & $\mathcal{E}_{el}(u,z) + E_{s}(z)  + E_{th}(\Theta)$ & \eqref{Eq:Theorem3:c}\\
			\noalign{\smallskip}\hline
		\end{tabular}
		\begin{tablenotes}\footnotesize
			\item [a] When a temperature $\Theta = \Theta(x) \in L^{2} (\Omega)$ is given. 
		\end{tablenotes}
	\end{threeparttable}
\end{table}

\linespread{1.}{
\begin{table}[!h]
	\begin{threeparttable}[b]
	\caption{List of  physical properties}
	\label{tab:3}
	\begin{tabular}{p{1.2cm}|p{3.6cm} p{.05cm} p{1.2cm}|p{4.0cm}}
		\cline{1-2} \cline{4-5}\noalign{\smallskip}
		Symbol & Physical meaning [unit] & & Symbol & Physical meaning [unit]  \\[4pt]		\noalign{\smallskip}\cline{1-2} \cline{4-5}\noalign{\smallskip}
		$u$ & Displacement [\si{\metre}] & & $\sigma^{*}[u,\Theta]$ & Stress tensor with thermal effect [\si{\pascal}] \\[18pt]
		$\Theta$  &  Temperature [\si{\kelvin}] &  & $e^{*}[u,\Theta]$ & Strain tensor with thermal effect [-]\\[18pt]
		$\Theta_{0}$ & Reference temperature [\si{\kelvin}] & &  $\beta$ & Stress thermal modulus [\si[]{\pascal\cdot\kelvin^{-1}}]\\[18pt]
		$z$ & Damage variable [-] & &  $\kappa_{0}$ & Thermal conductivity [\si[]{\watt\cdot\meter^{-1}\cdot\kelvin^{-1}}]\\[18pt]
		$\sigma[u]$ & Stress tensor [\si{\pascal}] & &  $\chi$& Volumetric heat capacity [\si[]{\joule\cdot\kelvin^{-1}\cdot\metre^{-3}}]\\[18pt]
		$e[u]$ & Strain tensor [-] &  & $a_{L}$ & Coefficient of linear thermal expansion [\si{\per\kelvin}]\\[18pt]
		$E_{\text{Y}}$ & Young's modulus [\si{\pascal}] &  & ${\delta}$& Thermoelastic coupling parameter [-]\\[18pt]
		$\nu_{\text{P}}$& Poisson ratio [-] &  &	$\gamma_{*}$ & Critical energy release rate\tnote{a}~ [\si{\pascal\cdot\metre}] \\[18pt]
		$\lambda$, $\mu$ & Lam\'{e}'s constants\tnote{b}~ [\si{\pascal}] &	&  $\epsilon$ & Length scale in F-PFM or TF-PFM [\si{\metre}]\\[18pt]
		$t$ & Time [\si{\second}] & &  $\alpha$ & Time regularization parameter in F-PFM or TF-PFM~[\si{\pascal\cdot\second}] \\ \cline{1-2} \cline{4-5} 
	\end{tabular}
	\begin{tablenotes}\footnotesize
	\item [a] $\gamma_{*}$ is usually denoted by $\text{G}_\text{c}$ \cite{Anderson2017,Kimura2021}.\\
	\item [b] $\lambda$ and $\mu$ are written as $\lambda = \cfrac{E_{\text{Y}}\nu_{\text{P}}}{(1+\nu_{\text{P}})(1-2\nu_{\text{P}})}$ and $\mu = \cfrac{E_{\text{Y}}}{2(1-\nu_{\text{P}})}$.
\end{tablenotes}
\end{threeparttable}
\end{table}
}

To easily understand the relevant notation and symbols in this paper, we introduce them in this section. Let $\Omega$ be a bounded domain in $\mathbb{R}^{d}$ ($d = 2$ or $3$). The position in $\mathbb{R}^{d}$ is denoted by $x = (x_{1}.\cdots,x_{d})^{T} \in \mathbb{R}^{d}$, where $~^{T}$ denotes the transposition of a vector or matrix. Let $\nabla$, $\mbox{div}$, and $\Delta$ be the gradient, divergence, and Laplacian operators with respect to $x$, respectively. For simplicity, we write $\dot{u}$, $\dot{\Theta}$, and $\dot{z}$ as the partial derivatives of $u$, $\Theta$ and $z$ with respect to $t$, respectively. For simplicity, we often denote $u(t) := u(\cdot,t)$, etc. The space of the real-valued (symmetric) $ d \times d$ matrix is denoted by $\mathbb{R}^{d\times d}$ ($\mathbb{R}^{d\times d}_{sym}$). The inner product of square matrices $A, B \in \mathbb{R}^{d\times d}$ is denoted by $A:B := \sum_{i,j= 1}^{d}A_{ij}B_{ij}$. Using $L^{2}(\Omega)$, we refer to the Lebesgue space on $\Omega$, while $H^{1}(\Omega,\mathbb{R}^{d})$ and $H^{\frac{1}{2}}(\Gamma_{D}^{u},\mathbb{R}^{d})$ represent the Sobolev space on $\Omega$ and its trace space on the boundary $\Gamma_{D}^{u}$, respectively. For more details on Sobolev spaces, we refer to the review in \cite{Girault1979}. In addition, we summarize the physical properties used in this paper in Table \ref{tab:3}.

\section{Thermoelasticity Model}\label{sec:2}
\subsection{Formulation of the problem}\label{subsec:2.1}
M.A. Biot \cite{Biot1956} proposed the following mathematical model for coupled thermoelasticity:
\begin{subequations}
	\begin{empheq}[left=\empheqlbrace]{align}
	&  -\mbox{div} \sigma{[u]}= \beta\nabla\Theta &\mbox{in}~ \Omega \times [0,T], \label{Biot-Model1}\\
	& \chi\dfrac{\partial}{\partial t}\Theta = \kappa_{0} \Delta\Theta- \Theta_{0}\beta\dfrac{\partial}{\partial t}(\mbox{div}u)  & \mbox{in} ~ \Omega \times (0,T], \label{Biot-Model2}
	\end{empheq}\label{BiotModel}
\end{subequations}\\[-3pt]
where $\Omega$ is a bounded domain in $\mathbb{R}^{d}$ ($d = 2$ or $3$). We suppose that $\Omega$ is an isotropic elastic body and consider the thermoelastic coupling between the mechanical deformation and the thermal expansion in $\Omega$. The constant $\beta$ is defined by $\beta := a_{L}(d\lambda + 2\mu)$ with $a_{L} > 0$ as the coefficient of linear thermal expansion and $\mu(> 0)$; $\lambda(> -\frac{2\mu}{d})$ are Lam\'{e}'s constants.

The unknown functions in \eqref{Biot-Model1} and \eqref{Biot-Model2} are the displacement $u(x,t) = (u_{1}(x,t)$ $,\cdots,u_{d}(x,t))^{T} \in \mathbb{R}^{d}$ and the temperature $\Theta(x,t) \in \mathbb{R}$. In addition, the constant $\Theta_{0} > 0$ is a fixed reference temperature. Similarly, strain $e[u]$ and stress tensors $\sigma[u]$ are defined as
\begin{subequations}
	\begin{empheq}[]{align}
	&e[u] := \frac{1}{2}\left(\nabla u^{T} + (\nabla u^{T})^{T}\right) \in \mathbb{R}_{sym}^{d\times d},\\
	& \sigma[u] := Ce[u] = \lambda(\mbox{div} u) {I} + 2\mu e[u] \in \mathbb{R}_{sym}^{d\times d},\label{sigma}
	\end{empheq}
\end{subequations}
where $C := (c_{ijkl}),~ c_{ijkl} = \lambda \delta_{ij}\delta_{kl} + \mu(\delta_{ik}\delta_{jl} + \delta_{il}\delta_{jk})$ is an isotropic elastic tensor and ${I}$ is the identity matrix of size $d$. From \eqref{sigma}, (\ref{Biot-Model1}) is also written in the form
\begin{align}
& -\mu\Delta u - (\lambda + \mu)\nabla(\mbox{div}u) = \beta\nabla\Theta.\notag
\end{align}
The term $\beta\nabla\Theta$ in \eqref{Biot-Model1} and the term $\Theta_{0}\beta\frac{\partial}{\partial t}(\mbox{div} u)$ in \eqref{Biot-Model2} represent the body force due to thermal expansion and the heat source due to the volume change rate, respectively. We remark that when $a_{L} = 0$, \eqref{Biot-Model1} and \eqref{Biot-Model2} are decoupled.

It is convenient to introduce the following strain and stress tensors, including the thermal effect.  
\begin{subequations}
	\begin{empheq}[]{align}
	& e^{*}[u,\Theta] := e[u] - a_{L}(\Theta(x,t) - \Theta_{0})I \in \mathbb{R}^{d\times d}_{sym}, \notag\\
	& \sigma^{*}[u,\Theta] := Ce^{*}[u,\Theta] = \sigma[u] - \beta(\Theta(x,t) - \Theta_{0})I \in \mathbb{R}^{d\times d}_{sym}. \notag
	\end{empheq}
\end{subequations}
Using the thermal stress tensor $\sigma^{*}[u,\Theta]$, (\ref{Biot-Model1}) can be written in the following form:
\begin{align}
& -\mbox{div} \sigma^{*}[u,\Theta] = 0.\notag
\end{align}
This means that the force $\sigma^{*}[u,\Theta]$ is in equilibrium in $\Omega$. In the preceding equation, \eqref{BiotModel} represents the force balance and the thermal diffusion in $\Omega$, respectively.

The system in  \eqref{BiotModel} is complemented by the following boundary and initial conditions:
\begin{subequations}
	\begin{empheq}[left=\empheqlbrace]{align}
	& u = u_{D}(x,t) & \mbox{on}~ \Gamma_{D}^{u} \times [0,T],\label{BC01}\\
	& \sigma^{*}[u,\Theta]n = 0 & \mbox{on}~ \Gamma_{N}^{u} \times [0,T], \label{BC02}\\
	& \Theta = \Theta_{D}(x,t) & \mbox{on}~ {\Gamma}_{D}^{\Theta} \times [0,T], \label{BC03}\\
	& \frac{\partial\Theta}{\partial n} = 0 & \mbox{on}~ {\Gamma}_{N}^{\Theta}  \times [0,T], \label{BC04}\\
	& \Theta(x,0) = \Theta_{*}(x) & \mbox{in} ~ \Omega,~\hspace{-2pt}\qquad\qquad \label{IC00}
	\end{empheq}\label{Boundary-Cond}
\end{subequations}\\[-5pt]
where $n$ is the outward unit normal vector along the boundary, $\Gamma = \Gamma_{D}^{u} ~\cup~ \Gamma_{N}^{u}~  (\Gamma = \Gamma_{D}^{\Theta}~ \cup~ \Gamma_{N}^{\Theta})$ with $\Gamma_{D}^{u} ~\cap ~\Gamma_{N}^{u} ~= \emptyset~({\Gamma}_{D}^{\Theta} ~\cap ~ {\Gamma}_{N}^{\Theta} =\emptyset)$. The boundaries  $\Gamma_{D}^{u}$ and $\Gamma_{N}^{u}$ ($\Gamma_{D}^{\Theta}$ and $\Gamma_{N}^{\Theta}$) are the Dirichlet and Neumann boundaries for $u$ (for $\Theta$), respectively. We  suppose that the $(d-1)$-dimensional volume of $\Gamma_{D}^{u}$ is positive  for the solvability of $u$. 
\begin{Remark}\label{Remark:1}
\normalfont Instead of boundary conditions \eqref{BC01} and \eqref{BC02}, we can also consider the following mixed-type condition. When $d=2$, on a part of the boundary  (which we denote by $\Gamma_{DN}^{u}$), $u=(u_{1},u_{2})^{T}$ and
\begin{align}
\left\{
\begin{array}{ll}
u_{1} = u_{D1} & \qquad \mbox{on}~\Gamma_{DN}^{u}, \\
(\sigma^{*}[u,\Theta]n)\cdot e_{2}= 0 & \qquad\mbox{on}~ \Gamma_{DN}^{u},
\end{array}
\right. \notag
\end{align}
or
\begin{align}
\left\{
\begin{array}{ll}
u_{2} = u_{D2} & \qquad \mbox{on}~\Gamma_{DN}^{u}, \\
(\sigma^{*}[u,\Theta]n)\cdot e_{1}= 0 & \qquad\mbox{on}~ \Gamma_{DN}^{u},
\end{array}
\right. \notag
\end{align}
where $u_{Di} := \Gamma_{DN}^{u} \mapsto \mathbb{R}$ is a given horizontal or vertical displacement and $e_{1}=(1,0)^{T}$, $e_{2} = (0,1)^{T}$. These types of mixed boundary conditions are considered in Section \ref{subsubsec:2.4.1}  and Section \ref{SubSec4.4.1}. Even for these mixed-type boundary conditions, we can easily extend the following arguments on weak solutions, variational principles, and energy equalities. 
\end{Remark}

\subsection{Variational principle and energy equality}\label{subsec:2.2}
This section aims to show a variational principle and provide an energy equality that implies the energy dissipation property for the system \eqref{BiotModel}. In linear elasticity theory, a weak form of the boundary value problem for $u_{D}\in H^{\frac{1}{2}}(\Gamma_{D}^{u}; \mathbb{R}^{d})$ is
\begin{align}
\left\{
\begin{array}{ll}
-\mbox{div} \sigma[u] = 0 & \qquad\mbox{in}~ \Omega, \\
u = u_{D} & \qquad\mbox{on} ~\Gamma_{D}^{u}, \\
\sigma[u]n = 0 & \qquad\mbox{on} ~\Gamma_{N}^{u},
\end{array}
\right.\label{Eq:LinElastic}
\end{align}
which is given by
\begin{align}
& u\in V^{u}(u_{D}), ~ \int_{\Omega} \sigma[u]:e[v]~dx = 0 \quad \text{for all}~ v\in V^{u}(0), \notag
\end{align}
where 
\begin{align}
& V^{u}(u_{D}) := \left\{u\in H^{1}(\Omega;\mathbb{R}^{d}); ~u\big|_{\Gamma_{D}^{u}} = u_{D}\right\}. \label{Original-Thermo-Elastic1}
\end{align}
A weak solution uniquely exists and is given by
\begin{equation}
u = \operatorname*{argmin}_{v \in V^{u}(u_{D})} {E_{el}(v)},\notag
\end{equation}
where 
\begin{align}
& E_{el}(v) := \frac{1}{2}\int_{\Omega} \sigma[v]:e[v]~dx \quad (v\in H^{1}\left(\Omega;\mathbb{R}^{d})\right) \label{original-elastic}
\end{align}
is an elastic energy. This is known as a variational principle \cite{Ciarlet2002,Duvaut1976}. For a fixed $\Theta(x)$, a weak form for $u$ of \eqref{Biot-Model1} and its variational principle are derived as follows.
\begin{Proposition} \label{Proposition:1}
	For $u \in H^{2}(\Omega ; \mathbb{R}^{d})$ and $\Theta \in H^{1}(\Omega)$,
	\begin{align}
	\left\{
	\begin{array}{ll}
	-\normalfont{\mbox{div}} \sigma^{*}[u,\Theta] = 0 & \qquad\mbox{in}~ \Omega, \\
	u = u_{D} & \qquad\mbox{on} ~\Gamma_{D}^{u}, \\
	\sigma^{*}[u,\Theta]n = 0 & \qquad\mbox{on} ~\Gamma_{N}^{u},
	\end{array}
	\right.\notag
	\end{align}
	is equivalent to the following weak form:
	\begin{align}
	\left\{
	\begin{array}{ll}
	\displaystyle\int_{\Omega}\sigma^{*}[u,\Theta]:e[v]~dx = 0 & \quad\mbox{for all} ~ v \in  V^{u}(0), \\[8pt]
	u\in V^{u}(u_{D}). &
	\end{array}
	\right. \label{Eq-Prop:1b}
	\end{align} 
	\normalfont\begin{Proof}\normalfont
		For $v\in V^{u}(0) $, we have
		\begin{eqnarray}
		\int_{\Omega} \left( -\mbox{div}\sigma^{*}[u,\Theta]\right)\cdot v~dx &=& \int_{\Omega} \sigma^{*}[u,\Theta]:e[v]~dx - \int_{\Gamma_{N}^{u}} (\sigma^{*}[u,\Theta]n)\cdot v~ds. \notag
		\end{eqnarray}
		The equivalency immediately follows from this equation.  \qed
	\end{Proof}
\end{Proposition}
\begin{Proposition}[Variational principle]\label{Proposition:2}
	  For a given $\Theta \in L^{2}(\Omega)$, $u_{D} \in H^{\frac{1}{2}} (\Gamma_{D}^{u} ; \mathbb{R}^{d})$, there exists a unique weak solution $u \in H^{1}(\Omega; \mathbb{R}^{d})$ that satisfies \eqref{Eq-Prop:1b}. Furthermore, the solution $u$ is a unique minimizer of the variational problem:
	\begin{equation}
	u = \operatorname*{argmin}_{v \in V^{u}(u_{D})} {E_{el}^{*}(v,\Theta)},\notag
	\end{equation}
	where
	\begin{align}
	& E_{el}^{*}(v,\Theta) = \frac{1}{2}\int_{\Omega} \sigma^{*}[v,\Theta]:e^{*}[v,\Theta]~dx. \label{Original-Thermo-Elastic}
	\end{align}
	We remark that $E_{el}^{*}(v,\Theta)$ represents thermoelastic energy. 
	\normalfont\begin{Proof}\label{Proof:1}
		\normalfont The unique existence of a weak solution for $u$ is shown by the Lax-Milgram theorem \cite{Ciarlet2002}  since \eqref{Eq-Prop:1b} is written as
		\begin{align}
		\left\{
		\begin{array}{ll}
		\displaystyle\int_{\Omega}\sigma[u]:e[v]~dx = \displaystyle\int_{\Omega} \beta(\Theta - \Theta_{0})\mbox{div}v~dx, \\[3pt]
		u \in V^{u}(u_{D}) \qquad(\mbox{for all}~ v \in V^{u}(0)).
		\end{array}
		\right.\notag
		\end{align}
		The coercivity of the above weak form is known as Korn's second inequality \cite{Ciarlet2002}:
		\begin{align}
		&^{\exists} a_{0} >0 ~\text{such that} ~ \int_{\Omega} \sigma[v]:e[v]~dx \geq a_{0}\| v\|^{2}_{H^{1}(\Omega; \mathbb{R}^{2})}, \quad \mbox{for all}~ v \in V^{u}(0). \notag
		\end{align} 
		For a weak solution $u$ and any $v \in V^{u}(0)$, using the equalities
		\begin{align}
		& \sigma^{*}[u+v,\Theta] = \sigma^{*}[u,\Theta] + \sigma[v],\notag\\
		& e^{*}[u+v,\Theta] = e^{*}[u,\Theta] + e[v],\notag\\		
		& \sigma^{*}[u,\Theta]:e[v] = e^{*}[u,\Theta]:\sigma[v] ,\notag
		\end{align}
		we have
		\begin{align}
		& E_{el}^{*}(u+v,\Theta) - E_{el}^{*}(u,\Theta)  \nonumber\\
		& \qquad \quad = \frac{1}{2}\int_{\Omega} \sigma^{*}[u+v,\Theta]:e^{*}[u+v,\Theta]~dx - \frac{1}{2}\int_{\Omega} \sigma^{*}[u,\Theta]:e^{*}[u,\Theta]~dx\qquad\nonumber\\
		& \qquad \quad = \int_{\Omega} \sigma^{*}[u,\Theta]:e[v]~dx + \frac{1}{2}\int_{\Omega} \sigma[v]:e[v]~dx\nonumber\\
		& \qquad \quad = \frac{1}{2}\int_{\Omega} \sigma[v]:e[v]~dx \geq 0.\notag
		\end{align}
		This shows that $u$ is a minimizer of $E^{*}_{el}(u,\Theta)$ among $V^{u}(u_{D})$.  
		
		On the other hand, if $u$ is a minimizer, the first variation of $E_{el}^{*}$ vanishes at $u$; i.e., for all $v\in V^{u}(0)$, we have
		\begin{align}
		& 0 = \frac{d}{ds} E^{*}_{el}(u+sv,\Theta)\big|_{s=0} = \int_{\Omega} \sigma^{*}[u,\Theta]:e[v]~dx.\notag
		\end{align}
		Hence, $u$ is a weak solution. Summarizing the above, there exists a unique weak solution to \eqref{original-elastic}, and $u$ is a weak solution if and only if it is a minimizer of $E^{*}_{el}$ among $V^{u}(u_{D})$. \qed
	\end{Proof}
\end{Proposition}
The next theorem represents a dissipation of the sum of the elastic and thermal energies during the thermomechanical process. We define thermal energy as
\begin{align}
	E_{th}(\Theta) := \frac{\chi}{2\Theta_{0}} \int_{\Omega} \left|\Theta(x) - \Theta_{0}\right|^{2}~dx. \label{ThermalEn}
\end{align}

\begin{Theorem}[Energy equality for Biot's model]\label{Theorem1}
Let $(u(x,t),\Theta(x,t))$ be a sufficiently smooth solution to \eqref{BiotModel}  and \eqref{Boundary-Cond}. In addition, we suppose that $u_{D}$ does not depend on $t$ and $\Theta_{D}=\Theta_{0}$. Then
\begin{align}
\frac{d}{dt}\Big(E_{el}(u(t)) + E_{th}(\Theta(t))\Big) &= -\frac{\kappa_{0}}{\Theta_{0}}\int_{\Omega} \left|\nabla\Theta(t)\right|^{2}~dx \leq 0. \label{EnergyEquality:1}
\end{align}

\normalfont\begin{Proof} \label{Proof:2}
\normalfont 
Since
		\begin{align}
		\frac{d}{dt}\left(\frac{1}{2}\sigma[u]:e[u]\right) &= \sigma[u]:e[\dot{u}]\notag\\
		& = (\sigma^{*}[u,\Theta]+\beta(\Theta-\Theta_{0})I):e[\dot{u}]\notag\\
		& = \sigma^{*}[u,\Theta] - \beta (\Theta-\Theta_{0})\mbox{div}\dot{u} \label{eq-proof1}
		\end{align}
we obtain
		\begin{align}
		\frac{d}{dt}E_{el}(u(t))  &= \frac{1}{2} \int_{\Omega} \frac{d}{dt} \left(\sigma[u]:e[u]\right)~dx\notag\\
		&  = \int_{\Omega} \sigma^{*}[u,\Theta]:e[\dot{u}]~dx + \int_{\Omega} \beta (\Theta - \Theta_{0})(\mbox{div}\dot{u})~dx \notag\\
		& =\int_{\Omega} \beta(\Theta - \Theta_{0})(\mbox{div}\dot{u})~dx. \label{Eq:Elastic01}
		\end{align}
		Substituting \eqref{Biot-Model2} into \eqref{Eq:Elastic01} and using the boundary conditions \eqref{BC03} and \eqref{BC04} for $\Theta$, we obtain
		\begin{align}
		\frac{d}{d t}E_{el}(u(t)) &= \int_{\Omega} \Big(\frac{1}{\Theta_{0}}(\Theta - \Theta_{0})\Big\{\kappa_{0}\Delta\Theta - \chi\frac{\partial \Theta}{\partial t}\Big\} \Big)dx\nonumber\\
		&=\frac{\kappa_{0}}{\Theta_{0}}\int_{\Gamma} (\Theta - \Theta_{0})\frac{\partial\Theta}{\partial n}ds - \frac{\kappa_{0}}{\Theta_{0}}\int_{\Omega} \big|\nabla\Theta\big|^{2} dx - \frac{d}{dt}\Big(\frac{\chi}{2\Theta_{0}}\int_{\Omega}\big|\Theta - \Theta_{0}\big|^{2}dx\Big)\notag\\
		&= -\frac{\kappa_{0}}{\Theta_{0}} \int_{\Omega} \big|\nabla\Theta\big|^{2}dx - \frac{d}{dt}E_{th}(\Theta(t)). \notag
		\end{align}
This gives the energy equality for \eqref{BC-IC}. \qed
\end{Proof}
\end{Theorem}

As shown in Proposition \ref{Proposition:2}  and Theorem \ref{Theorem1}, Biot's thermoelasticity model is related to both energies $E_{el}(u)$ and $E_{el}^{*}(u,\Theta)$. We denote their energy densities as follows:
\begin{align}
&{W}(u) := \sigma[u]:e[u],\\
&{W}^{*}(u,\Theta):= \sigma^{*}[u,\Theta]:e^{*}[u,\Theta],
\end{align}
where ${W}(u)$ and ${W}^{*}(u,\Theta)$ are the elastic and thermoelastic energy densities, respectively.

\subsection{Numerical Experiment}\label{subsec:2.3}
\subsubsection{Nondimensional setting}
In the following numerical examples, we introduce a nondimensional form of Biot's model. We consider the following scaling for $x$, $t$, $u$, $C$ (or $\lambda$, $\mu$), and $\Theta$:
\begin{eqnarray}
\tilde{x} = \frac{x}{c_{x}}, ~ \tilde{t} = \frac{t}{c_{t}}, ~ \tilde{u} = \frac{u}{c_{u}}, ~ \tilde{C} = \frac{C}{c_{e}}, ~ \tilde{\Theta} = \frac{\Theta - \Theta_{0}}{c_{\Theta}},~ \tilde{a}_{L} = \frac{c_{x}c_{\Theta}}{c_{u}}a_{L},~ \tilde{\beta} = 1, \label{ND}
\end{eqnarray}
where $c_{x}$,  $c_{t}$, $c_{u}$, $c_{e}$, and $c_{\Theta}> 0$ are the scaling parameters. Let $c_{x}$  [\si{\meter}], $c_{e}$ [\si{\pascal}], and $c_{\Theta}$ [\si{\kelvin}] be characteristic scales for the length of the domain, the size of the elastic tensor and the temperature, respectively. The parameters $c_{t}$ and $c_{u}$ are defined as 
\begin{eqnarray}
c_{t} := \frac{c_{x}^{2}{\chi}}{\kappa_{0}} ~[\si{\second}], \qquad c_{u} := \frac{c_{\Theta}c_{x}\beta}{c_{e}} ~[\si{\meter}], \label{ND2}
\end{eqnarray}
where $\chi$ [\si{\pascal\cdot\kelvin^{-1}}], $\kappa_{0}$ [\si{\pascal\cdot\meter^{2}\cdot\second^{-1}\cdot\kelvin^{-1}}] and $\beta = a_{L}(d\lambda + 2\mu)$ [\si{\pascal\cdot\kelvin^{-1}}]. Then, \eqref{BiotModel} is written in the following nondimensional form:
\begin{subequations}
	\begin{empheq}[left=\empheqlbrace]{align}
	&  -\widetilde{\mbox{div}} \tilde{\sigma}{[\tilde{u}]}= \tilde{\nabla}\tilde{\Theta} &\mbox{in}~ \tilde{\Omega} \times [0,\tilde{T}], \label{ND-Biot-Model1}\\
	&\dfrac{\partial}{\partial \tilde{t}}\tilde{\Theta} = \tilde{\Delta}\tilde{\Theta}- \delta\dfrac{\partial}{\partial \tilde{t}}(\widetilde{\mbox{div}}\tilde{u})  & \mbox{in} ~ \tilde{\Omega} \times (0,\tilde{T}]. \label{ND-Biot-Model2}
	\end{empheq}\label{ND-BiotModel}
\end{subequations}\\[-3pt]
The system \eqref{ND-BiotModel} has only three parameters, $\tilde{\lambda}$, $\tilde{\mu}$, and $\delta$. The parameter $\delta$ is a nondimensional thermoelastic coupling parameter defined by 
\begin{eqnarray}
{\delta} &:=& \frac{\Theta_{0}\beta^{2}}{c_{e}\chi}~ [-], \notag
\end{eqnarray}
and $\delta > 0$. If we choose $\delta=0$, \eqref{ND-Biot-Model2} is decoupled from \eqref{ND-Biot-Model1}, and the temperature field $\tilde{\Theta}$ in \eqref{ND-Biot-Model1} is essentially a given function. In the following example, the case $\delta = 0$ is referred to as the uncoupled case. 

Under the above scaling, we denote the (thermo)elastic strain, stress tensors, and (thermo)elastic energy densities as follows:
\begin{subequations}
\begin{empheq}[]{align}
	& \tilde{e}[\tilde{u}] : = \frac{1}{2}\left(\frac{\partial\tilde{u}_{i}}{\partial\tilde{x}_{j}}+ \frac{\partial\tilde{u}_{j}}{\partial\tilde{x}_{i}}\right) = \frac{c_{x}}{c_{u}}e[u],\\
	&\tilde{\sigma}[\tilde{u}] := \tilde{C}\tilde{e}[\tilde{u}] = \frac{c_{x}}{c_{u}c_{e}}\sigma[u],\\
	&\widetilde{{W}}(\tilde{u}) := \tilde{\sigma}[\tilde{u}]:\tilde{e}[\tilde{u}] = \frac{c_{e}}{(\beta c_{\Theta})^{2}}{W}[u],\\
	& \tilde{\sigma}^{*}[\tilde{u},\tilde{\Theta}] := \tilde{\sigma}[\tilde{u}] - \tilde{\Theta}I = \frac{1}{\beta c_{\Theta}}\sigma^{*}[u,\Theta],\\
	& \tilde{e}^{*}[\tilde{u},\tilde{\Theta}] := \tilde{e}[\tilde{u}] - \tilde{a}_{L}\tilde{\Theta}I = \frac{c_{x}}{c_{u}}\sigma^{*}[u,\Theta],\\
	&\widetilde{{W}}^{*}(\tilde{u},\tilde{\Theta}) := \tilde{\sigma}^{*}[\tilde{u},\tilde{\Theta}]:\tilde{e}^{*}[\tilde{u},\tilde{\Theta}] = \frac{c_{e}}{(\beta c_{\Theta})^{2}}{W}^{*}[u,\Theta].
\end{empheq}\label{Non-Scaling}
\end{subequations}\\[-3pt]
In the following section, we apply these nondimensional forms and omit $\sim$ for simplicity. 

\subsubsection{Numerical setup and time discretization}
\begin{figure}[!h]
	\begin{tabular}{cc}
		{\includegraphics[trim=1cm 9cm 1cm 7.5cm,width=0.5\textwidth]{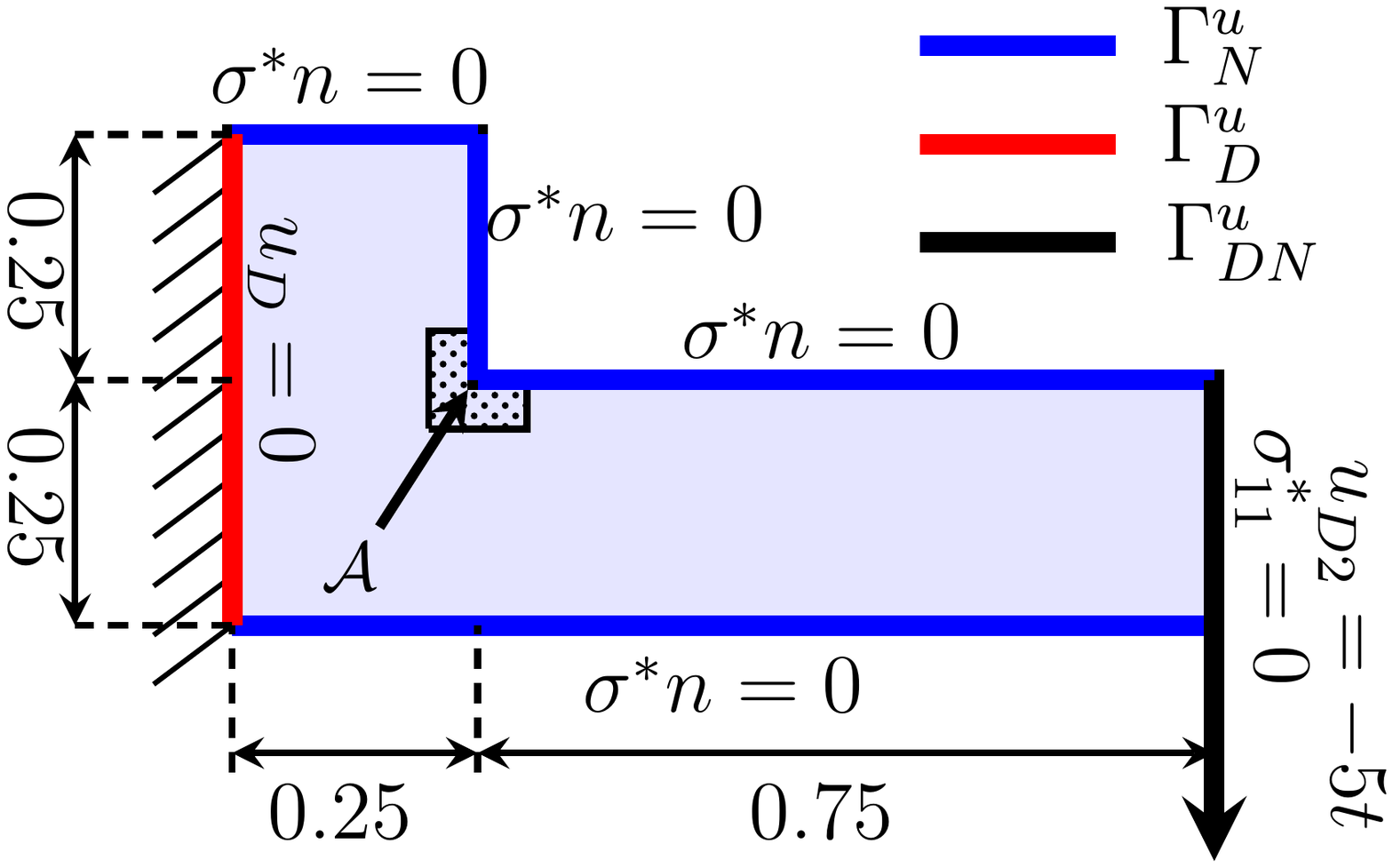}}& \qquad
		{\includegraphics[trim=1cm 7.5cm 1cm 7.5cm,width=0.35\textwidth]{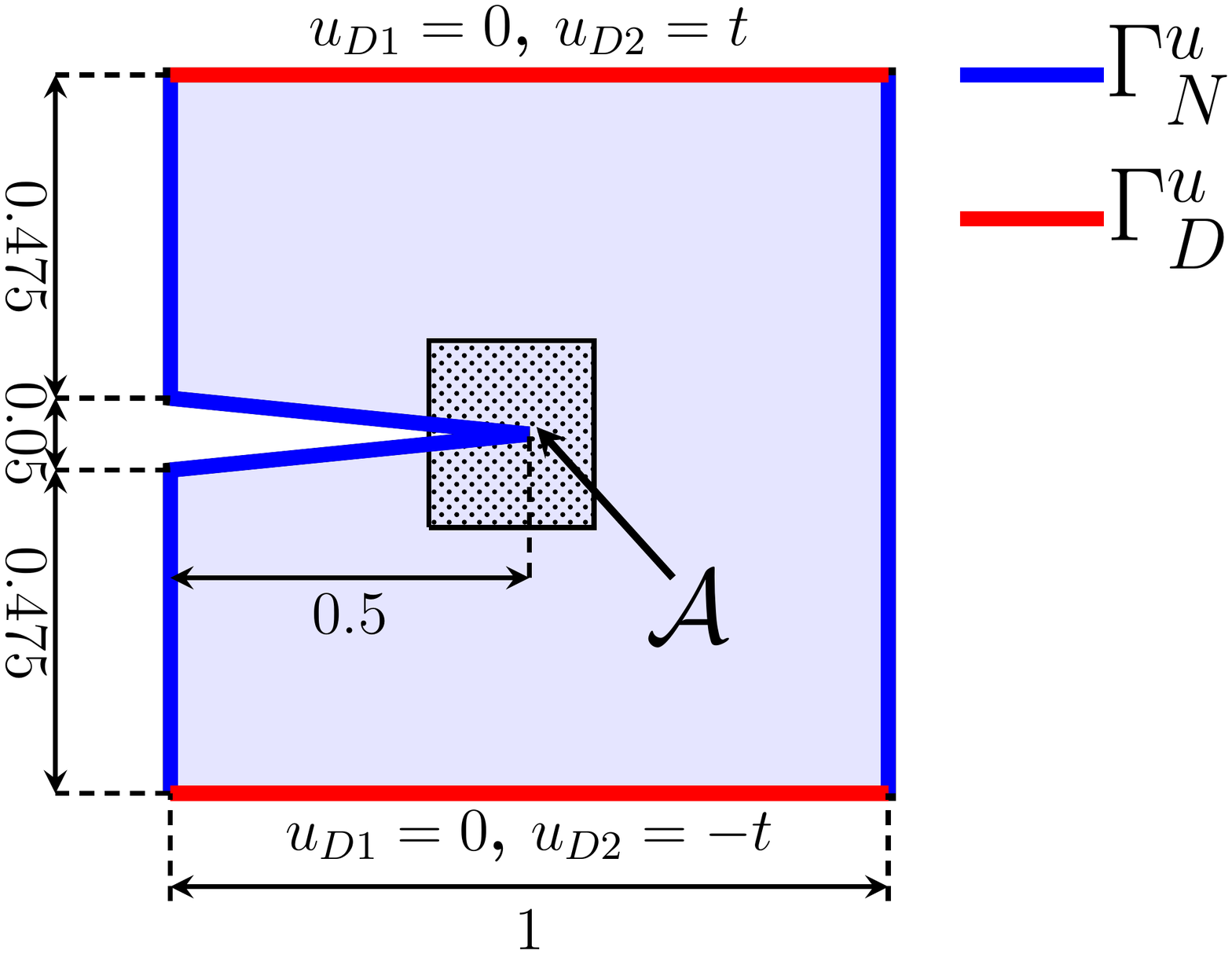}} 
	\end{tabular}
	\caption{An L-shaped cantilever (left) and a cracked domain (right) with the subdomain $\mathcal{A}$ as an observation area.} 
	\label{fig:3}
\end{figure}

In the following examples, we set Young's modulus $E_{\text{Y}} = 1$, Poisson's ratio $\nu_{\text{P}} = 0.32$, the coefficient of linear thermal expansion $a_{L} = 0.475$ and the thermoelasticity  coupling parameter ${\delta} = 0.0,~0.1,~0.5$ in the nondimensional form of \eqref{ND-BiotModel}. We consider two numerical examples for \eqref{ND-BiotModel}, an L-shaped cantilever domain and a square domain with a crack (more precisely, a very sharp notch), as illustrated in Figure \ref{fig:3}.

We apply the following implicit time discretization for \eqref{ND-BiotModel}:
\begin{align}
&\left\{
\begin{array}{ll}
-\displaystyle{{\text{div}}{\sigma}^{*}[{u}^{k},\Theta^{k-1}]= 0}& \displaystyle\mbox{in}~\Omega,\\[5pt]
\displaystyle\frac{{\Theta}^{k}- {\Theta}^{k-1}}{\Delta{t}}  - {\Delta}{\Theta}^{k} + {\delta}{\text{div}}\left(\frac{{u}^{k} - {u}^{k-1}}{\Delta{t}}\right) = 0 & \displaystyle\mbox{in}~\Omega,		    
\end{array}
\right. \label{SD}
\end{align}
where ${u}^{k}$ and ${\Theta}^{k}$ are approximations to ${u}$ and ${\Theta}$ at ${t}=k\Delta{t}~(k = 0,1,2,\cdots)$. At each  time step $k=1,2,\cdots$, we solve \eqref{SD} with given boundary and initial conditions \eqref{Boundary-Cond} using the finite element method.  The details of the weak forms for \eqref{SD} and their unique solvability are described in \ref{AppendixA}.

In observation area $\mathcal{A}$ illustrated in Figure \ref{fig:3}, we define the average of (thermo)elastic energy densities in $\mathcal{A}$ as follows:
\begin{align}
\mathcal{W}(\mathcal{A}) &:= \frac{1}{|\mathcal{A}|}\int_{\mathcal{A}} {W}(u)~dx,\notag \\
\mathcal{W}^{*}(\mathcal{A}) &:= \frac{1}{|\mathcal{A}|}\int_{\mathcal{A}} {W}^{*}(u,\Theta)~dx, \notag
\end{align}
and the differences between $\mathcal{W}(\mathcal{A})$ and $\mathcal{W}^{*}(\mathcal{A})$ for each $\delta > 0$ and for $\delta = 0$ are defined by
\begin{subequations}
	\begin{empheq}[]{align}
	&\Delta \mathcal{W}(\mathcal{A}) := \mathcal{W}(\mathcal{A})\big|_{\delta} - \mathcal{W}(\mathcal{A})\big|_{\delta=0}~, \notag\\
	&  \Delta \mathcal{W}^{*}(\mathcal{A}) := \mathcal{W}^{*}(\mathcal{A})\big|_{\delta} - \mathcal{W}^{*}(\mathcal{A})\big|_{\delta=0}~. \notag
	\end{empheq}
\end{subequations}\\[-3pt] 
In the following examples, we use the software FreeFEM \cite{Hecht2012} with P2 elements and unstructured meshes. For the time interval and time step, we use $0 \leq t \leq 0.1$ and  $\Delta t = 1\times 10^{-4}$, respectively.

\subsubsection{L-shape cantilever}\label{subsubsec:2.4.1}
Here, we consider the L-shaped cantilever whose left side is fixed, and the vertical displacement $u_{2}$ is given on the right side, as illustrated in Figure \ref{fig:3} (left). We denote the left and right boundaries by $\Gamma_{D}^{u}$ and $\Gamma_{DN}^{u}$, respectively, and define $\Gamma_{N}^{u} := \Gamma \setminus (\Gamma_{D}^{u} \cup \Gamma_{DN}^{u})$. The boundary conditions for $u$ are
\begin{align}
& u = 0,~ \mbox{on}~ \Gamma_{D}^{u}, \qquad \left\{
\begin{array}{ll}
\sigma_{11}^{*}[u,\Theta]n=0,\\
u_{2} = -0.1t
\end{array}
\right. ~ \mbox{on}~ \Gamma_{DN}^{u}, \qquad \sigma^{*}[u,\Theta]n = 0~ \mbox{on}~ \Gamma_{N}^{u}. \notag
\end{align}
For $\Theta$, we suppose $\frac{\partial\Theta}{\partial n}=0$ on $\Gamma$ and the initial temperature $\Theta_{*} = 0$. Although we adopt the above slightly modified boundary conditions in this example, the previous arguments are valid with small modifications, and we omit their details.

We apply the finite element method to \eqref{SD}. The  total number of triangular meshes $= 18215$ and the number of nodes (the vertices of the triangles) $=9301$. 

\begin{figure}[!h]
	\begin{center}
		\begin{tabular}{ccc}
			{\includegraphics[trim=8cm 0.5cm 2.5cm 0cm, scale=0.041]{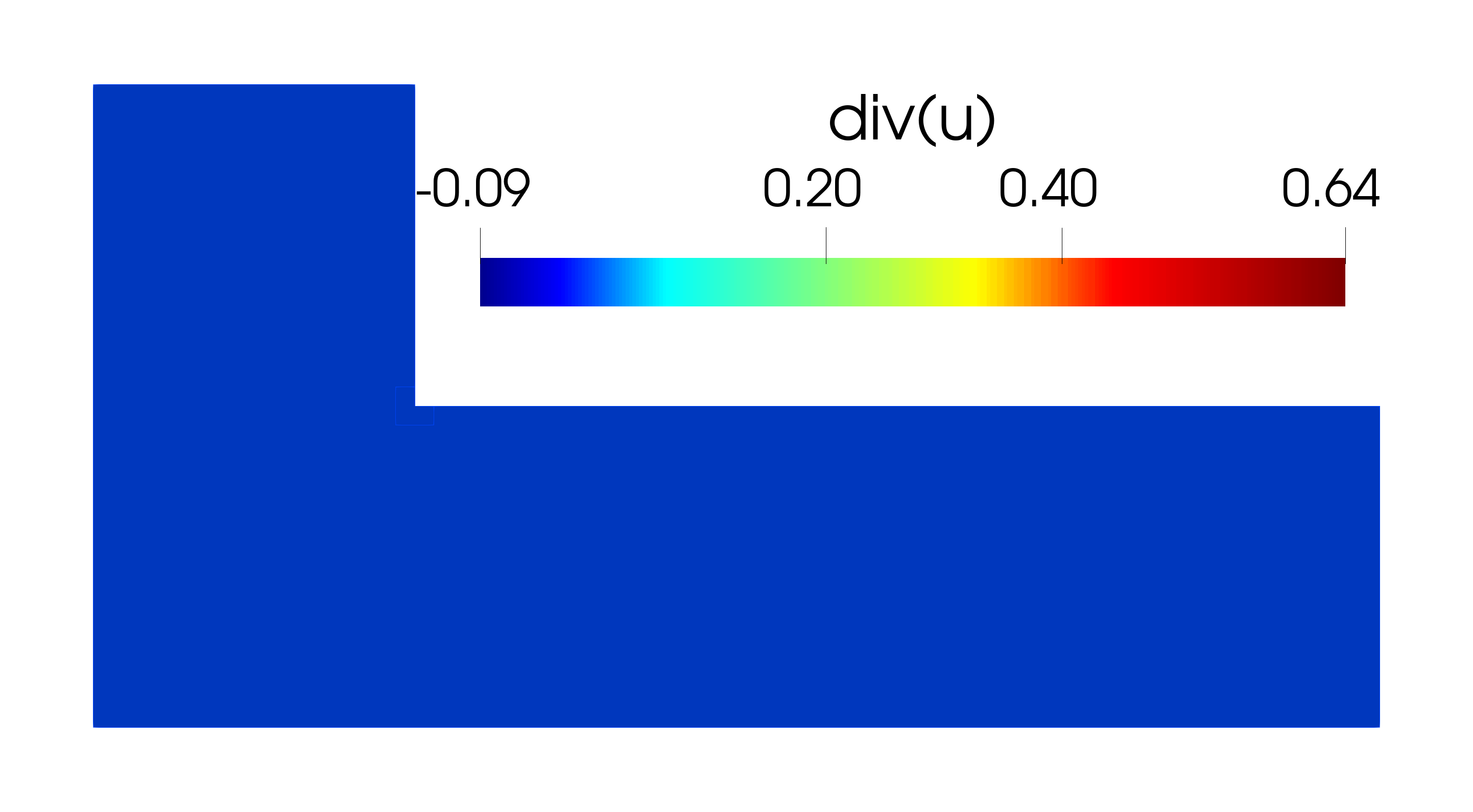}} &\hspace{-10pt}
			{\includegraphics[trim=8cm 0.5cm  2.5cm 0cm, scale=0.041]{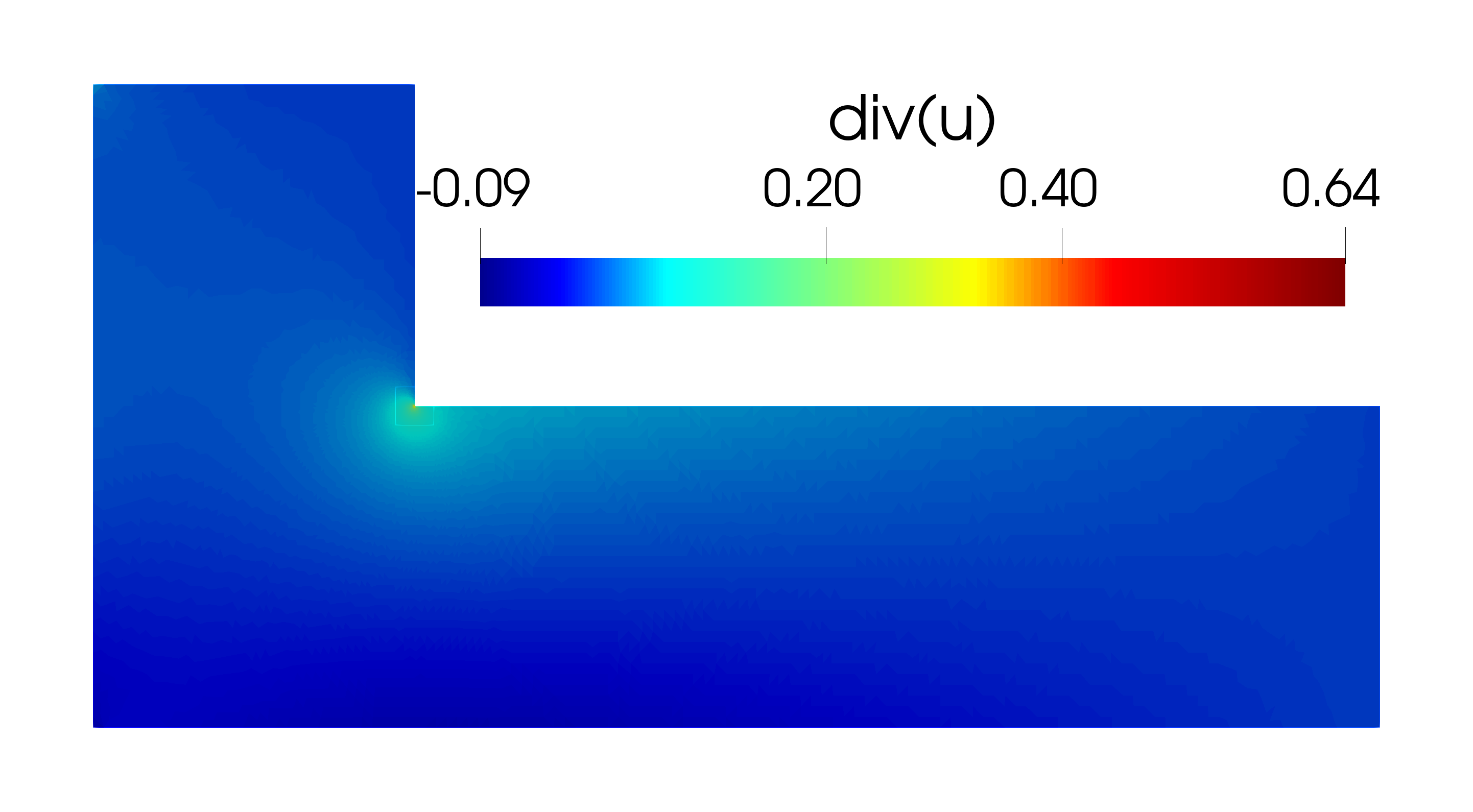}} &\hspace{-10pt}
			{\includegraphics[trim=8cm 0.5cm  2.5cm 0cm, scale=0.041]{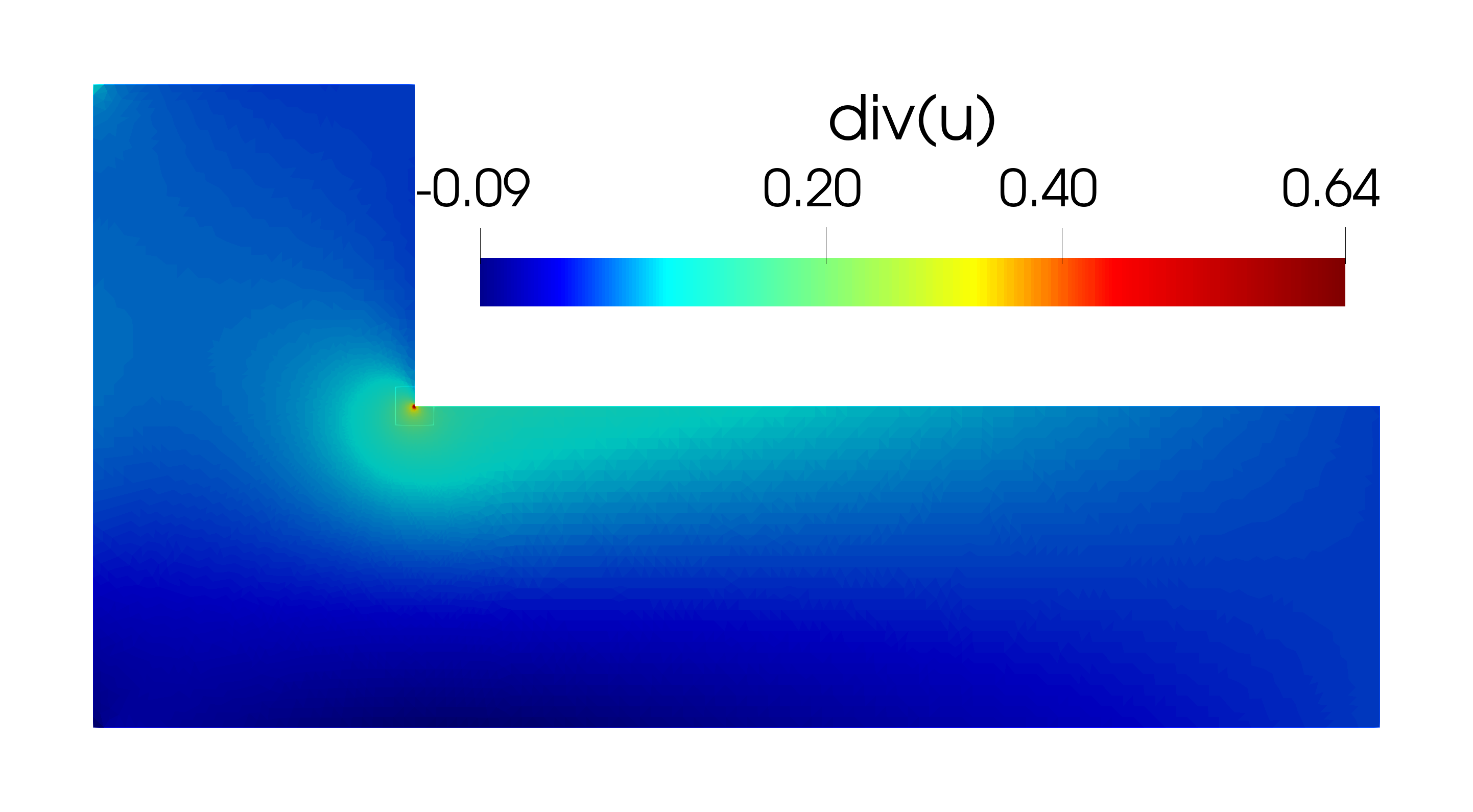}}
			\\[-10pt]
			{\includegraphics[trim=8cm 0.5cm 2.5cm 0cm, scale=0.04]{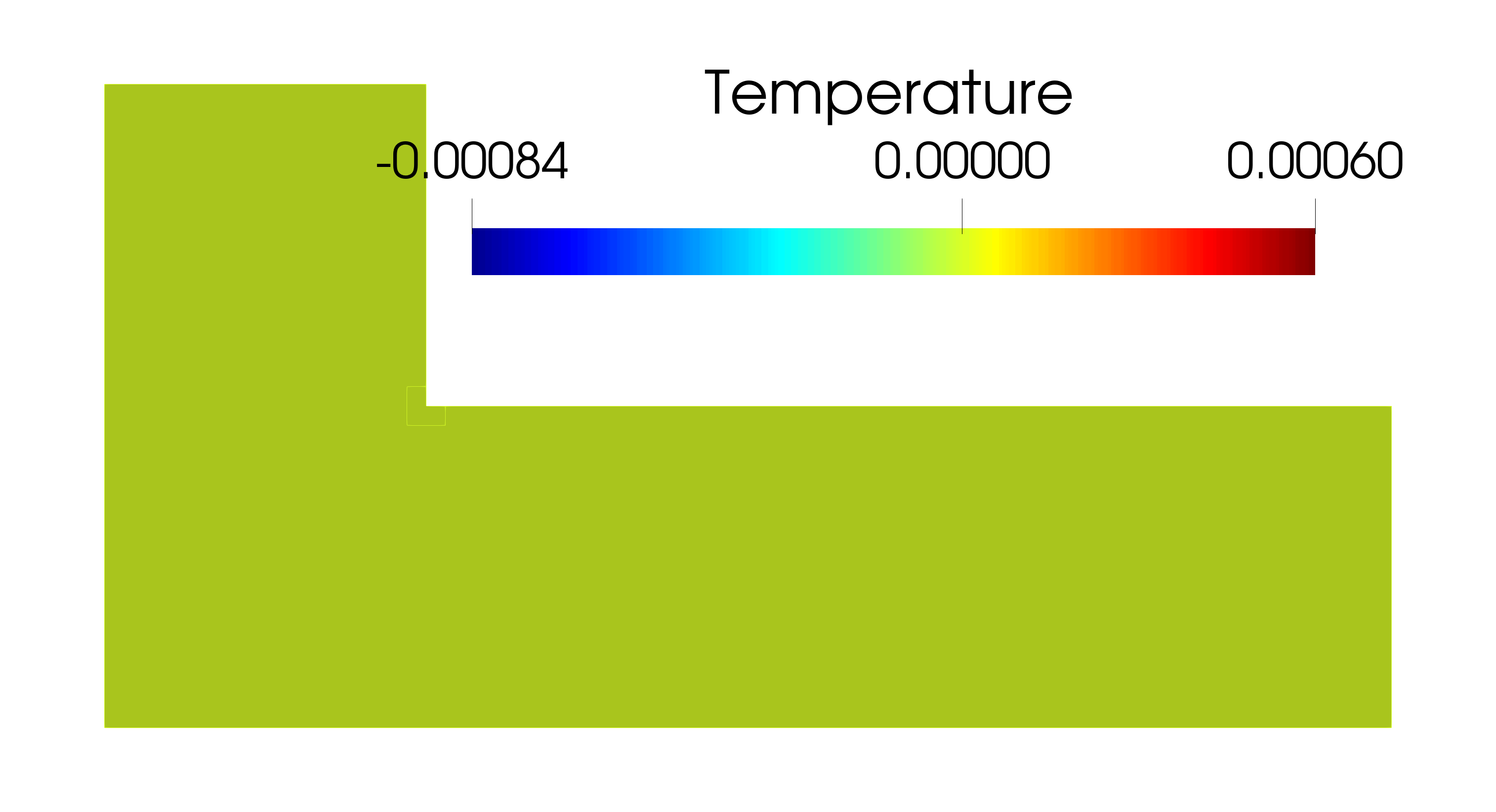}} &\hspace{-10pt}
			{\includegraphics[trim=8cm 0.5cm 2.5cm 0cm, scale=0.04]{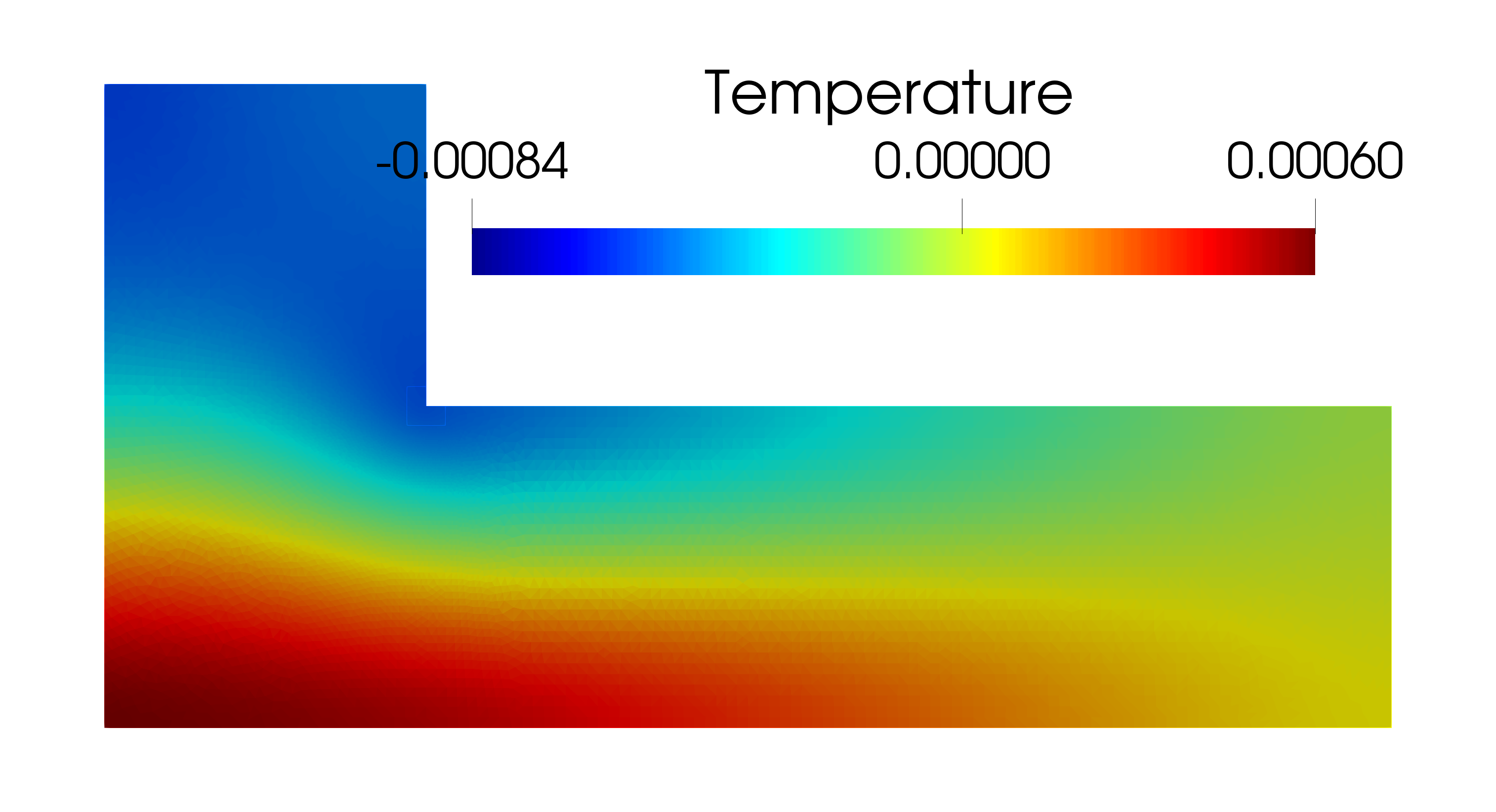}} &\hspace{-10pt}
			{\includegraphics[trim=8cm 0.5cm 2.5cm 0cm, scale=0.04]{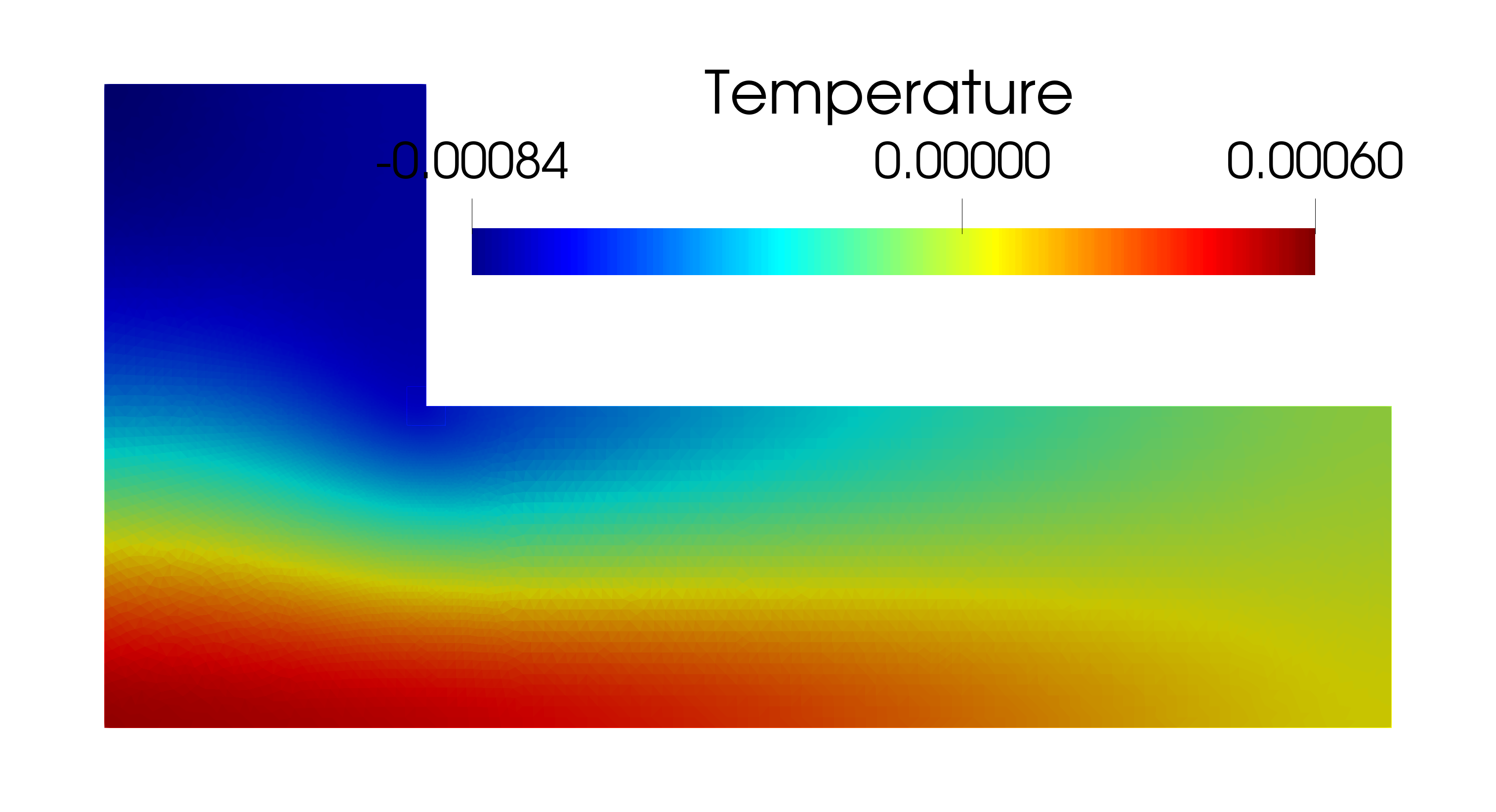}}\\[-8pt]
			${{t} = 0}$ & ${{t} = 0.05}$ & ${{t} = 0.1}$
		\end{tabular}\\[-16pt]
		\caption{Snapshots of $\mbox{div}u$ (upper) and the temperature (lower) of the L-shape cantilever for $t = 0, ~0.05, ~0.1$ using $\delta = 0.1$. Near the re-entrant corner, the domain is expanded ($\mbox{div}u > 0$), and the temperature decreases. On the other hand, near the bottom boundary, the domain is compressed ($\mbox{div}u < 0$), and the temperature increases.}\label{fig5}
	\end{center}
\end{figure}
\begin{figure}[!h]
	\centering%
	\floatbox{figure}{%
		\begin{subfloatrow}[3]
			\ffigbox[\FBwidth]{%
				\includegraphics[trim=4cm 4cm 4cm 7cm, width=0.9\linewidth]{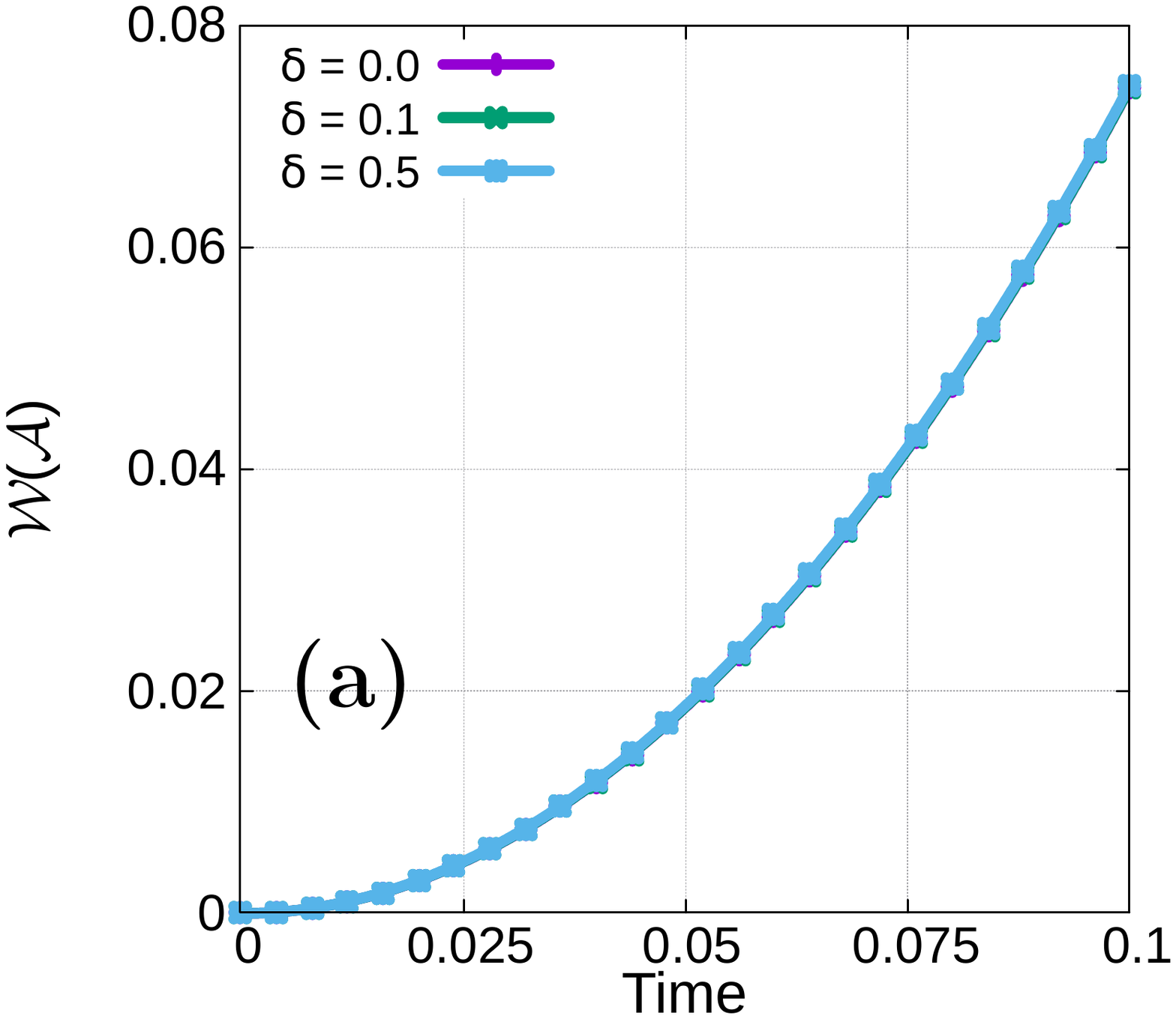}}{\caption*{}}~
			\ffigbox[\FBwidth]{%
				\includegraphics[trim=4cm 4cm 4cm 7cm, width=0.9\linewidth]{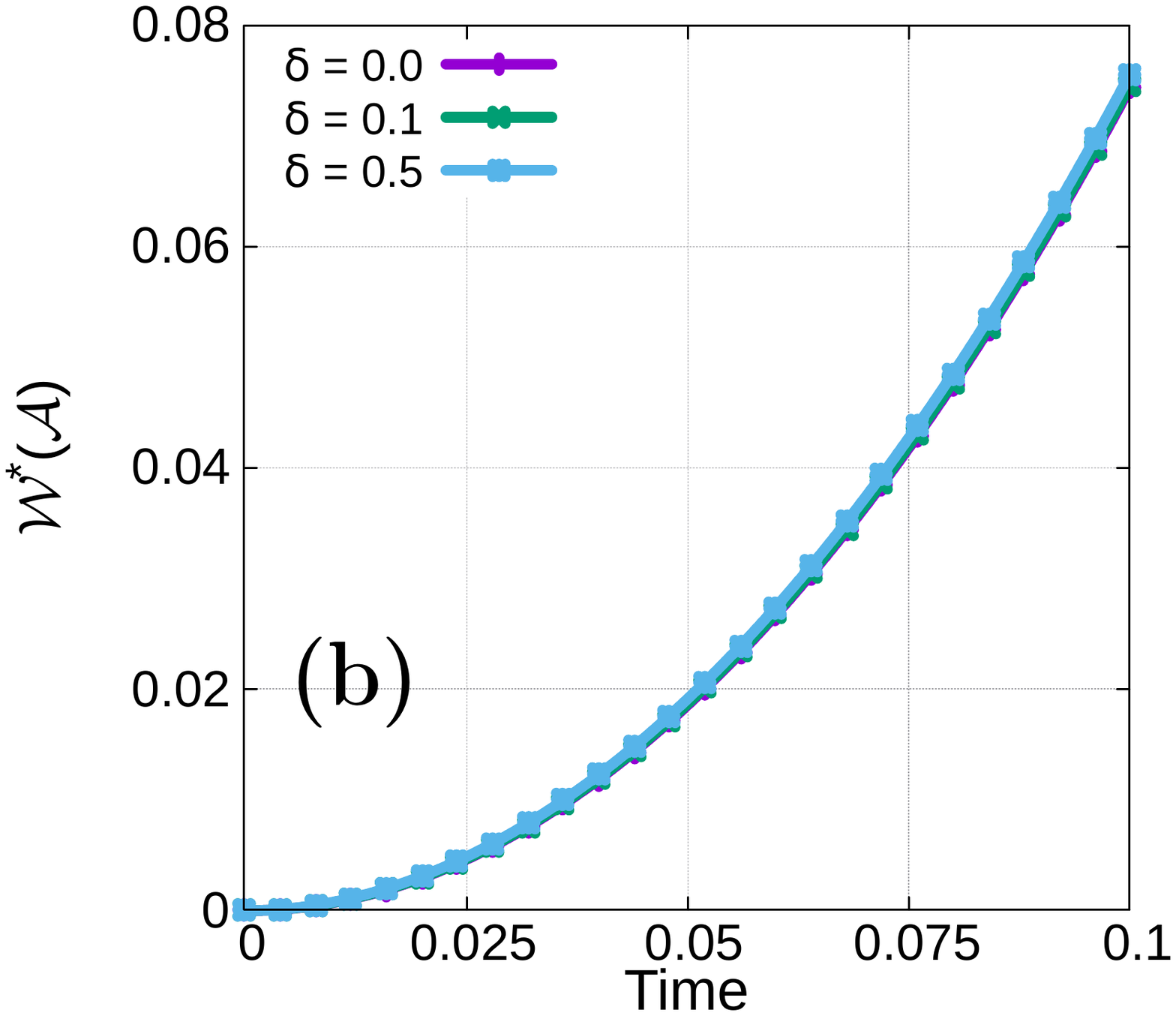}}{\caption*{}}~
			\ffigbox[\FBwidth]{%
				\includegraphics[trim=4cm 4cm 4cm 7cm, width=0.9\linewidth]{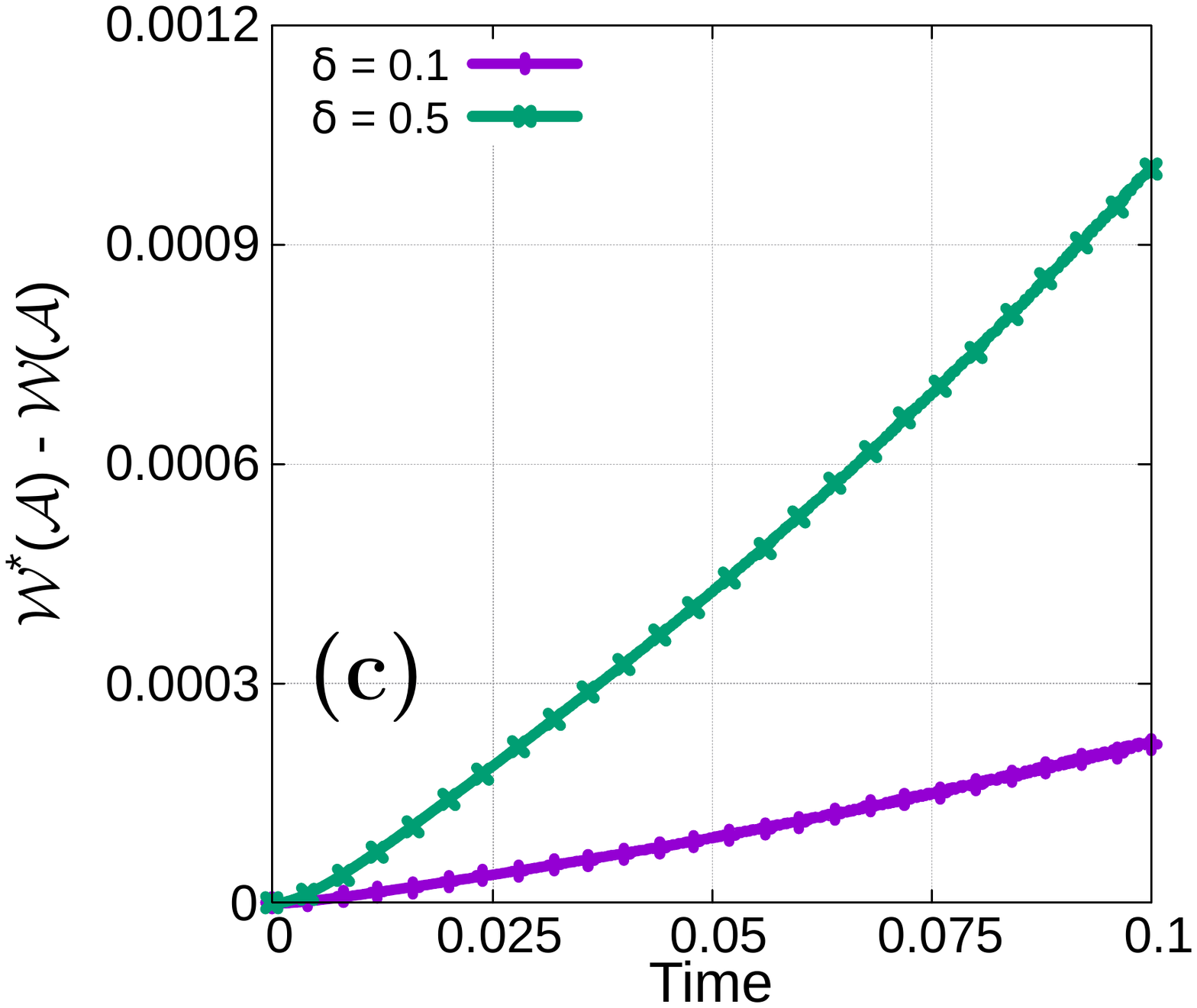}}{\caption*{}}
		\end{subfloatrow}
		\begin{subfloatrow}[2]
			\ffigbox[\FBwidth]{%
				\includegraphics[trim=4cm 7cm 4cm 10cm, width=0.73\linewidth]{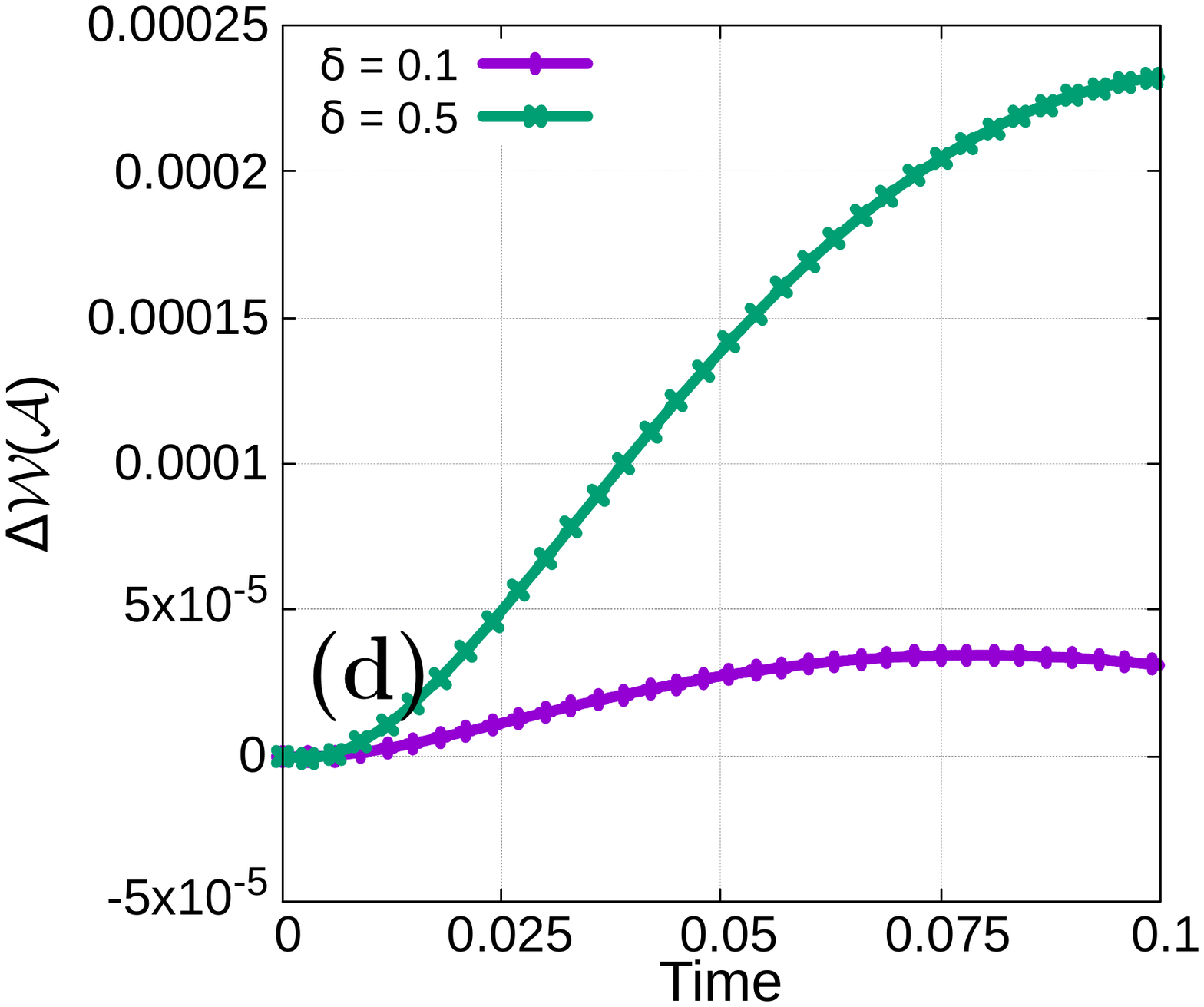}}{\caption*{}}
			\ffigbox[\FBwidth]{%
				\includegraphics[trim=4cm 7cm 4cm 10cm, width=0.73\linewidth]{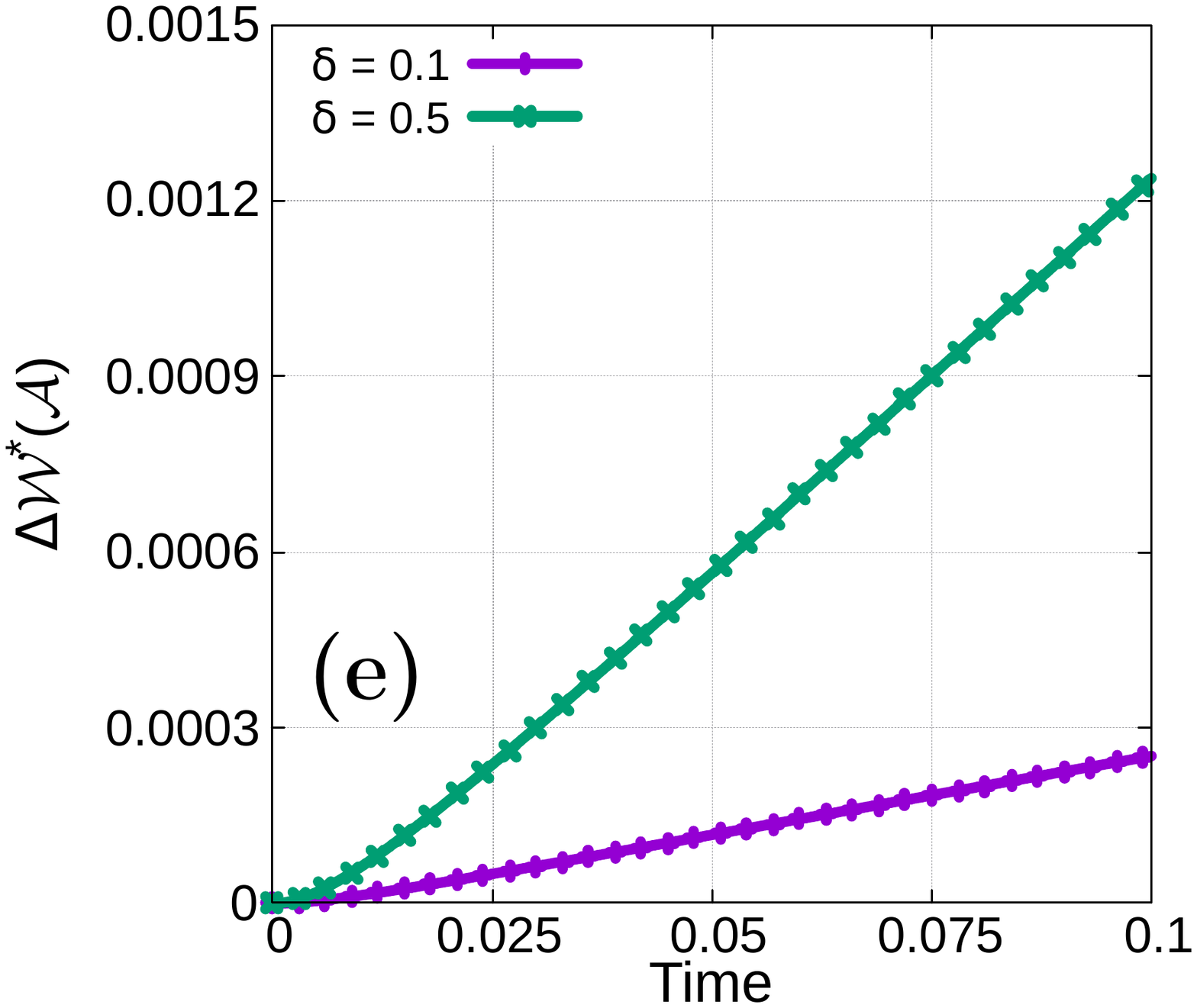}}{\caption*{}}
		\end{subfloatrow}\\[-30pt]
	}
	{%
		\caption{Profiles of (a) $\mathcal{W}(\mathcal{A})$, (b) $\mathcal{W}^{*}(\mathcal{A})$, (c) $\mathcal{W}^{*}(\mathcal{A})-\mathcal{W}(\mathcal{A})$, (d) $\Delta \mathcal{W}(\mathcal{A})$ and (e) $\Delta \mathcal{W}^{*}(\mathcal{A})$ in an L-shaped cantilever during the loading process.}
		\label{fig:sigmaustar}
	}
\end{figure}

As shown in the lower part of Figure \ref{fig5}, we observe that the highest temperature is in the contracting area and the lowest is in the expanding area. Furthermore, there exists a contribution $\delta$ for each $\delta > 0$ during  the loading process. Although the disparity is small, the thermoelastic coupling parameter $\delta$ contributes to the variations in $\mathcal{W}(u)$ and $\mathcal{W}^{*}(u)$, as shown in Figure \ref{fig:sigmaustar} (a)-(b). Here, a larger $\delta$ value implies larger $\mathcal{W}(\mathcal{A})$ and $\mathcal{W}^{*}(\mathcal{A})$ values (Figure \ref{fig:sigmaustar} (d)-(e)). In addition, we also observe that $\mathcal{W}^{*}(\mathcal{A})$ is larger than $\mathcal{W}(\mathcal{A})$ for each $\delta > 0$ (Figure \ref{fig:sigmaustar} (c)).

In the L-shape cantilever case for each $\delta > 0$, we conclude that the thermal coupling parameter enhances the singularity of (thermo)elastic energy in the expanding area. The (thermo)elastic energy plays a role in the driving force in the phase field model \cite{Miehe2015}, which means that the parameter $\delta$ can accelerate crack growth in the expanding area. 

\subsubsection{Cracked domain}\label{subsubsec:2.4.2}
Here, we consider a cracked domain with vertical displacements on the top and bottom sides, and the other sides are free traction, as shown in Figure \ref{fig:3} (right). The boundary conditions for $u$ are 
\begin{align}
\left\{
\begin{array}{ll}
u_{1}=0,\\
u_{2} = \pm t
\end{array}
\right. ~ \mbox{on}~ \Gamma_{\pm D}^{u}, \qquad \sigma^{*}[u,\Theta]n = 0~ \mbox{on}~ \Gamma_{N}^{u}, \notag
\end{align}
where $\Gamma_{+D}^{u}$ and $\Gamma_{-D}^{u}$ denote the top and bottom boundaries of $\Omega$, respectively, and $\Gamma_{N}^{u}:=\Gamma\setminus(\Gamma_{+D}^{u} \cup \Gamma_{-D}^{u})$. For $\Theta$, we suppose $\frac{\partial\Theta}{\partial n}=0$ on $\Gamma_{N}^{\Theta}=\Gamma$ and the initial temperature $\Theta_{*} = 0$. 

We use the finite element method to solve \eqref{SD}. Therefore, the  total number of triangular meshes and the number of nodes (the vertices of the triangles) are $11176$ and $5722$, respectively. 
\begin{figure}[!h]
	\begin{center}
		\begin{tabular}{cc}
			\includegraphics[trim=4cm 6cm 0cm 0cm, width=0.45\linewidth]{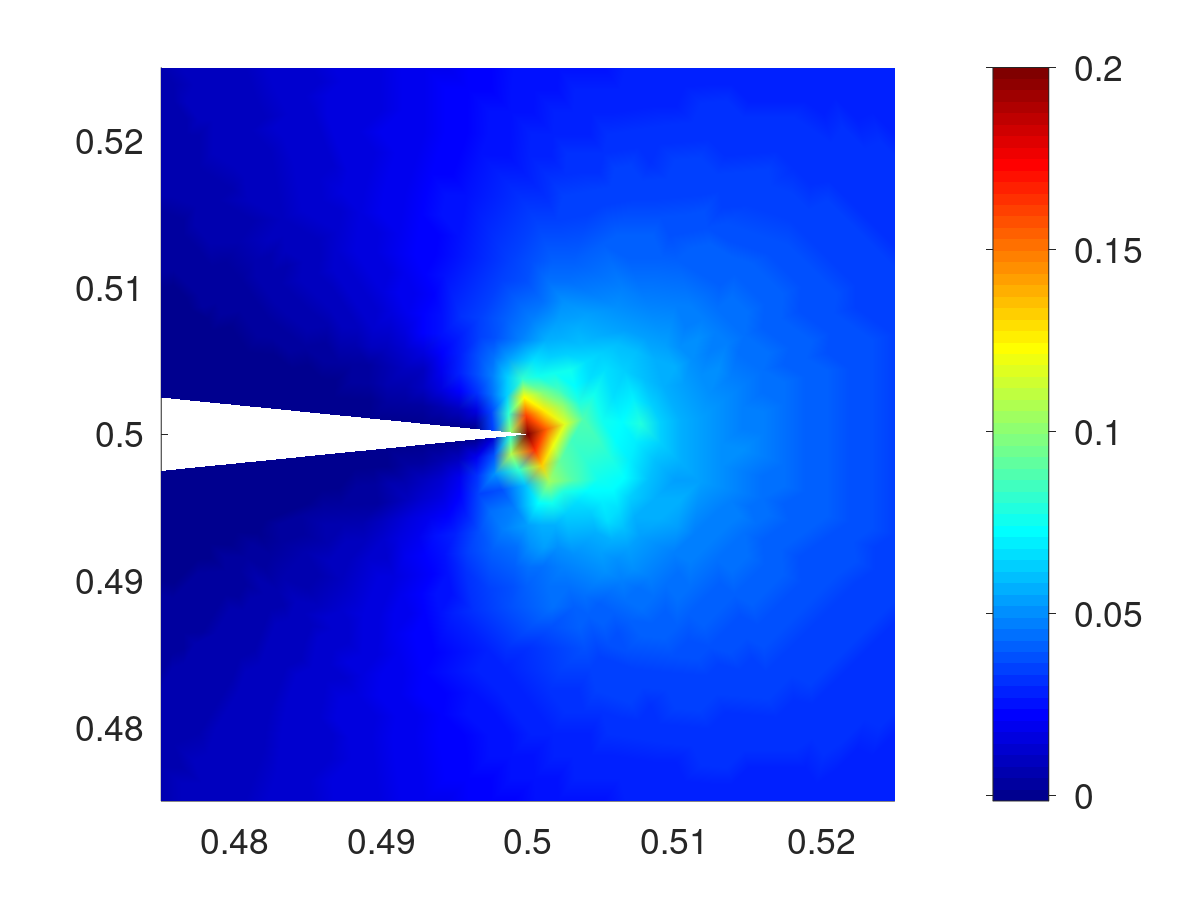}&
			\hspace{8pt}\includegraphics[trim=4cm 6cm 0cm 0cm, width=0.45\linewidth]{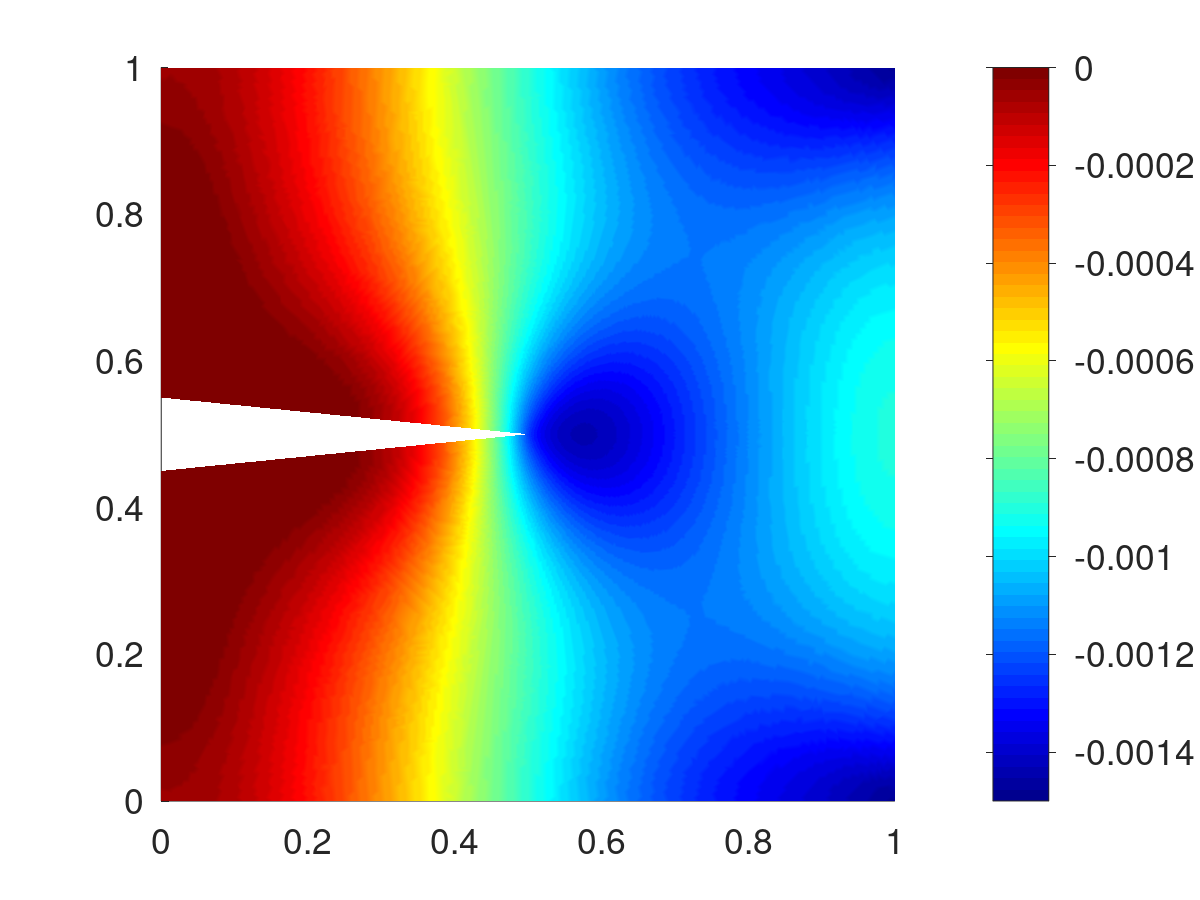}
		\end{tabular}\\[8pt]
		\caption{Snapshot of $\mbox{div}u$ on the subdomain $\mathcal{A}$ (left) and temperature $\Theta$ in $\Omega$ (right) using $\delta = 0.1$  at $t = 0.1$.}
		\label{fig:3Dblunt0}
	\end{center}
\end{figure}
\begin{figure}[!h]
	\begin{center}
		\begin{tabular}{cc}
			\includegraphics[trim=2cm 3cm 0.1cm 0cm, width=0.5\linewidth]{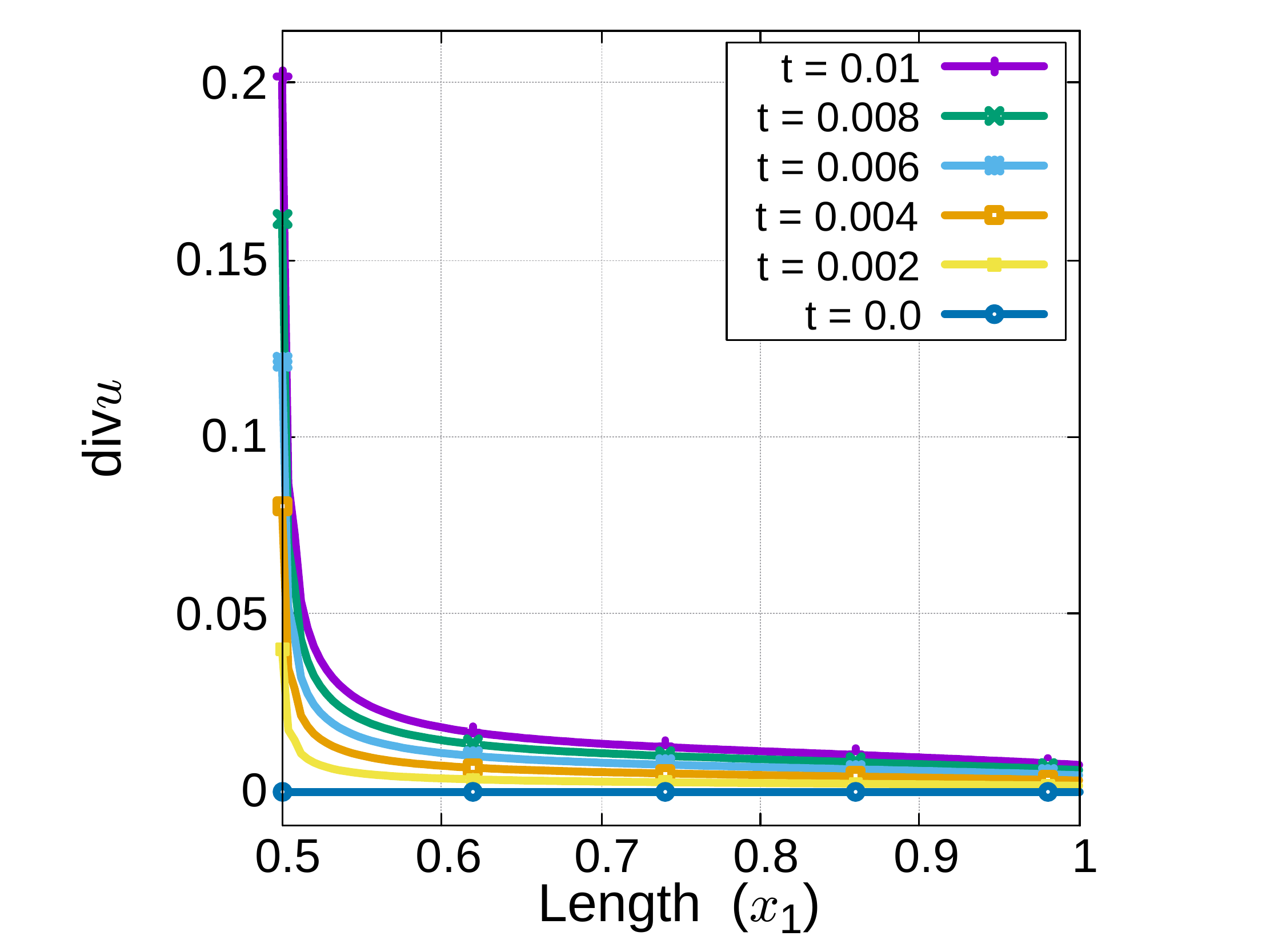}&
			\includegraphics[trim=2cm 3cm 0.1cm 0cm, width=0.5\linewidth]{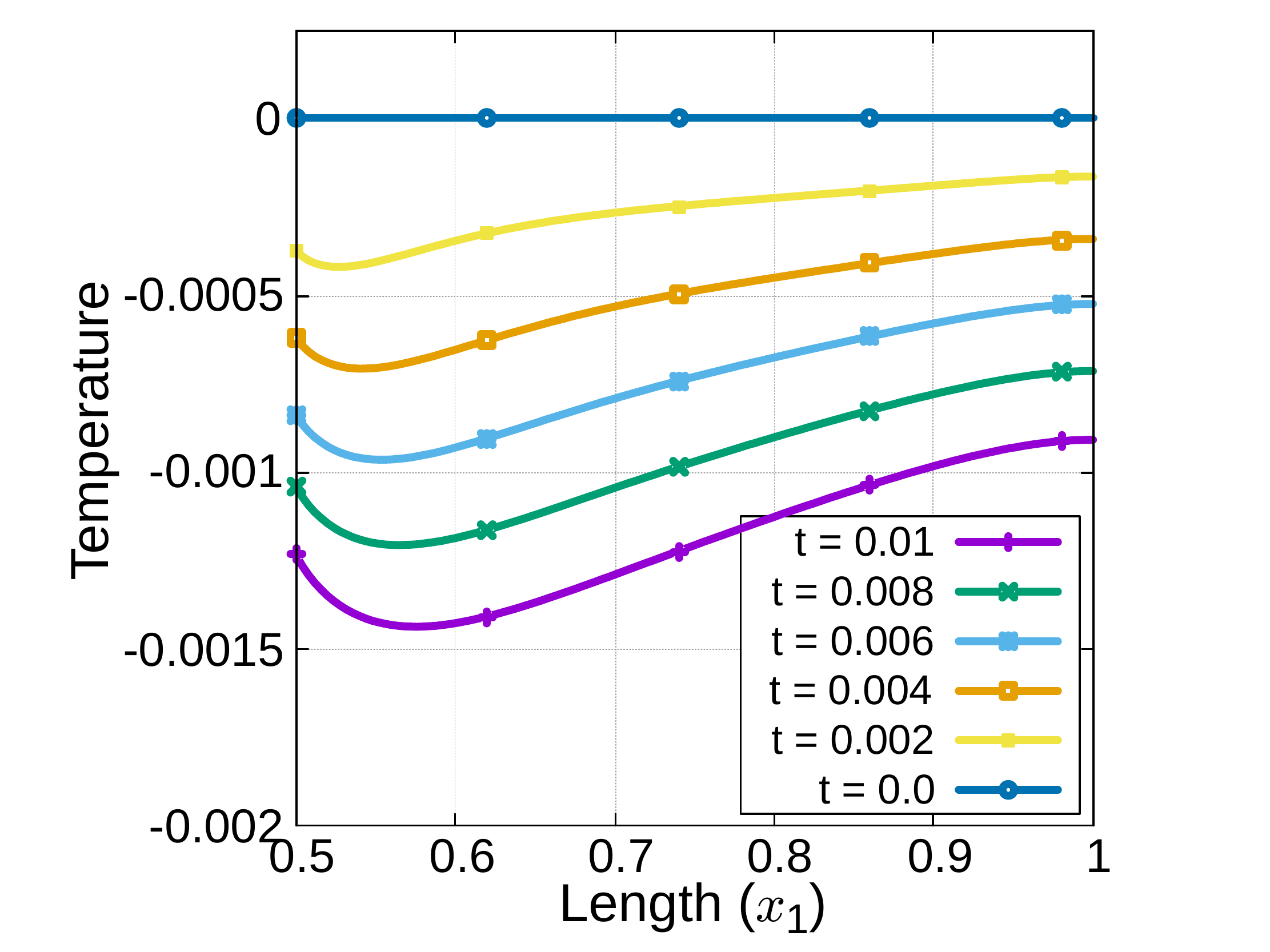}
		\end{tabular}\\[8pt]
		\caption{Profile of $\mbox{div}u$ (left) and temperature $\Theta$ (right) using $\delta = 0.1$ along the $x_{1}$ axis, i.e., $x_{2}=0$, $0.5 \leq x_{1} \leq 1$, during the loading process.}
		\label{fig:3Dblunt}
	\end{center}
\end{figure}

From Figure \ref{fig:3Dblunt0} (left), we conclude that the area that expands the most (i.e., $\mbox{div}u$ is largest) appears near the crack tip. This can be compared with the analytical solution for the linear elasticity in a cracked domain in \ref{AppendixB}. We also observe that the region with the lowest temperature appears to the right of the crack tip in Figure \ref{fig:3Dblunt0} (right). From the temporal change in the temperature along the $x_{1}$ axis plotted in Figure \ref{fig:3Dblunt} (right), we also observe that the lowest temperature region appears in $0.5 < x_{1} < 0.6$ and that the temperature decreases over time. This is shown in Figure \ref{fig:3Dblunt} (left), where the value of $\mbox{div}u$ is plotted along the $x_{1}$ axis and $\mbox{div}u$ is increasing over time; i.e., the heat source term $\mbox{div}\dot{u}$ in \eqref{Biot-Model2} is positive. 

\begin{figure}[!h]
	\begin{center}
		\begin{tabular}{cc}
			\includegraphics[trim=2cm 2cm 1cm 0cm, width=0.5\linewidth]{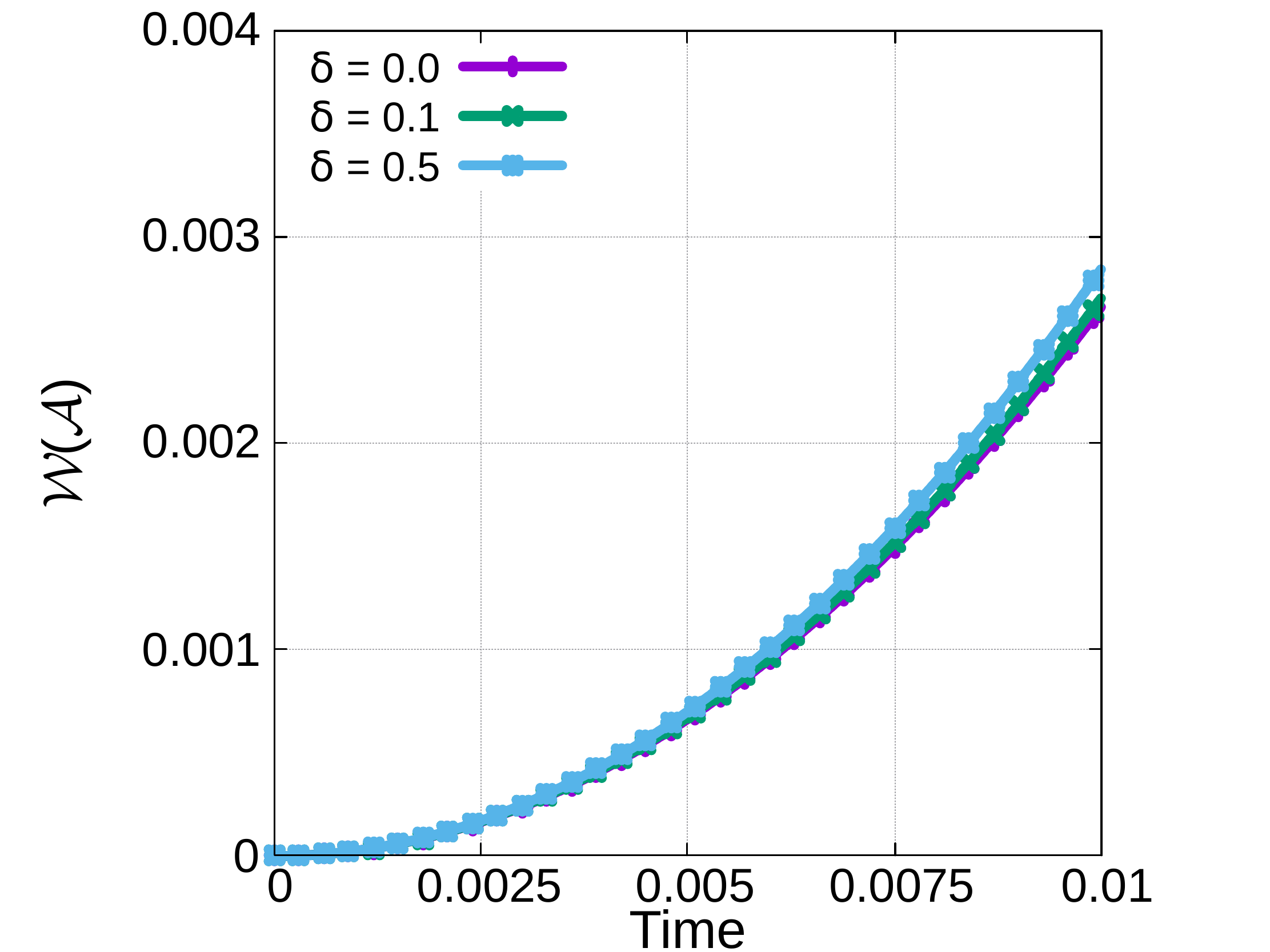}
			&
			\includegraphics[trim=2cm 2cm 1cm 0cm, width=0.5\linewidth]{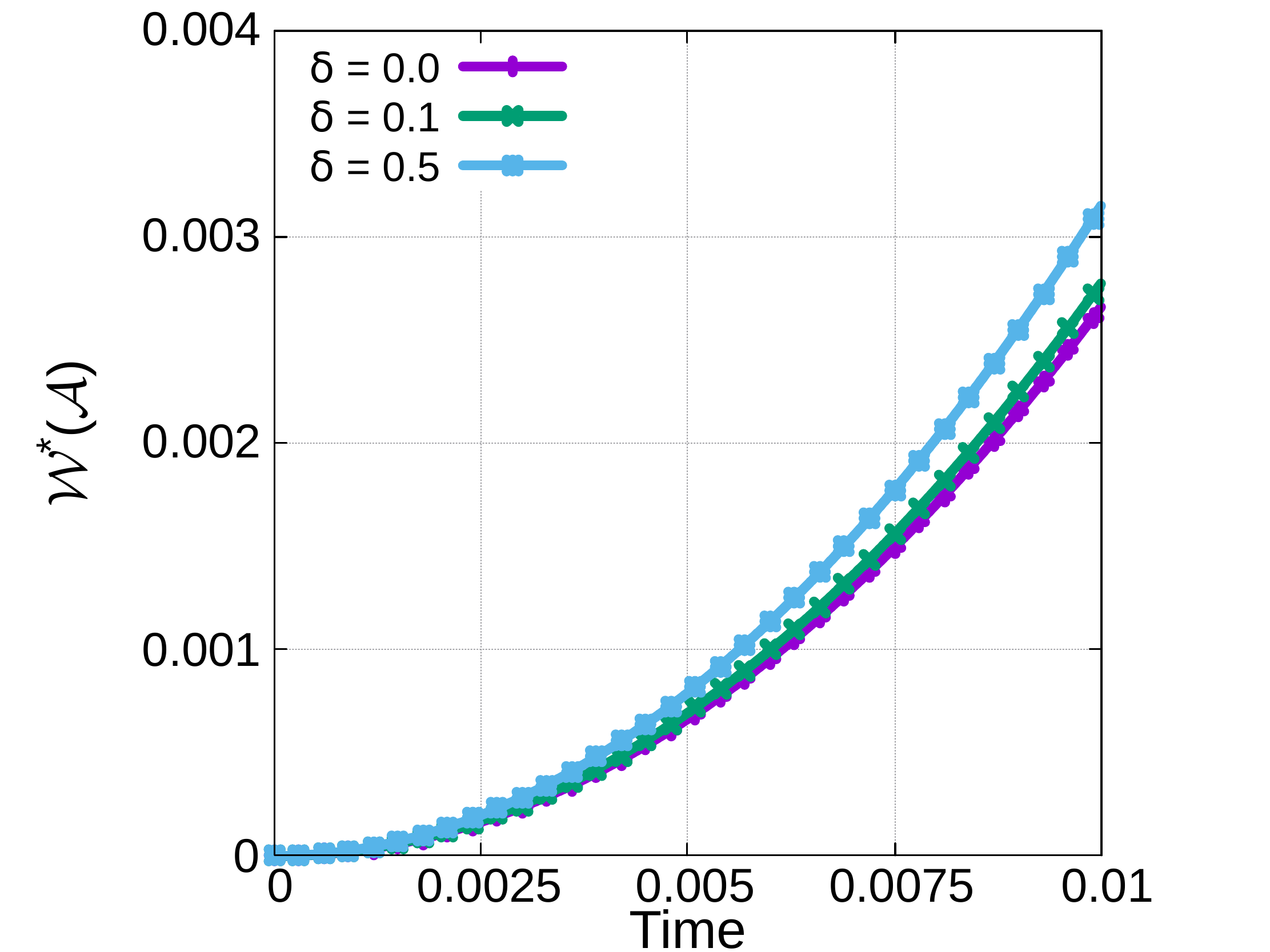}
		\end{tabular}
		\caption{Profiles of $\mathcal{W}(\mathcal{A})$ (left)  and  $\mathcal{W}^{*}(\mathcal{A})$ (right) in subdomain $\mathcal{A}$ during the loading process.}
		\label{fig:sigmaustar2}
	\end{center}
\end{figure}

Similar to Section \ref{subsubsec:2.4.1}, for each $\delta > 0$, we obtain variations of $\mathcal{W}(\mathcal{A})$ and $\mathcal{W}^{*}(\mathcal{A})$ in subdomain $\mathcal{A}$ (Figure \ref{fig:sigmaustar2}), where the subdomain $\mathcal{A}$ corresponds to the area that expands the most. From Figure \eqref{fig:sigmaustar2}, it is observed that $\mathcal{W}^{*}(\mathcal{A})$ is larger than $\mathcal{W}(\mathcal{A})$. This suggests that the thermoelastic energy density ${W}^{*}(u,\Theta)$ has a higher value than the elastic energy density ${W}(u)$. These observations are confirmed by the comparison of our thermal fracturing phase field models. 

\section{Crack Propagation under Thermal Stress} \label{sec:3}
This section is devoted to the phase field models for thermal fracturing, which are the main purpose of this paper.
\subsection{Fracturing phase field model (F-PFM)}\label{subsec:3.1}
\begin{figure}[!h]
	\begin{tabular}{cc}
		\includegraphics[trim=4cm 4cm 4cm 4cm,width=0.39\textwidth]{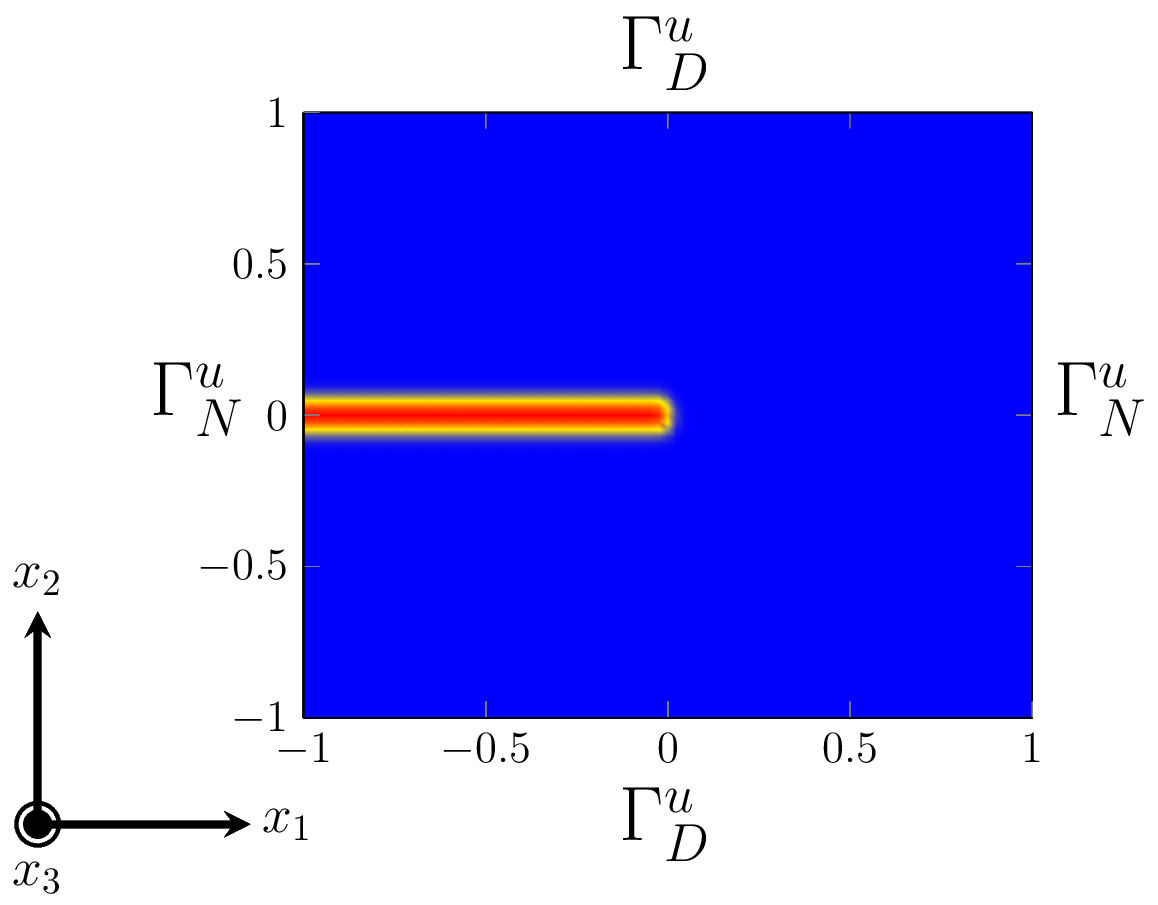}
		&\qquad\quad
		{\includegraphics[trim=4cm 4cm 4cm 4cm,width=0.42\textwidth]{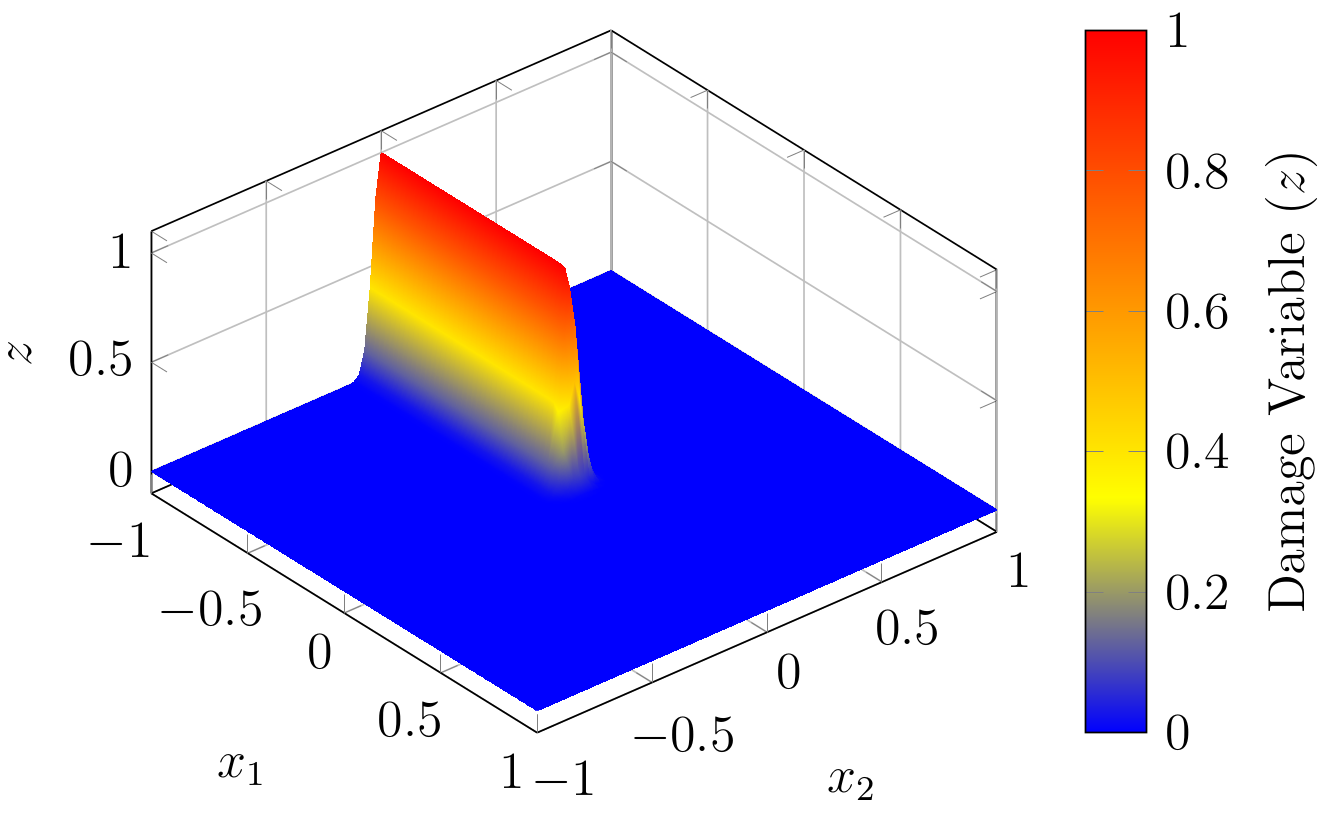}}
	\end{tabular}\\[12pt]
	\caption{Illustration of the phase field approximation of the cracked surface in an elastic body.} 
	\label{PhaseFiled:1}
\end{figure}

According to the works \cite{Kimura2009,Kimura2021}, we introduce fracturing PFM (we call it F-PFM) in this section. Let $\Omega$ be a bounded (uncracked) domain in $\mathbb{R}^{d}$ and $\Gamma := \partial\Omega = \Gamma_{D}^{u} \cup \Gamma_{N}^{u}$, similar to Section \ref{sec:2}. In F-PFM, a crack in $\Omega$ at time $t$ is described by a damage variable $z(x,t) \in [0,1]$ for $x \in \overline{\Omega}$ with space regularization. The cracked and uncracked regions are represented by $z \approx 1$ and $z\approx 0$, respectively, and $z\in (0,1)$ indicates slight damage.  A typical example of a straight crack in a square domain is illustrated in Figure \ref{PhaseFiled:1}. 

The F-PFM is described as: 
\begin{subequations}
	\begin{empheq}[left=\empheqlbrace]{align}
	& -{\mbox{div}}\left((1-z)^{2}{\sigma}[{u}]\right)	= 0  &  \mbox{in}~ {\Omega}\times[0,{T}],  \label{PML1}\\
	&  {\alpha} \frac{\partial z}{\partial {t}} = \left({\epsilon}~{\mbox{div}}\left({\gamma_{*}} {\nabla} z\right) - \frac{{\gamma_{*}}}{{\epsilon}}z + (1-z){W}(u)\right)_{+} & \mbox{in}~ {\Omega}\times[0,{T}], \label{PML2}
	\end{empheq}\label{PML}
\end{subequations} \\[-3pt]
with the following boundary and initial conditions: 
\begin{subequations}
	\begin{empheq}[left=\empheqlbrace]{align}
	&  u = u_{D}(x,t) &  &\mbox{on}~ \Gamma_{D}^{u}\times[0,{T}], \label{BC1}\\
	& \sigma[u]n = 0 &  & \mbox{on}~ \Gamma_{N}^{u}\times[0,{T}], \label{BC2}\\
	&  \frac{\partial z}{\partial n} = 0 & & \mbox{on}~\Gamma\times[0,{T}], \label{BC3}\\
	& z({x,0}) = z_{*}(x) &  & \mbox{in} ~ \Omega,\label{IC2}
	\end{empheq}\label{BC-IC}
\end{subequations}\\[-3pt]
where the displacement $u:\overline{\Omega} \times [0,T] \mapsto \mathbb{R}^{d}$ and the damage variable $z: \overline{\Omega} \times [0,T] \mapsto [0,1]$ are unknowns. The parameters $\alpha > 0$ and $\epsilon > 0$ are small numbers related to regularization in time and space, respectively. The critical energy release rate is denoted by $\gamma_{*}$ (which is often denoted by $\text{G}_{\text{c}}$), and the elastic energy density is defined by ${W} = {W}(u) := \sigma[u]:e[u]$. In \eqref{PML2}, the term ${W}$ works as a driving force for $z$. 

The symbol $(~)_{+}$ on the right-hand side in \eqref{PML2} denoted the positive part $(s)_{+} := \text{max}(s,0)$, and it represents the irreversible property of crack growth. 

F-PFM is derived as a unidirectional gradient flow of the total energy $\mathcal{E}_{el}(u,z) + E_{s}(z)$, where
\begin{align}
\mathcal{E}_{el}(u,z) & := \frac{1}{2}\int_{\Omega} (1-z)^{2}\sigma[u]:e[u]~dx, \label{Modify-Elastic}\\
E_{s}(z) & := \frac{1}{2} \int_{\Omega}\gamma_{*}\left(\epsilon\big|\nabla z\big|^{2} + \frac{|z|^{2}}{\epsilon}\right)dx.\label{Surface-Energy}
\end{align}
More precisely, $u(t)$ obeys the following variational principle:
\begin{equation}
u(t) = \operatorname*{argmin}_{u \in V(u_{D}(t))} {\mathcal{E}_{el}(u,z(t))}, \label{GradientFlow1}
\end{equation}
and \eqref{PML2} becomes a gradient flow of the energy \(\min\limits_{u} \mathcal{E}_{el}(u,z) + E_{s}(z) \). 

We remark that $\mathcal{E}_{el}(u,z)$ is a modified elastic energy, which corresponds to the elastic energy with a damaged Young's modulus $\tilde{E}_{\text{Y}} = (1-z)^{2}E_{\text{Y}}$.  The energy $E_{s}(z)$ is regularized surface energy, which approximates the crack area ($d=3$) or length ($d=2$)  as $\epsilon \mapsto 0$. Please see \cite{Kimura2021} for more details.  The following energy equality for F-PFM is shown in \cite{Kimura2021} (\cite{Kimura2009} for the antiplane setting).

\begin{Theorem}[Energy equality for F-PFM]\label{Theorem2:a}
Let $(u(x,t), z(x,t))$ be a sufficiently smooth solution to \eqref{PML} and \eqref{BC-IC}. If $u_{D}$ is  independent of $t$, then we have
\begin{eqnarray}
\frac{d}{dt}\left(\mathcal{E}_{el}(u(t),z(t)) + E_{s}(z(t))\right) = -\alpha\int_{\Omega} \left|\dot{z}\right|^{2}dx \leq 0. \label{Eq:Theorem2:a}
\end{eqnarray}
\normalfont	\begin{Proof}
\normalfont Differentiating the total energy in $t$ and applying integration by parts,  we obtain
\begin{align}
		& \hspace{-12pt}\frac{d}{dt}\left(\mathcal{E}_{el}(u(t),z(t)) + E_{s}(z(t))\right) \nonumber\\
		&= \int_{\Omega}(1-z)^{2}\sigma[u]:{e}[\dot{u}]~dx + \int_{\Omega}\left(\gamma_{*}\epsilon\nabla z\cdot\nabla\dot{z} + \big(\frac{\gamma_{*}}{\epsilon}z - (1-z){W}(u)\big)\dot{z}\right)~dx\notag\\
		& = \int_{\Gamma} (1-z)^{2}\underbrace{\left(\sigma[u]n\right)}_{0}\cdot\dot{u}~ds - \int_{\Omega} \underbrace{\mbox{div}\left((1-z)^{2}\sigma[u]\right)}_{0}\cdot\dot{u}~dx \notag\\
		& \quad + \int_{\Gamma} \gamma_{*}\epsilon\underbrace{\frac{\partial z}{\partial n}}_{0}\dot{z}~ds - \int_{\Omega} \mathcal{H}\dot{z}~dx, \label{EnergyDisipation21a}
\end{align}
where we define $\mathcal{H} := \epsilon \mbox{div}\left(\gamma_{*}\nabla z\right) - \frac{\gamma_{*}}{\epsilon}z + (1-z){W}(u)$. Since \eqref{PML2} is written as $\alpha\dot{z} = (\mathcal{H})_{+}$, using the equality $\mathcal{H}(\mathcal{H})_{+} = (\mathcal{H})_{+}^{2}$, we conclude that
\begin{align}
		& \frac{d}{dt}\left(\mathcal{E}_{el}(u(t),z(t)) + E_{s}(z(t))\right) \notag\\
		& \qquad = - \int_{\Omega} \mathcal{H}\dot{z}~dx = -\int_{\Omega} \mathcal{H}\frac{(\mathcal{H})_{+}}{\alpha}~dx = -\int_{\Omega}\frac{(\mathcal{H})_{+}^{2}}{\alpha}~dx =  -\int_{\Omega}\alpha\left|\dot{z}\right|^{2}dx \notag. \qquad \text{\qed}
\end{align}
\end{Proof}
\end{Theorem}

\subsection{Thermal fracturing phase field model 1 (TF-PFM1)} \label{subsec:3.2}
To combine the Biot model in \eqref{BiotModel} and F-PFM in \eqref{PML}, their variational principles for $u$,  Proposition \ref{Proposition:2} and \eqref{GradientFlow1},  suggest that we consider the following modified thermoelastic energy: 
\begin{eqnarray}
\mathcal{E}_{el}^{*}(u,\Theta,z) &:=& \frac{1}{2}\int_{\Omega} (1-z)^{2}{\sigma}^{*}[u,\Theta]:{e}^{*}[u,\Theta]~dx,\label{Thermo-Elastic}
\end{eqnarray}
and a variational principle:
\begin{equation}
u(t) = \operatorname*{argmin}_{u \in V(u_{D}(t))} {\mathcal{E}_{el}^{*}(u,\Theta(t),z(t))}.\label{GradientFlow2}
\end{equation}
From the definition of the modified thermoelastic energy \eqref{Thermo-Elastic}, it is natural to replace the driving force term ${W}(u) = \sigma[u]:e[u]$ in \eqref{PML2} by the thermoelastic energy density ${W}^{*}(u,\Theta) := \sigma^{*}[u,\Theta]:e^{*}[u,\Theta]$. 

For heat equation \eqref{Biot-Model2}, since $\beta = a_{L}(d\lambda + 2\mu)$ and Lam\`{e}'s constants ($\lambda,~\mu$) are replaced by damaged constants ($(1-z)^{2}\lambda$, $(1-z)^{2}\mu$), $\beta$ should also be replaced by damaged constant $(1-z)^{2}\beta$. The thermal conductivity $\kappa_{0}$ is also considered  to be modified by $z$, because the heat is usually insulated across the crack. We suppose $\kappa = \kappa(z) > 0$ in this section, and we set it as $\kappa(z) = (1-z)^{2}\kappa_{0}$ in Section \ref{sec-4}.

Summarizing the above statements, we obtain the following thermal fracturing model, PFM 1 (TF-PFM1).
\begin{subequations}
	\begin{empheq}[left=\empheqlbrace]{align}
	&-{\mbox{div}}\left((1-z)^{2}{\sigma}^{*}[{u},\Theta]\right)	= 0  &\hspace{-16pt}\mbox{in}~ {\Omega}\times[0,{T}],  \label{TMF1a}\\
	&{\alpha} \frac{\partial z}{\partial {t}} = \left({\epsilon}~{\mbox{div}}({\gamma_{*}} {\nabla} z) - \frac{{\gamma_{*}}}{{\epsilon}}z + (1-z){W}^{*}(u,\Theta)\right)_{+} &\hspace{-16pt}\mbox{in}~ {\Omega}\times[0,{T}], \label{TMF1b}\\
	&\chi\frac{\partial {\Theta}}{\partial {t}} = ~{\mbox{div}}\left({\kappa(z)}{\nabla} {\Theta}\right) - \Theta_{0}(1-z)^{2}\beta\frac{\partial}{\partial {t}}({\mbox{div}} {u}) &\hspace{-16pt}\mbox{in}~ {\Omega}\times(0,{T}], \label{TMF1c}
	\end{empheq}\label{TMF1}
\end{subequations} \\[-3pt]
Similar to \eqref{BiotModel} and \eqref{PML}, the boundary and the  initial conditions to solve \eqref{TMF1} are presented as follows:
\begin{subequations}
	\begin{empheq}[left=\empheqlbrace]{align}
	& u = u_{D}(x,t)  & & \mbox{on}~ \Gamma_{D}^{u}\times[0,{T}],\label{TTBC1}\\
	& \sigma^{*}[u,\Theta]n = 0 & & \mbox{on}~ \Gamma_{N}^{u}\times[0,{T}], \label{TTBC2}\\
	& \Theta = {\Theta}_{D}(x,t) & &\mbox{on}~ {\Gamma}_{D}^{\Theta}\times[0,{T}], \label{TTBC3}\\
	& \frac{\partial \Theta}{\partial n} = 0& &\mbox{on}~ {\Gamma}_{N}^{\Theta}\times[0,{T}], \label{TTBC4}\\
	& \frac{\partial z}{\partial n} = 0 & & \mbox{on}~\Gamma \times[0,{T}], \label{TTBC5}\\
	& z(x,0) = z_{*}(x) & &\mbox{in} ~ \Omega, \label{TTIC1}\\
	& \Theta(x,0) = \Theta_{*}(x) & & \mbox{in} ~ \Omega.\label{TTIC2}
	\end{empheq}\label{TBC-IC1}
\end{subequations}\\[-3pt]
In the following, for simplicity, we define 
\begin{align}
 \sigma_{z}^{*}[u,\Theta] :=(1-z)^{2}\sigma^{*}[u,\Theta]. \notag
\end{align}

As a natural extension of Proposition \ref{Proposition:2} and Theorem \ref{Theorem1}, we obtain the following "partial" energy equality for TF-PFM1.
\begin{Theorem}[Energy equality for TF-PFM1]\label{Theorem2:b}
We suppose that $u_{D} \in H^{\frac{1}{2}}(\Gamma_{D}^{u};\mathbb{R}^{2})$ and $\Theta \in L^{2}(\Theta)$ are given and do not depend on $t$.  If $u(x,t)$ and $z(u,t)$ are sufficiency smooth and satisfy \eqref{TMF1a}, \eqref{TMF1b}, \eqref{TTBC1}, \eqref{TTBC2}, \eqref{TTBC5}, and \eqref{TTIC1}, the following energy equality holds:
\begin{eqnarray}
	\frac{d}{dt}\left(\mathcal{E}_{el}^{*}(u(t),\Theta,z(t)) + E_{s}(z(t))\right) = -\alpha \int_{\Omega} \left|{\dot{ z}}\right|^{2}dx \leq 0. \label{Eq:Theorem2:b}
\end{eqnarray}
\normalfont\begin{Proof}
\normalfont Under this condition, let us derive $\mathcal{E}_{el}^{*}(u(t),\Theta,z(t))$ and $E_{s}(z(t))$ with respect to $t$. 
\begin{align}
		& \hspace{-5pt}\frac{d}{dt}\left(\mathcal{E}_{el}^{*}(u(t),\Theta,z(t)) + E_{s}(z(t))\right) \nonumber\\
		& \hspace{-5pt}= \frac{1}{2}\frac{d}{dt} \int_{\Omega}\left( \sigma^{*}_{z}[u,\Theta]:e^{*}[u,\Theta]\right)dx  + \frac{1}{2}\frac{d}{dt} \int_{\Omega}\gamma_{*}\left(\epsilon\big|\nabla z\big|^{2} + \frac{|z|^{2}}{\epsilon}\right)~dx\nonumber\\
		&\hspace{-5pt}= \int_{\Omega}\sigma_{z}^{*}[u,\Theta]:e[\dot{u}]~dx + \int_{\Omega} \left(\gamma_{*}\epsilon\nabla z\cdot \nabla\dot{z} + \big(\frac{\gamma_{*}}{\epsilon}z - (1-z){W}^{*}(u,\Theta)\big)\dot{z}\right)dx \nonumber\\
		&\hspace{-5pt}= \int_{\Gamma} \underbrace{\sigma_{z}^{*}[u,\Theta]n}_{0}\cdot e[\dot{u}]~ds- \int_{\Omega}\underbrace{\mbox{div}\sigma_{z}^{*}[u,\Theta]}_{0}\cdot e[\dot{u}]~dx + \gamma_{*}\epsilon\int_{\Gamma} \underbrace{\frac{\partial z}{\partial n}}_{0}\dot{z}~ds \nonumber\\[-8pt]
		& \qquad \qquad - \int_{\Omega} \mathcal{H}^{*}\dot{ z} ~dx, \label{EnergyDisipation2a}
\end{align}
where we also define $\mathcal{H}^{*}:= \epsilon\mbox{div}(\gamma_{*}\nabla z) - \frac{\gamma_{*}}{\epsilon}z + (1-z){W}^{*}(u,\Theta)$. Since \eqref{TMF1b} is changed to $\alpha\dot{ z} = (\mathcal{H}^{*})_{+}$, similar to that in Section \ref{subsec:3.1}, we conclude that 
\begin{align}
		& \frac{d}{dt}\left(\mathcal{E}_{el}^{*}(u(t),\Theta,z(t)) + E_{s}(z(t))\right) = -\alpha\int_{\Omega}\left|\dot{z}\right|^{2}~dx \leq 0, \notag
\end{align}
which is equivalent to (\ref{Eq:Theorem2:b}). \qed
\end{Proof}
\end{Theorem}

\subsection{Thermal fracturing phase field model 2 (TF-PFM2)}\label{subsec:3.3}
In the previous section, we proposed TF-PFM1 based on the thermoelastic energy $E_{el}^{*}(u,\Theta)$. We proved a variational principle but proved only partial energy equality.
As shown in Section \ref{subsec:2.2}, the Biot model is related to both energies $E^{*}_{el}(u,\Theta)$ and $E_{el}(u)$. The variational principle holds for $E_{el}^{*}(u,\Theta)$ (Proposition \ref{Proposition:2}), and the energy equality holds for $E_{el}(u)$ (Theorem \ref{Theorem1}). This motivates us to consider another type of thermal fracturing PFM based on elastic energy $E_{el}(u)$. 
We call the following thermal fracturing model TF-PFM2:
\begin{subequations}
	\begin{empheq}[left=\empheqlbrace]{align}
	&  -{\mbox{div}}\left((1-z)^{2}{\sigma}^{*}[{u},{\Theta}]\right)	= 0 & \hspace{-16pt}\mbox{in}~ {\Omega}\times[0,{T}],  \label{TMF2a}\\
	&  {\alpha} \frac{\partial z}{\partial {t}} = \left({\epsilon}~{\mbox{div}}({\gamma_{*}} {\nabla} z) - \frac{{\gamma_{*}}}{{\epsilon}}z + (1-z){W}(u)\right)_{+} & \hspace{-16pt} \mbox{in}~ {\Omega}\times[0,{T}], \label{TMF2b}\\
	& \chi\frac{\partial {\Theta}}{\partial {t}} = ~{\mbox{div}}\left({\kappa(z)}{\nabla} {\Theta}\right) - \Theta_{0}(1-z)^{2}\beta\frac{\partial}{\partial {t}}({\mbox{div}} {u}) & \hspace{-16pt} \mbox{in}~ {\Omega}\times(0,{T}].\label{TMF2c}
	\end{empheq}\label{TMF2}
\end{subequations}\\[-3pt]
The associated boundary and initial conditions are given by \eqref{TBC-IC1}. For this model, we can show the following energy equality. 
\begin{Theorem}[Energy equality for TF-PFM2] \label{Theorem3}
We suppose that $(u(x,t),$ $\Theta(x,t), z(x,t))$ is a sufficiently smooth solution for \eqref{TMF2} and \eqref{TBC-IC1}. If $u_{D}$ is independent of $t$ and $\Theta_{D}=\Theta_{0}$, then the following energy equality holds:
\begin{align}
	& \frac{d}{dt}\left( \mathcal{E}_{el}(u(t),z(t)) + E_{s} (z(t)) + E_{th}(\Theta(t))\right) \nonumber\\
	& \hspace{2cm}= - \frac{1}{\Theta_{0}} \int_{\Omega} \kappa(z) \left|\nabla\Theta\right|^{2}~dx - \alpha \int_{\Omega} \left| \dot{z}\right| ^{2}dx \leq 0. \label{Eq:Theorem3:c}
\end{align}
\normalfont\begin{Proof}\label{Pr:TFM2} 
\normalfont Since the relation in \eqref{eq-proof1} is written as
\begin{align}
& \frac{d}{dt}\left(\frac{1}{2}W(u)\right)= \sigma^{*}[u,\Theta] - \beta (\Theta-\Theta_{0})\mbox{div}\dot{u} ,\notag
\end{align} 
we obtain
\begin{align}
& \frac{d}{dt}\left(\frac{1}{2}(1-z)^{2}W(u)\right)= \sigma_{z}^{*}[u,\Theta]:e[\dot{u}] + \beta(1-z)^{2} (\Theta-\Theta_{0})\mbox{div}\dot{u} - (1-z)\dot{z}W(u).\notag
\end{align}
Hence, we have
\begin{align}
		&\frac{d}{d t} \mathcal{E}_{el}(u(t),z(t)) + \frac{d}{d t} E_{s}(z(t))\notag \\
		&\quad \quad= \int_{\Omega} \frac{d}{dt}\left(\frac{1}{2}(1-z)^{2}W(u)\right)~dx + \int_{\Omega}\left(\epsilon\mbox{div}(\gamma_{*}\nabla z)-\frac{\gamma_{*}}{\epsilon}z\right)\dot{z}~dx\notag\\
		& \quad \quad= \underbrace{\int_{\Omega}\sigma_{z}^{*}[u,\Theta]:e[\dot{u}]~dx}_{0} + \int_{\Omega}\beta(1-z)^{2} (\Theta-\Theta_{0})\mbox{div}\dot{u}~dx - \int_{\Omega} \mathcal{H}\dot{z}~dx\notag\\
		& \quad \quad= \int_{\Omega}\beta(1-z)^{2}(\Theta-\Theta_{0})\mbox{div}\dot{u}~dx - \int_{\Omega}\alpha|\dot{z}|^{2}~dx,\label{eq:profTFM2a}
\end{align}
where $\mathcal{H} = \epsilon\mbox{div}(\gamma_{*}\nabla z)-\frac{\gamma_{*}}{\epsilon}z + (1-z)W(u)$.

On the other hand,
\begin{align}
& \frac{d}{dt} E_{th}(\Theta(t)) = \frac{\chi}{\Theta_{0}}\int_{\Omega} (\Theta-\Theta_{0})\dot{\Theta}~dx\notag\\
& \quad \quad \quad= \frac{1}{\Theta_{0}}\int_{\Omega} (\Theta-\Theta_{0})\left\{\mbox{div}(\kappa(z)\nabla \Theta) - \Theta_{0}\beta(1-z)^{2}\mbox{div}\dot{u}\right\}~dx\notag\\
&\quad \quad \quad= -\frac{1}{\Theta_{0}}\int_{\Omega} \kappa(z)|\nabla \Theta|^{2}~dx - \int_{\Omega} \beta(1-z)^{2}(\Theta-\Theta_{0})\mbox{div}\dot{u}~dx.\label{eq:profTFM2b}
\end{align}
Taking a sum of these equalities \eqref{eq:profTFM2a}-\eqref{eq:profTFM2b}, we obtain the energy equality \eqref{Eq:Theorem3:c}. \qed
\end{Proof}
\end{Theorem}

\section{Numerical Experiments} \label{sec-4}
In this section, we conduct numerical experiments to test F-PFM, TF-PFM1, and TF-PFM2, which were derived in Section \ref{sec:3}, and report the numerical results. Through the numerical experiments, we observe the effect of thermal coupling on the crack speed and the crack path during its growth process.

\subsection{Nondimensional setting}
In the following numerical examples, we suppose $\kappa(z) = (1-z)^{2}\kappa_{0}$. For convenience, we consider the nondimensional form with \eqref{ND}, \eqref{ND2}, \eqref{Non-Scaling},  and
\begin{eqnarray}
\tilde{\epsilon} = \frac{\epsilon}{c_{x}},  ~ \tilde{\gamma}_{*} = \frac{{c_{e}\gamma_{*}}}{c_{x}(\beta c_{\Theta})^{2}}, ~\tilde{\alpha} =\frac{c_{e} \alpha}{c_{t}(\beta c_{\Theta})^{2}}, ~\tilde{a}_{L} = \frac{c_{x}c_{\Theta}}{c_{u}}{a}_{L}, ~\tilde{\beta} = 1.\notag
\end{eqnarray}
Then, TF-PFM1 in \eqref{TMF1} is expressed in the following nondimensional form:
\begin{subequations}
	\begin{empheq}[left=\empheqlbrace]{align}
	&-{\mbox{div}}\left((1-z)^{2}{\sigma}[{u}]\right)	= (1-z)^{2}\nabla\Theta  &\hspace{-16pt}\mbox{in}~ {\Omega}\times[0,{T}],  \label{TMF12aND}\\
	&{\alpha} \frac{\partial z}{\partial {t}} = \left({\epsilon}~{\mbox{div}}({\gamma_{*}} {\nabla} z) - \frac{{\gamma_{*}}}{{\epsilon}}z + (1-z){W}^{*}(u,\Theta)\right)_{+} &\hspace{-16pt}\mbox{in}~ {\Omega}\times[0,{T}], \label{TMF12bND}\\
	&\frac{\partial {\Theta}}{\partial {t}} = ~{\mbox{div}}\left((1-z)^{2}{\nabla} {\Theta}\right) - (1-z)^{2}\delta\frac{\partial}{\partial {t}}({\mbox{div}} {u}) &\hspace{-16pt}\mbox{in}~ {\Omega}\times(0,{T}]. \label{TMF12dND}
	\end{empheq}\label{TMF12ND}
\end{subequations} \\[-3pt]
For TF-PFM2, we change \eqref{TMF12bND} to:
\begin{align}
&{\alpha} \frac{\partial z}{\partial {t}} = \left({\epsilon}~{\mbox{div}}({\gamma_{*}} {\nabla} z) - \frac{{\gamma_{*}}}{{\epsilon}}z + (1-z){W}(u)\right)_{+} & \hspace{-16pt}\mbox{in}~ {\Omega}\times[0,{T}]. \label{TMF12cND}
\end{align}

\subsection{Time discretization}\label{SubSec4.2}
\begin{figure}[!h]
	\begin{tabular}{cc}
		\includegraphics[trim=1cm 5.5cm 1cm 4.5cm,width=0.4\textwidth]{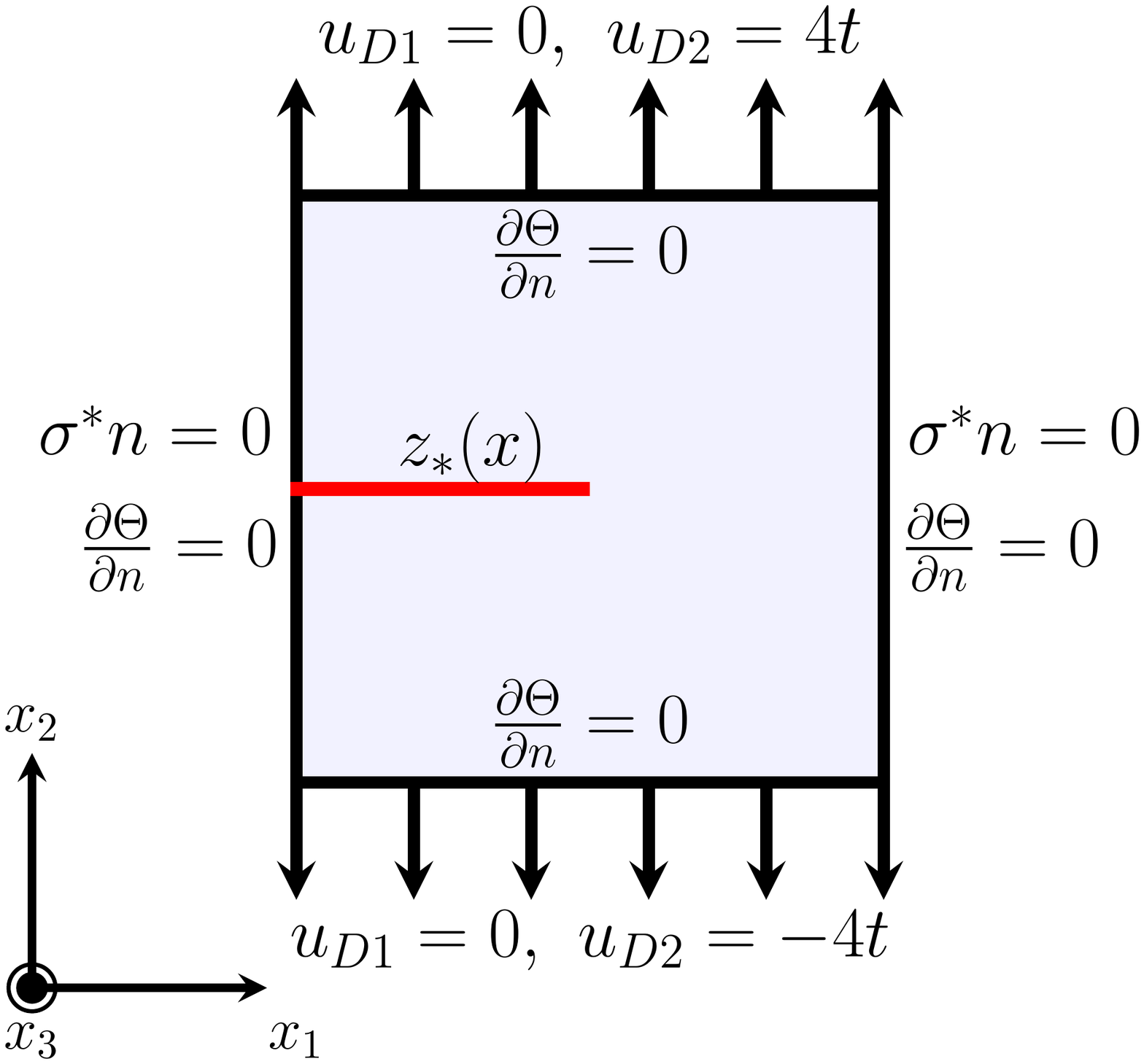}
		&
		{\includegraphics[trim=0cm -5cm 1cm 0cm,width=0.25\textwidth]{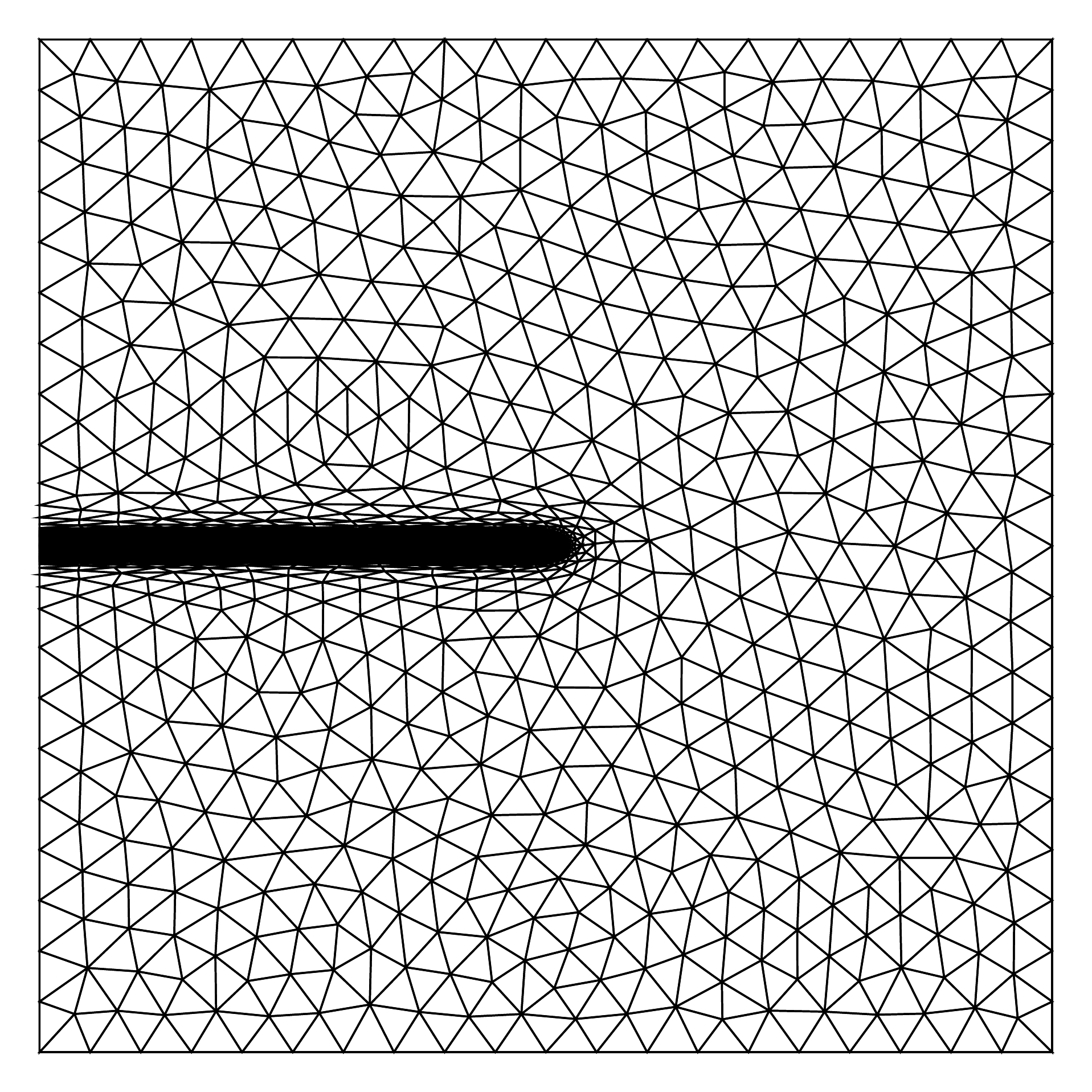}}
	\end{tabular}
	\caption{Domain for Section \ref{SubSec4.3} with $z_{*}(x)$ as the initial crack (left) and the adaptive mesh for the initial crack (right).} 
	\label{dom:2}
\end{figure}

To solve problem \eqref{TMF12ND}, we adopt the following semi-implicit time discretization scheme \cite{Kimura2021,Kimura2009}.
\begin{subequations}
	\begin{empheq}[left=\empheqlbrace]{align}
	& -{\mbox{div}}\left((1-z^{k-1})^{2}{\sigma}^{*}[{u}^{k},{\Theta}^{k-1}]\right)	= 0, \label{Disk1}\\
	& {\alpha} \frac{\tilde{z}^{k} - z^{k-1}}{\Delta{t}} = {\epsilon}~{\mbox{div}}\left({\gamma_{*}} {\nabla} \tilde{z}^{k}\right) - \frac{{\gamma_{*}}}{{\epsilon}}\tilde{z}^{k} + \left(1-\tilde{z}^{k}\right){W}^{*}({u}^{k-1},{\Theta}^{k-1}),\label{Disk2}\\
	& z^{k} := \mbox{max}\left(\tilde{z}^{k}, z^{k-1}\right), \label{Disk5}\\
	& \dfrac{{\Theta}^{k} - {\Theta}^{k-1}}{\Delta{t}} = ~{\mbox{div}}\left((1- z^{k-1}){\nabla} {\Theta}^{k}\right) - (1-z^{k-1})\delta{\mbox{div}}\Big(\dfrac{ {u}^{k} - {u}^{k-1}}{\Delta {t}}\Big).\label{Disk4}
	\end{empheq}\label{Disk}
\end{subequations}\\[-3pt]
For TF-PFM2, \eqref{Disk2} is replaced by
\begin{align}
& {\alpha} \frac{\tilde{z}^{k} - z^{k-1}}{\Delta{t}} = {\epsilon}~{\mbox{div}}\left({\gamma_{*}} {\nabla} \tilde{z}^{k}\right) - \frac{{\gamma_{*}}}{{\epsilon}}\tilde{z}^{k} + \left(1-\tilde{z}^{k}\right){W}({u}^{k-1}), \label{Disk3}
\end{align}
where  $u^{k}$, $z^{k}$, and $\Theta^{k}$ are the approximations of $u$, $z$,  $\Theta$, respectively, at time ${t}_{k} := k\Delta t (k = 1, 2, 3,\cdots)$. Since the adaptive mesh technique in the FEM is often effective and accurate in numerical experiments with phase field models,  problems \eqref{Disk} - \eqref{Disk3} are calculated using adaptive finite elements with P2 elements with a minimum mesh size of $h_{min} = 2\times 10^{-3}$ and a maximum mesh size of $h_{max} = 0.1$. The adaptive mesh control at each time step is performed by the adaptmesh() command in FreeFEM based on the variable $z$. An example of the adaptive mesh is illustrated in Figure \ref{dom:2} (right).  In addition, the code for the following numerical experiments in the current study is written on FreeFEM \cite{Hecht2012} and executed on a desktop  with an Intel(R) Core i7$-$7820X CPU@3.60 GHz, 16 core processor, and 64 GB RAM. 

\subsection{Thermoelastic effect on the crack speed} \label{SubSec4.3}
We set a square domain $\Omega := (-1,1)^{2} \subset \mathbb{R}^{2}$ with the initial crack $z_{*}(x) := \exp{(-(x_{2}/\eta)^{2})}/(1 + \exp{(x_{1}/\eta)})$ and $\eta = 1.5\times 10^{-2}$. The initial mesh is adapted to $z_{*}(x)$, as illustrated in Figure \ref{dom:2} (right). The material constants for the following examples in the nondimensional form are listed in Table \ref{tab:4}.
\begin{table}[!h]
	\caption{List of the nondimensional parameters for Sections \ref{SubSec4.3} and \ref{SubSec4.4}}
	\label{tab:4}       
	\begin{tabular}{p{1.65cm}m{1.0cm}p{1.0cm}m{1.0cm}p{1.0cm}m{1.0cm}p{1.0cm}p{1.0cm}}
		\hline\noalign{\smallskip}
		Parameter & $E_{Y}$ & $\nu_{P}$ & $a_{L}$ & $\alpha$ & $\epsilon$ & $\gamma_{*}$ & $\Theta_{*}$  \\
		\noalign{\smallskip}\hline\noalign{\smallskip}
		Value & 1 & 0.3 & 0.7 & 0.001 & 0.01 & 5.08& 0\\
		\noalign{\smallskip}\hline
	\end{tabular}
\end{table}

The boundary conditions for $u$ and $\Theta$ are illustrated in Figure \ref{dom:2} (left). For $z$, we set $\frac{\partial  z}{\partial n} = 0$ on $\Gamma$.

\begin{figure}[!h]
	\begin{center}
		\begin{tabular}{c}
			{\includegraphics[trim=1cm 6cm 1cm 6cm, scale=0.46]{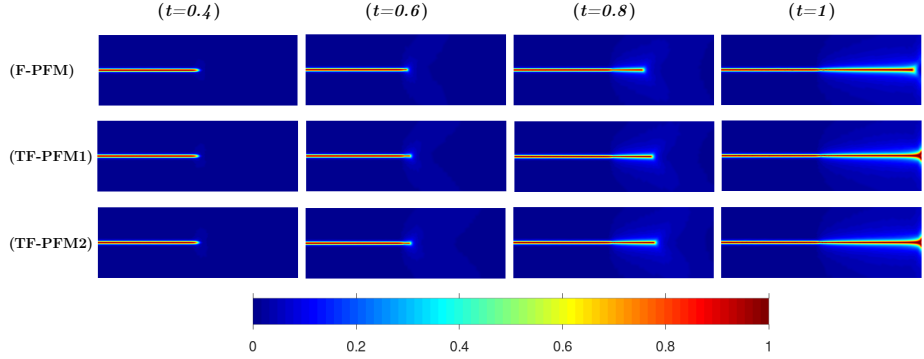}}
		\end{tabular}
		\caption{Snapshots of crack propagation with F-PFM, TF-PFM1, and TF-PFM2 in  $(-1,1)\times (-0.35,0.35)$  at $t = 0.4, ~0.6, ~0.8, ~1$ (left to right). For TF-PFM1 and TF-PFM2, we use the thermoelasticity coupling parameter $\delta = 0.5$, and the color represents the value of $z$.}\label{fig1}
	\end{center}
\end{figure}

\begin{figure}[!h]
	\begin{tabular}{cc}
		\includegraphics[trim=2cm 6cm 2cm 6cm, scale=0.28]{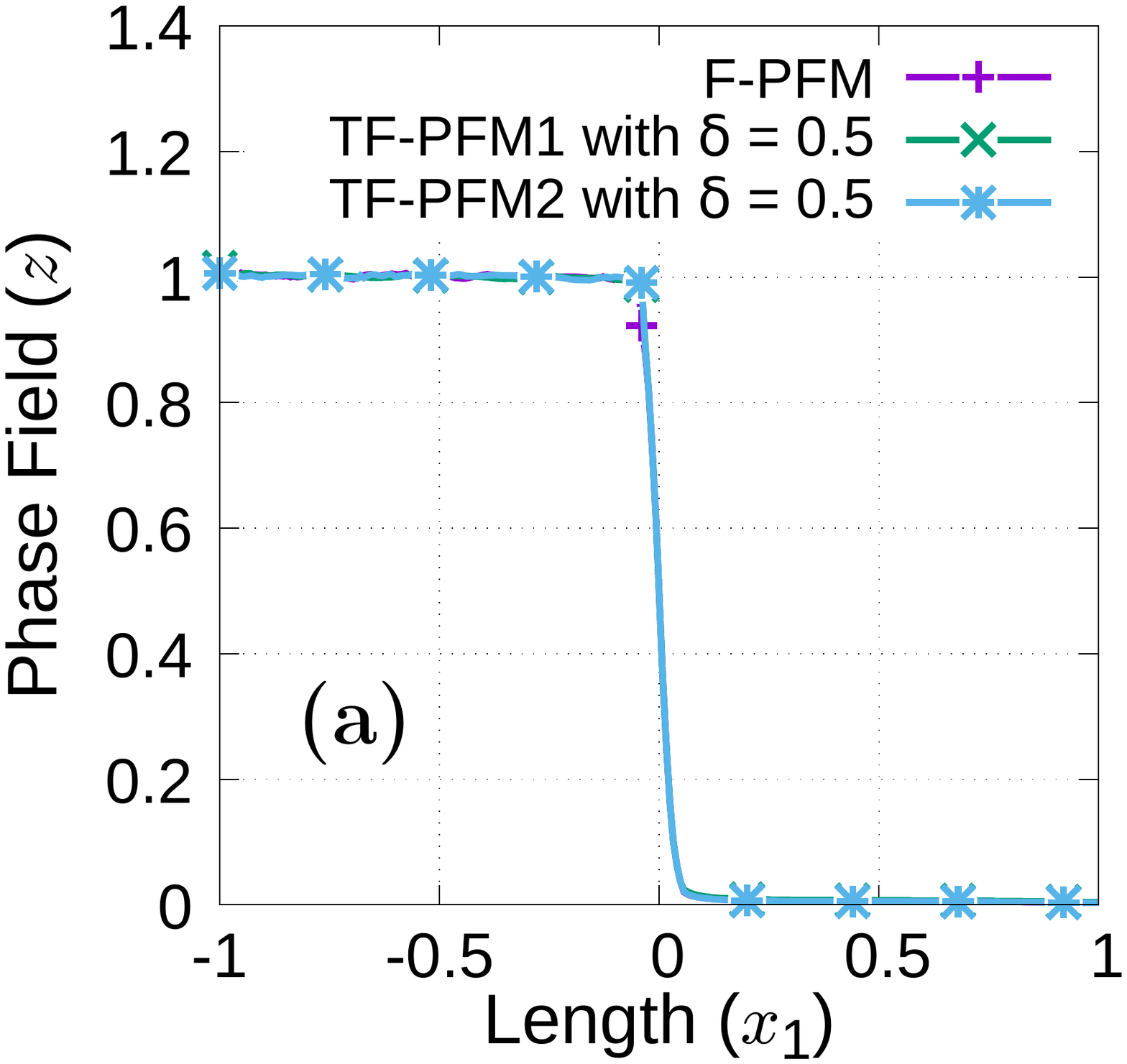}
		&
		\includegraphics[trim=2cm 6cm 2cm 6cm, scale=0.28]{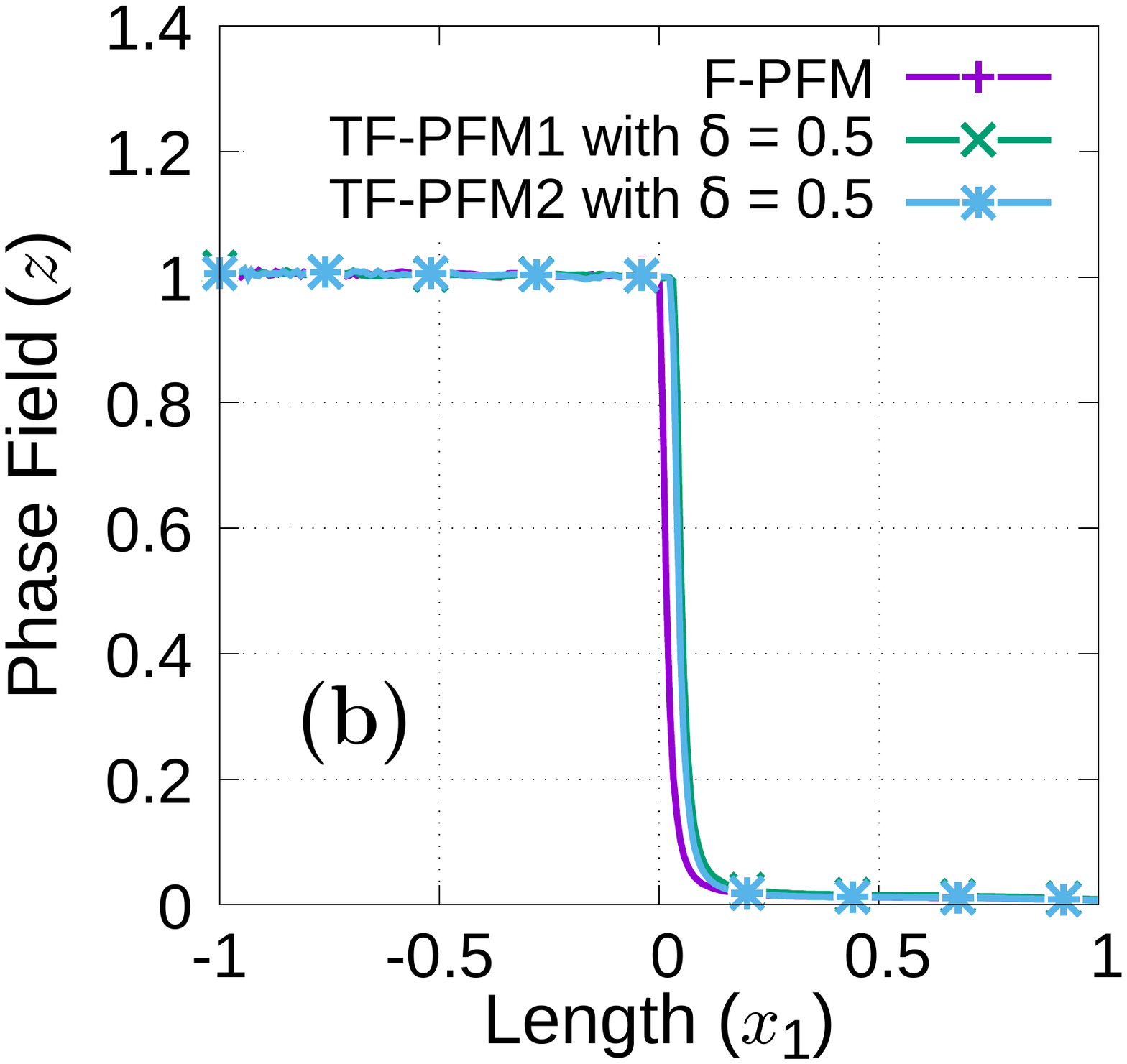}
		\\[5pt]
		\includegraphics[trim=2cm 6cm 2cm 6cm, scale=0.28]{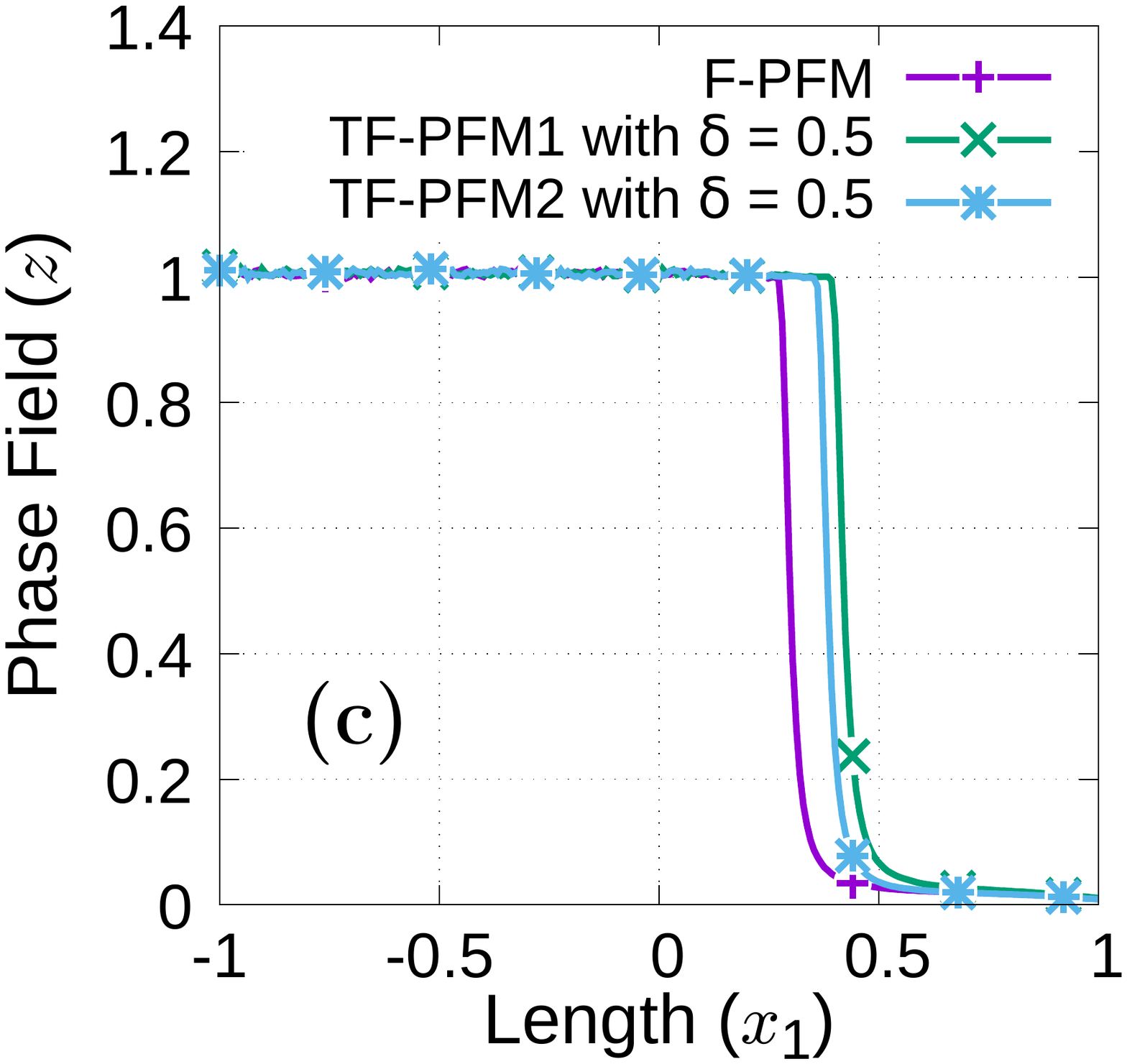}
		&
		\includegraphics[trim=2cm 6cm 2cm 6cm, scale=0.28]{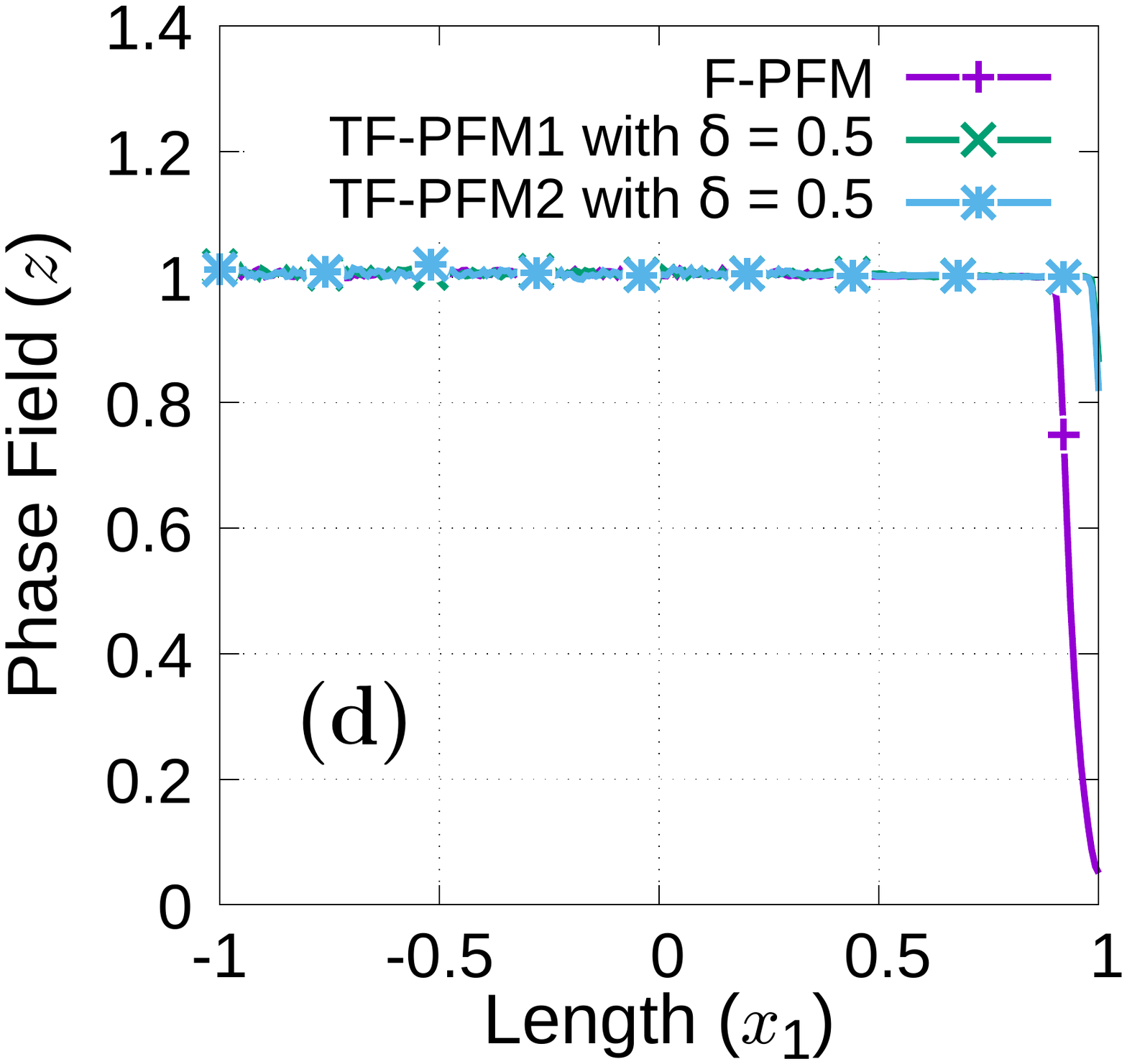}
	\end{tabular}
	\caption{Comparison of the profiles of $z$ obtained by F-PFM, TF-PFM1, and TF-PFM2 along the line $x_{2}=0$ at (a) $t=0.4$, (b) $t=0.6$, (c) $t=0.8$, and (d) $t=1$.}\label{Compare01}
\end{figure}

In Figure \ref{fig1}, the numerical results obtained by F-PFM, TF-PFM1, and TF-PFM2 are shown in the upper, middle, and bottom parts, respectively, where we set $\delta = 0.5$ for TF-PFM1 and TF-PFM2. In addition, the profile of $z$ on line $x_{2}=0$ is shown in Figure \ref{Compare01}. From Figures \ref{fig1} and \ref{Compare01}, we observe that the crack propagation rate obtained by F-PFM is slower than that obtained by the others, and that the crack propagation rate obtained by TF-PFM1 is slightly faster than that obtained by TF-PFM.  

\begin{figure}[!h]
	\begin{center}
		\begin{tabular}{cccc}
			{\includegraphics[trim=2.5cm 0cm 1cm 0cm, scale=0.145]{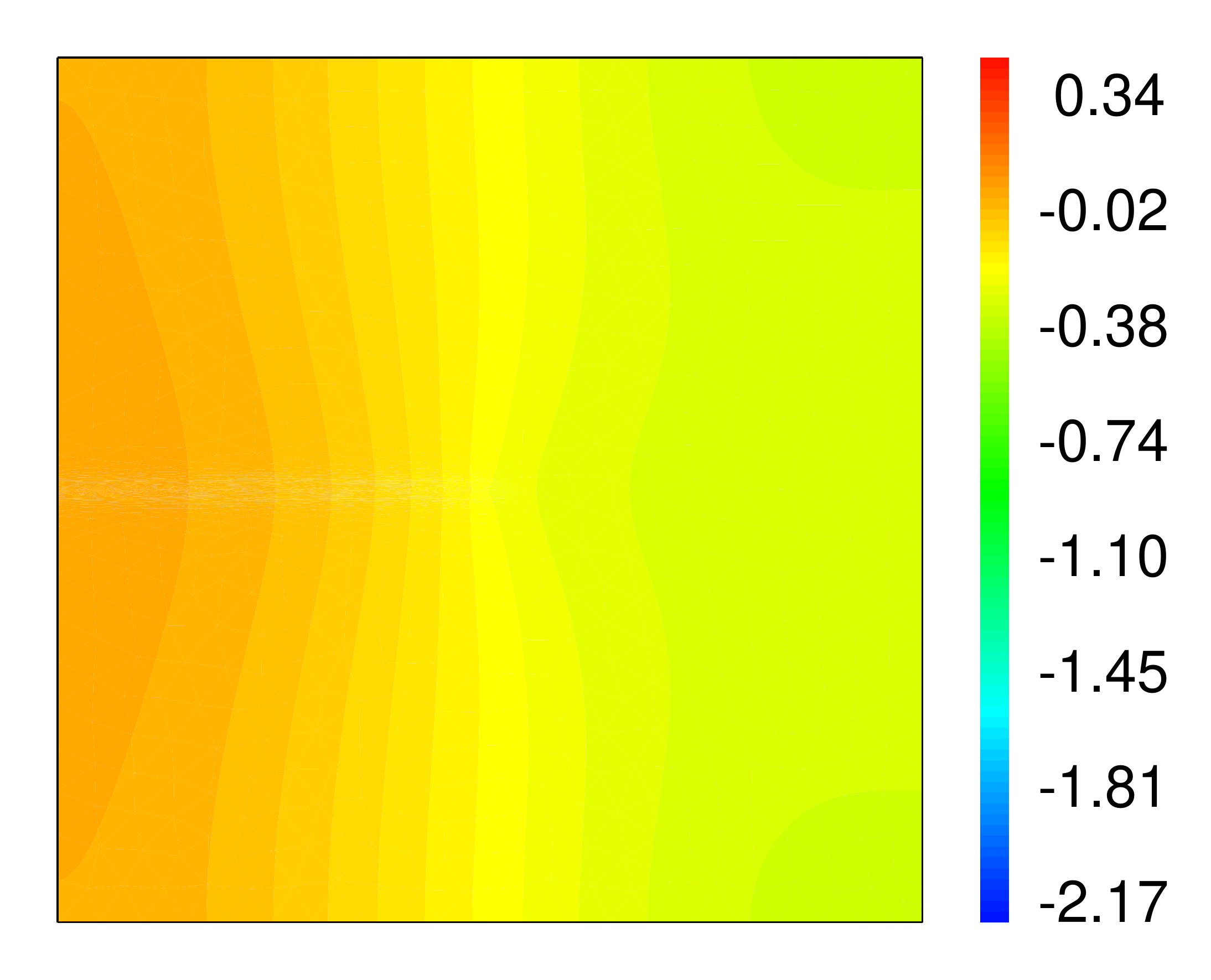}} & \hspace{-10pt}
			{\includegraphics[trim=1.5cm 0cm 1cm 0cm, scale=0.145]{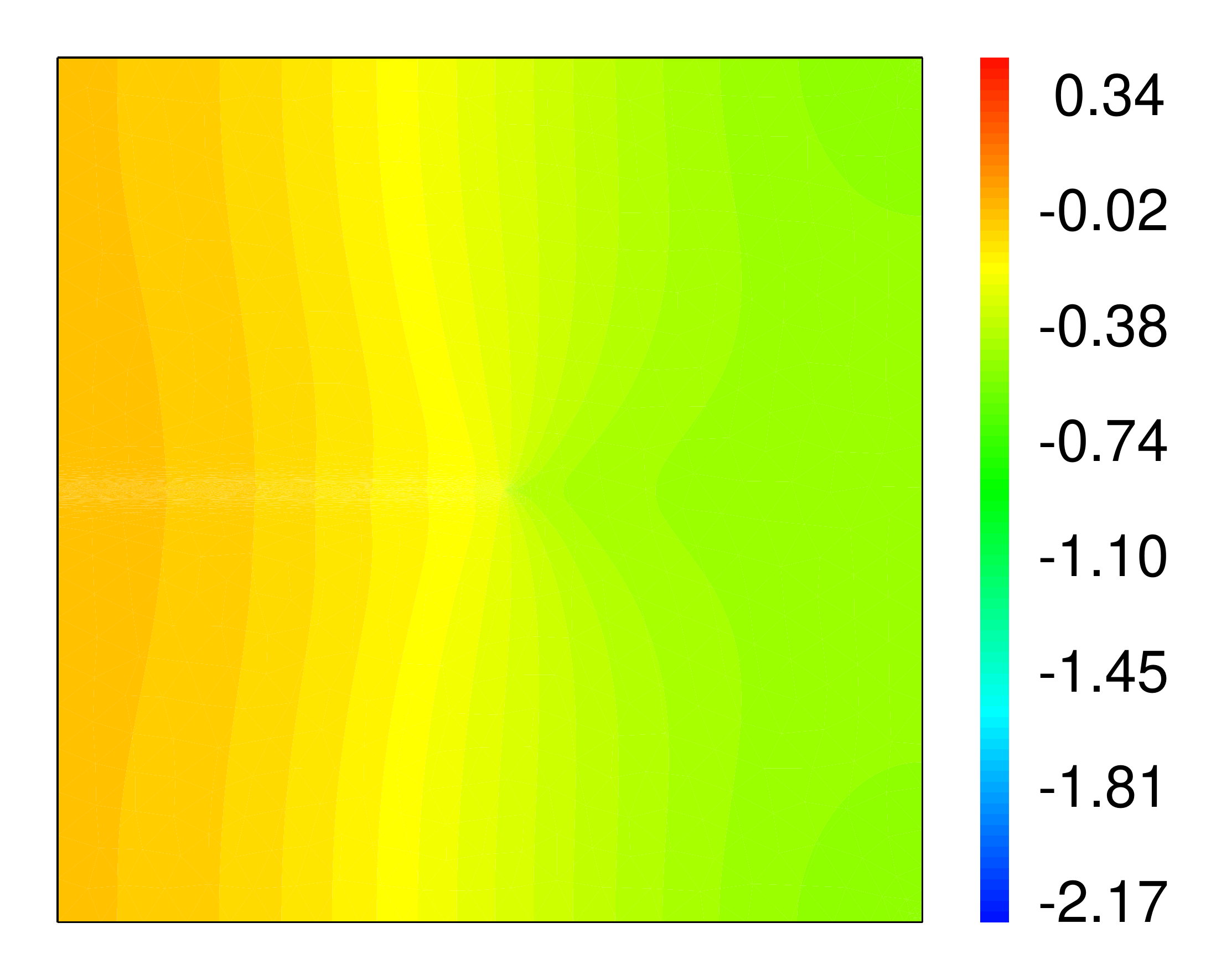}} & \hspace{-10pt}
			{\includegraphics[trim=1.5cm 0cm 1cm 0cm, scale=0.145]{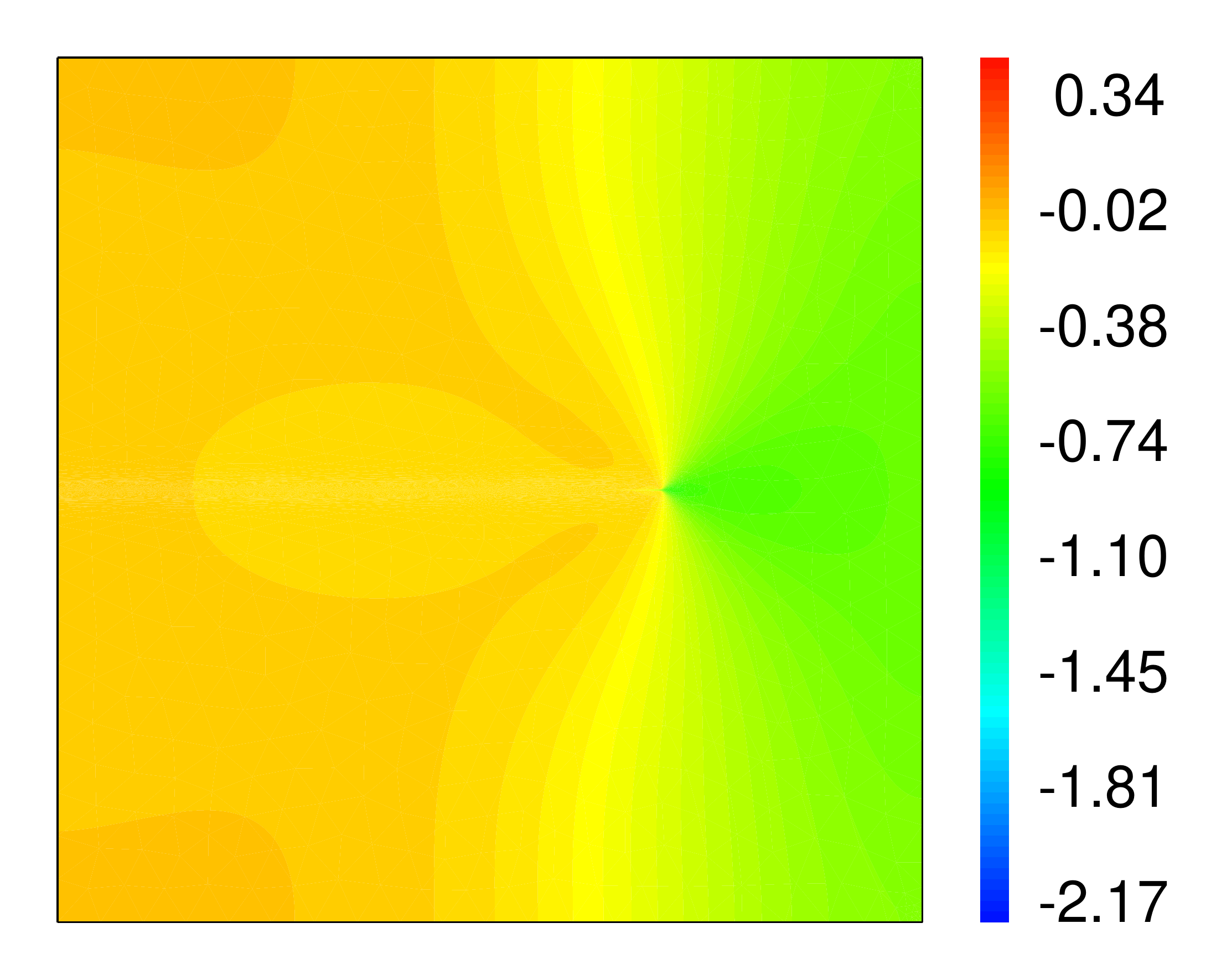}} & \hspace{-10pt}
			{\includegraphics[trim=1.5cm 0cm 1cm 0cm, scale=0.145]{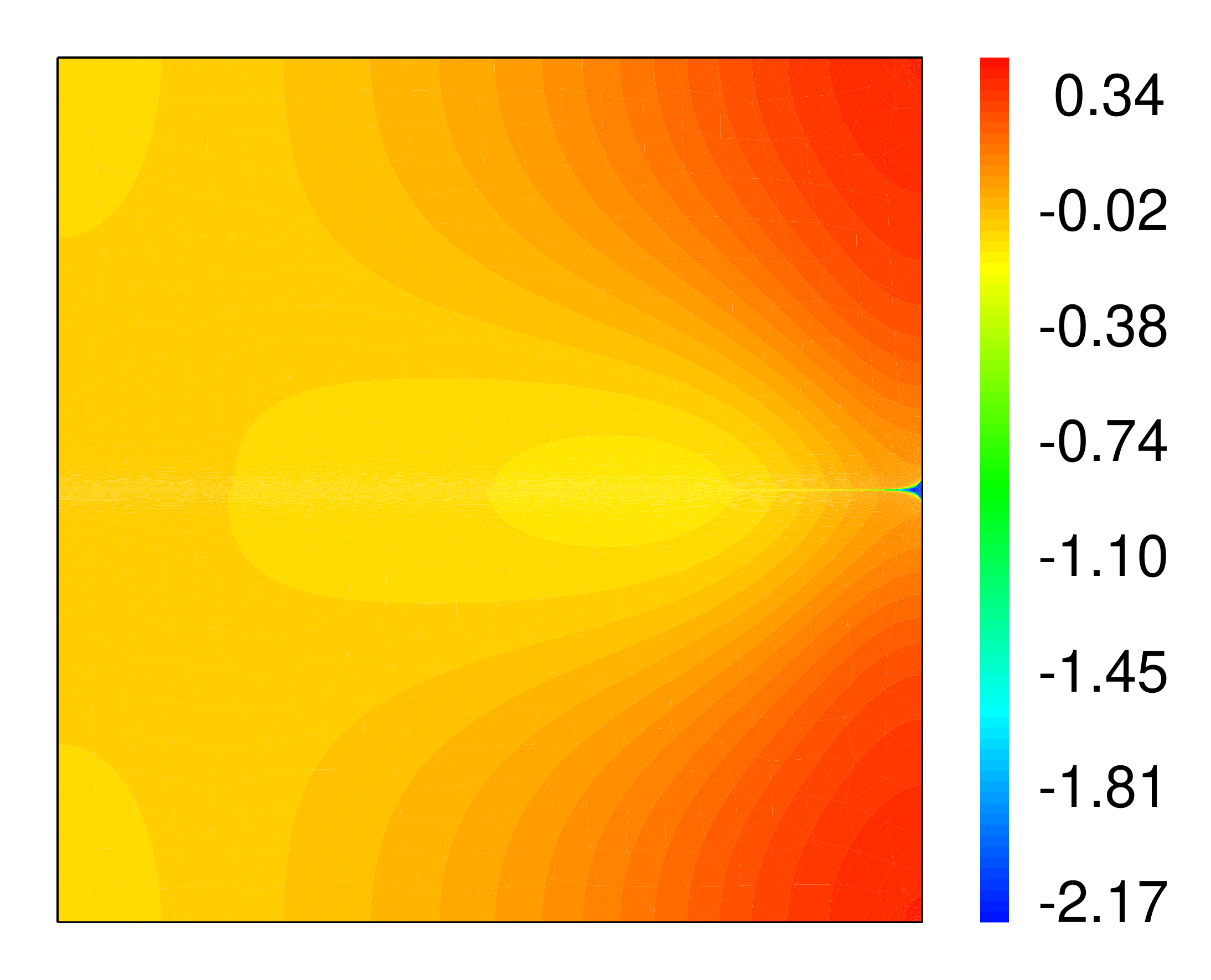}} \\[-8pt]
			{\includegraphics[trim=2.5cm 0cm 1cm 0cm, scale=0.145]{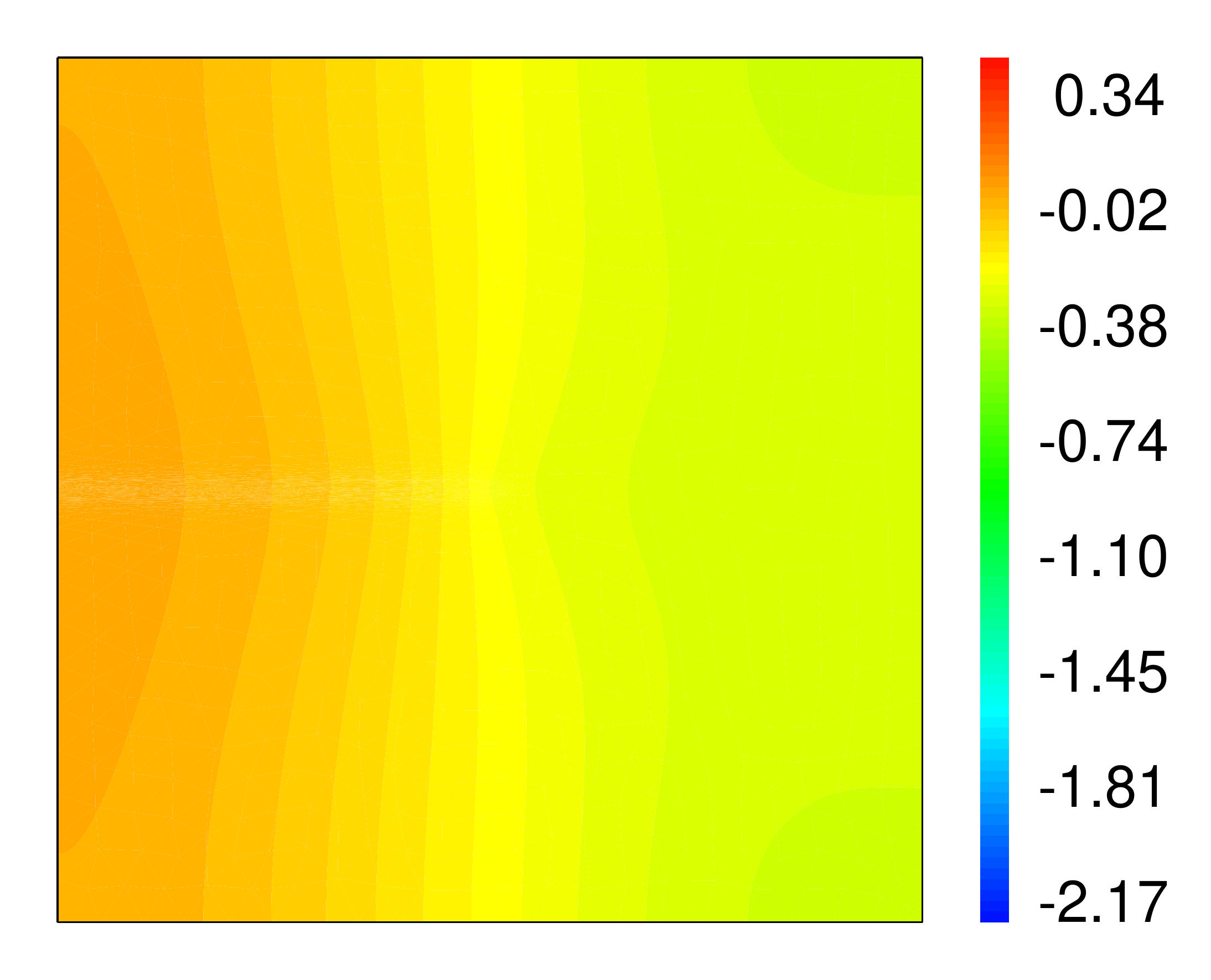}} & \hspace{-10pt}
			{\includegraphics[trim=1.5cm 0cm 1cm 0cm, scale=0.145]{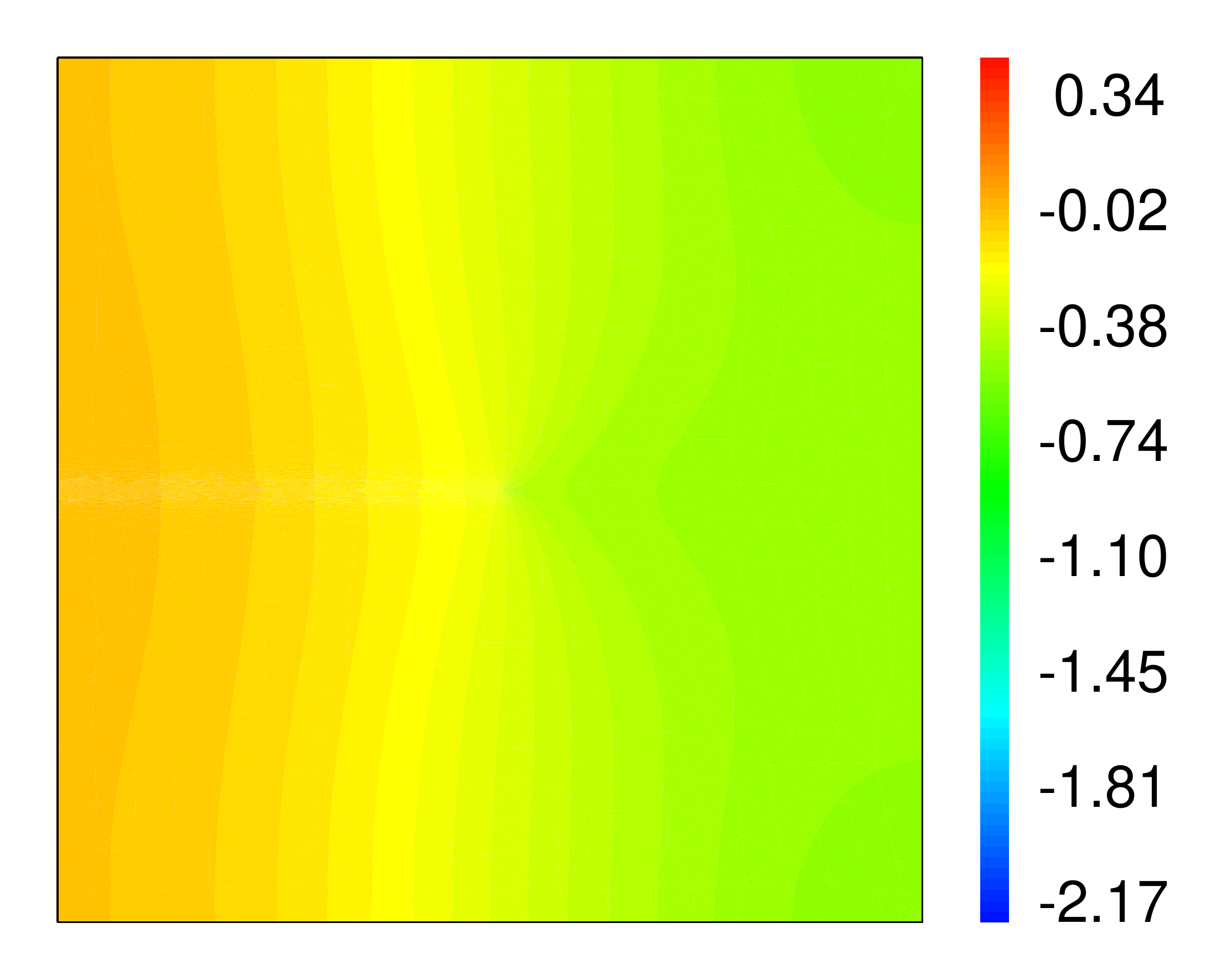}} & \hspace{-10pt}
			{\includegraphics[trim=1.5cm 0cm 1cm 0cm, scale=0.145]{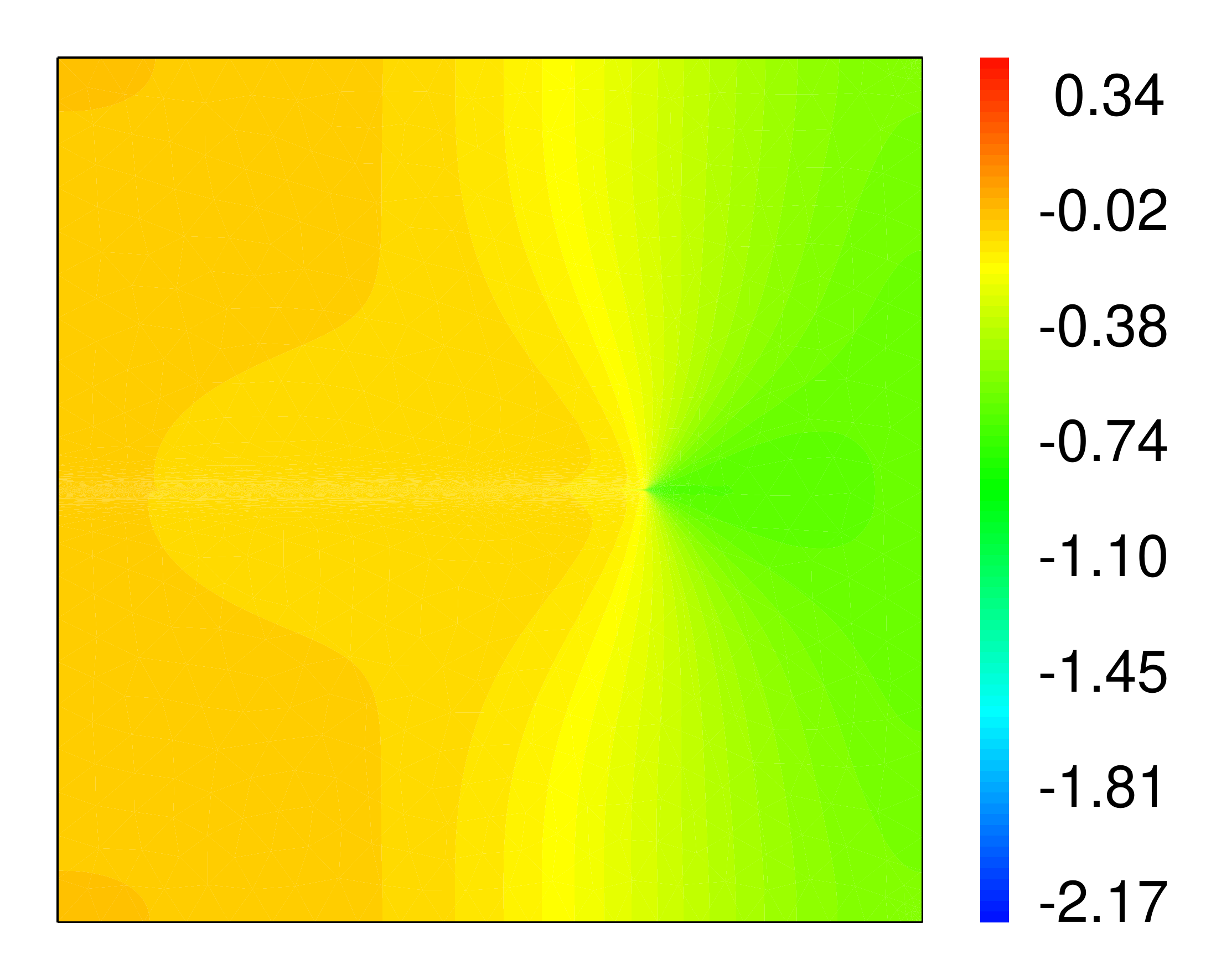}} & \hspace{-10pt}
			{\includegraphics[trim=1.5cm 0cm 1cm 0cm, scale=0.145]{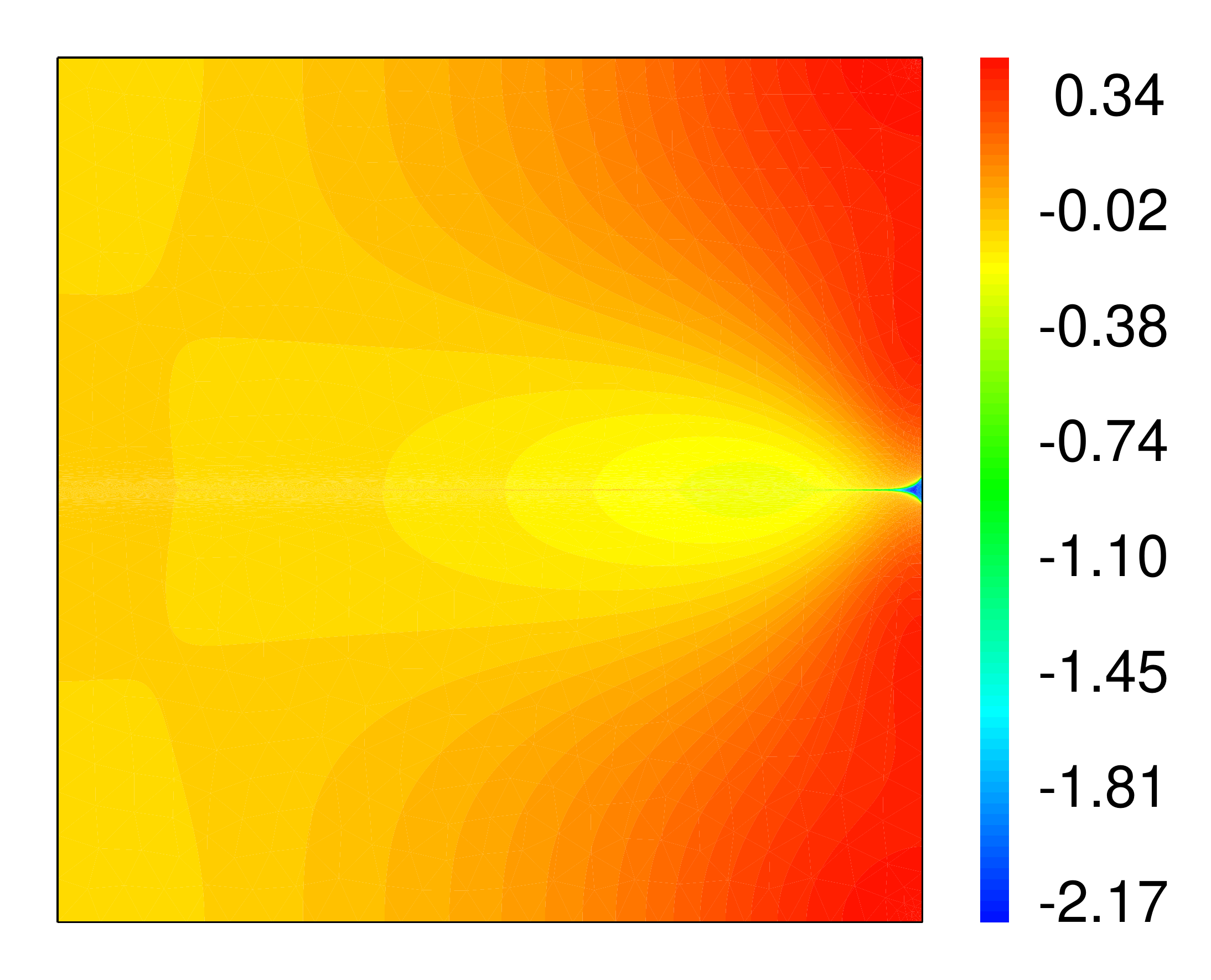}}
		\end{tabular}\\[-10pt]
		\caption{Snapshots of the temperatures obtained by TF-PFM1 (upper) and TF-PFM2 (lower) at $t = 0.4, ~0.6, ~0.8, ~1$ (left to right); the color represents the value of $\Theta$.}\label{fig2}
	\end{center}
\end{figure}

The temperature distributions obtained by TF-PFM1 and TF-PFM2 are shown in Figure \ref{fig2}. In the equation for $\Theta$, the heat resource is given by $-(1-z)^{2}\delta\frac{d}{dt}(\mbox{div}u)$. During crack propagation $(0.4 \leq t \leq 0.8)$, the areas near the crack tip, the upper-right corner, and lower-right corner are continuously expanding when $\mbox{div}u > 0$ and $\frac{\partial}{\partial t}(\mbox{div}u) > 0$. Therefore, due to the negative source $-\frac{\partial}{\partial t}(\mbox{div}u)$, lower temperatures are observed in those areas. On the other hand, at $t=1$, due to the sudden compression caused by the total fracture, positive heat is generated, and a higher temperature is observed, especially near the upper-right and lower-right corners. 


\begin{figure}[!h]
	\begin{tabular}{cc}
		\hspace{-5pt}\includegraphics[trim=1cm 0cm 1cm 0cm,width=0.5\linewidth]{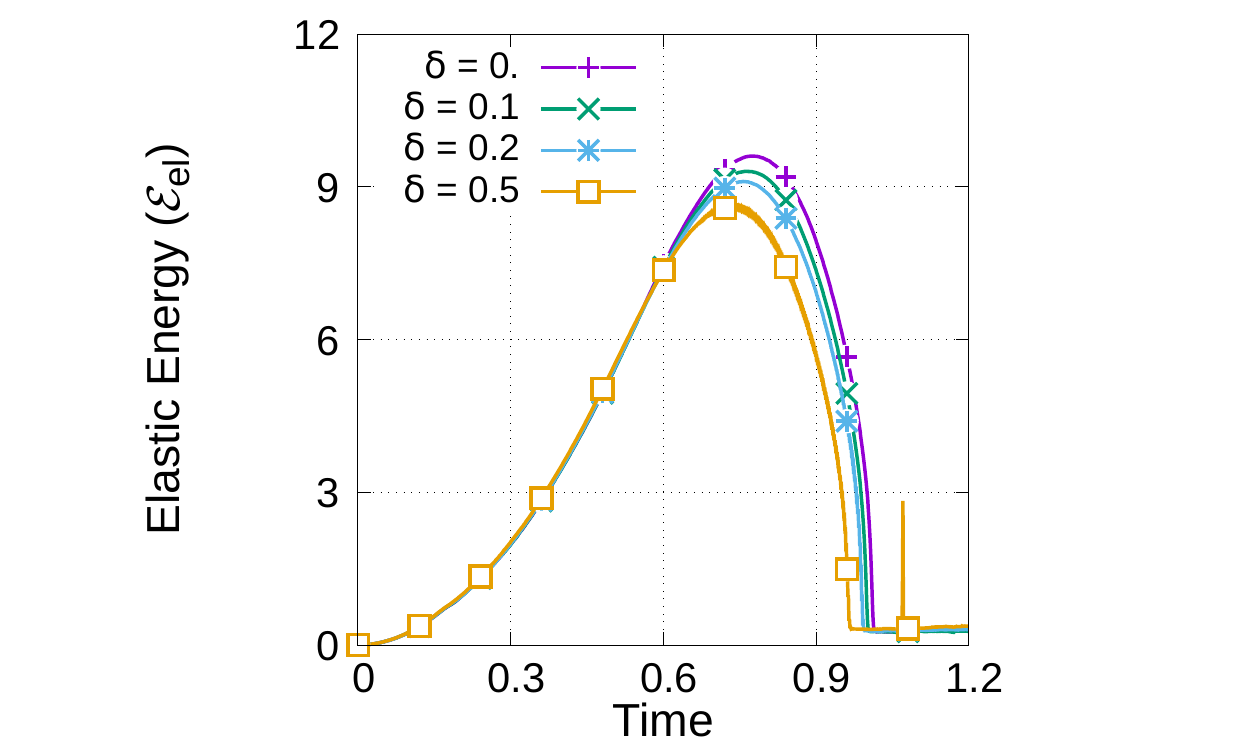}&\hspace{-1cm}
		\includegraphics[trim=1cm 0cm 1cm 0cm,width=0.5\linewidth]{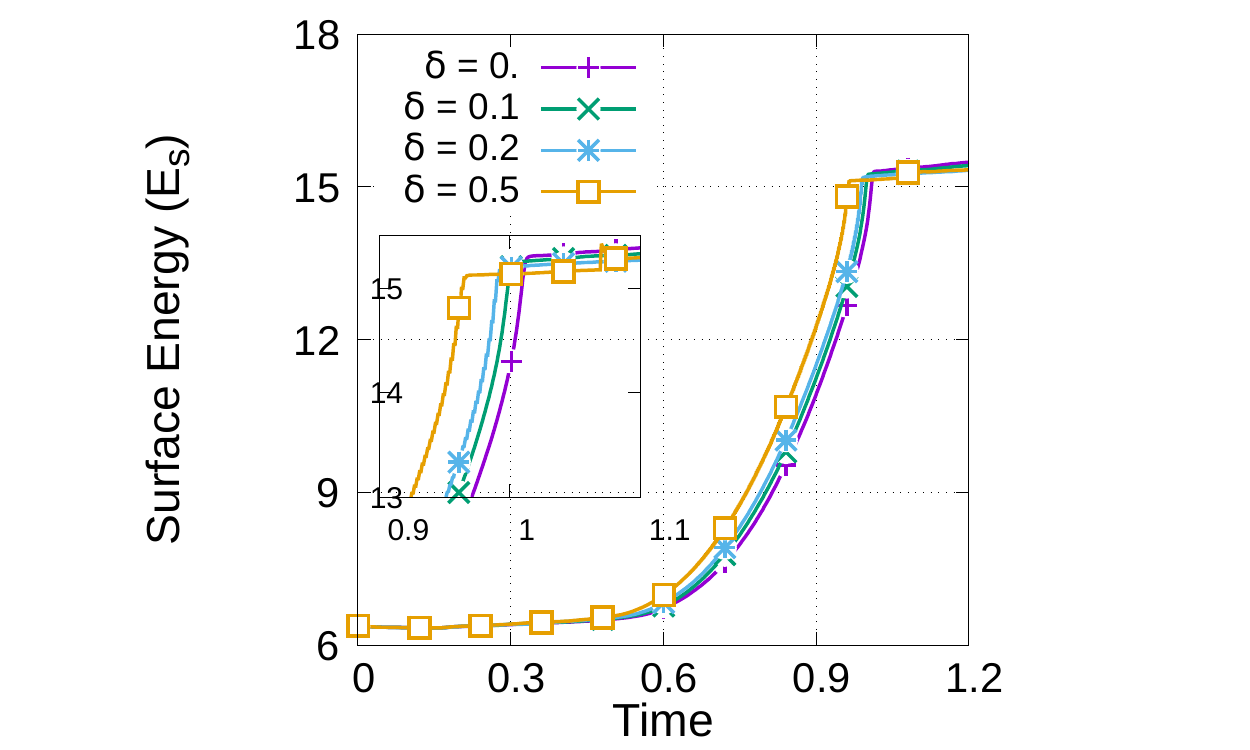}\\
		\hspace{-5pt}\includegraphics[trim=1cm 0.5cm 1cm 0cm,width=0.5\linewidth]{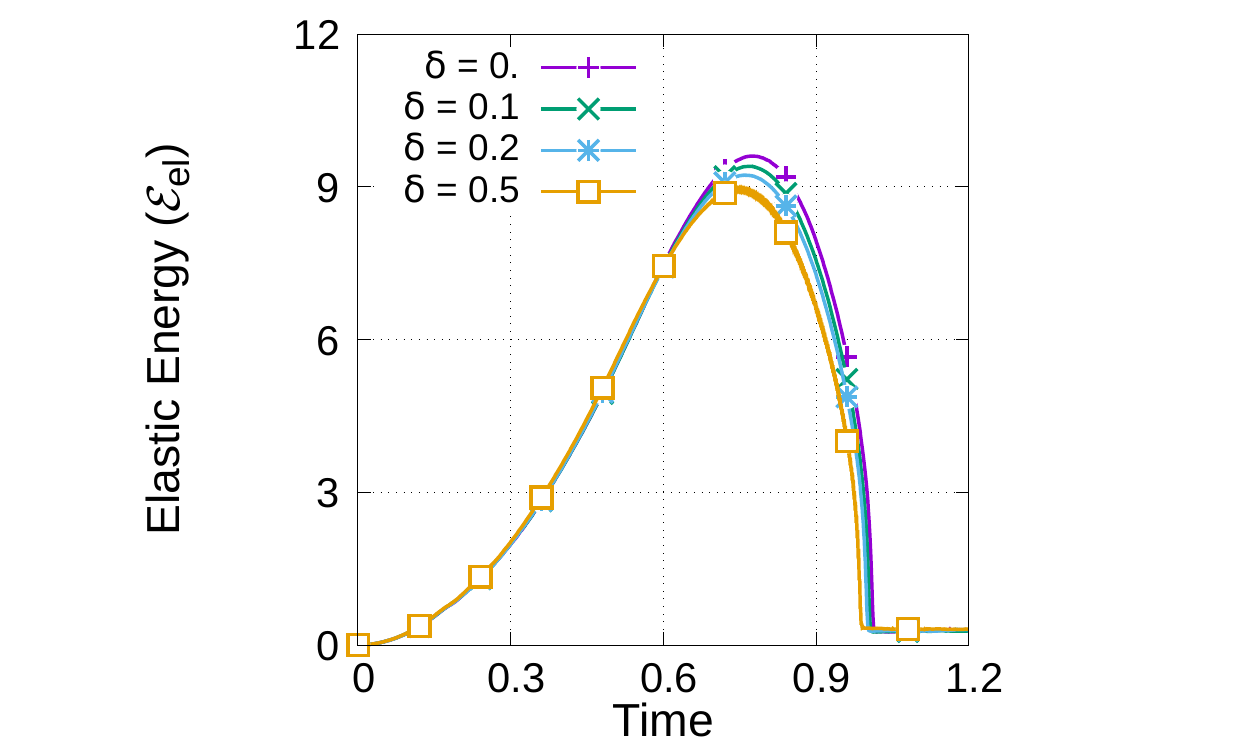}&\hspace{-1cm}
		\includegraphics[trim=1cm 0.5cm 1cm 0cm,width=0.5\linewidth]{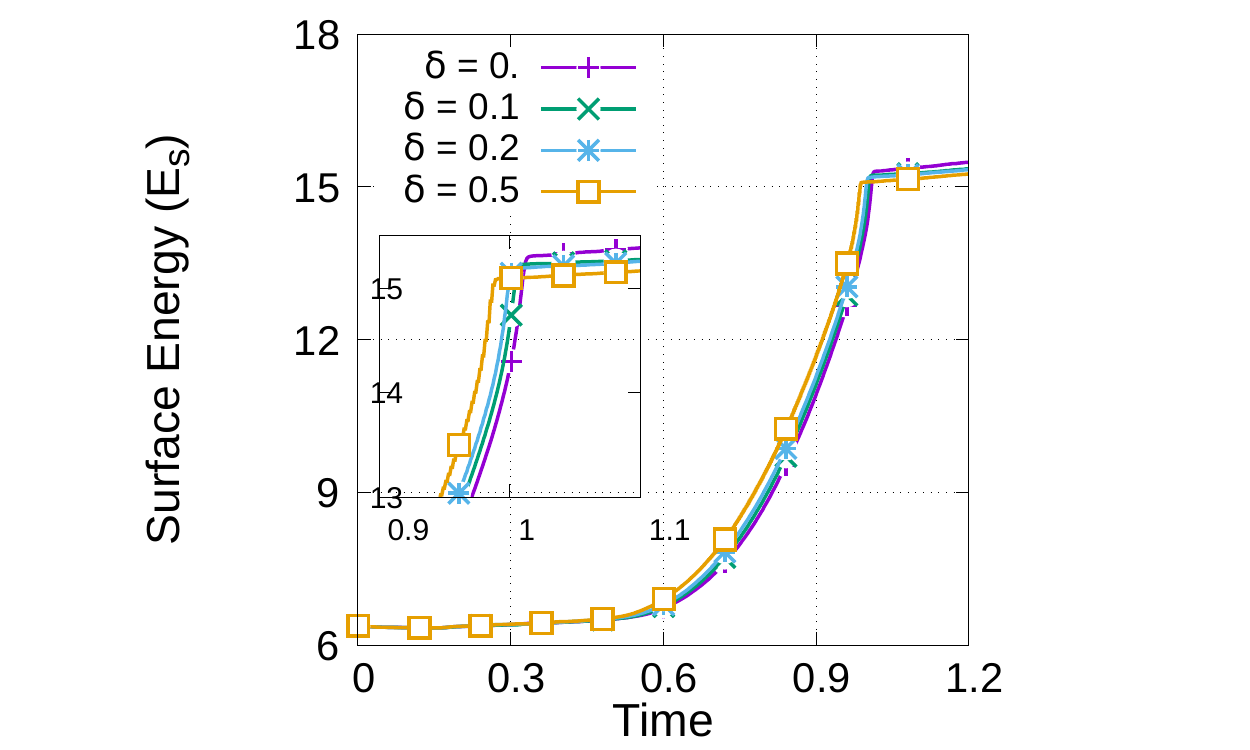}
	\end{tabular}
	\caption{Profile of the elastic (left) and surface energy  (right) under thermal expansion during crack propagation using TF-PFM1 (top) and TF-PFM2 (bottom).}\label{fig8}
\end{figure}

To see how the thermoelastic coupling parameters contribute to enhanced crack propagation, we consider $\delta = 0, ~0.1, ~0.2, ~0.5$ for TF-PFM1 and TF-PFM2, and their elastic and surface energies are plotted in Figure \ref{fig8}.  From Figure \ref{fig8}, we observe that faster crack propagation occurs with a larger coupling parameter. The figure also shows that crack propagation using TF-PFM1 is faster than that using TF-PFM2.

\begin{figure}[!h]
	\begin{tabular}{cc}
		\includegraphics[trim=1cm 7.25cm 1cm 9cm,width=0.435\textwidth]{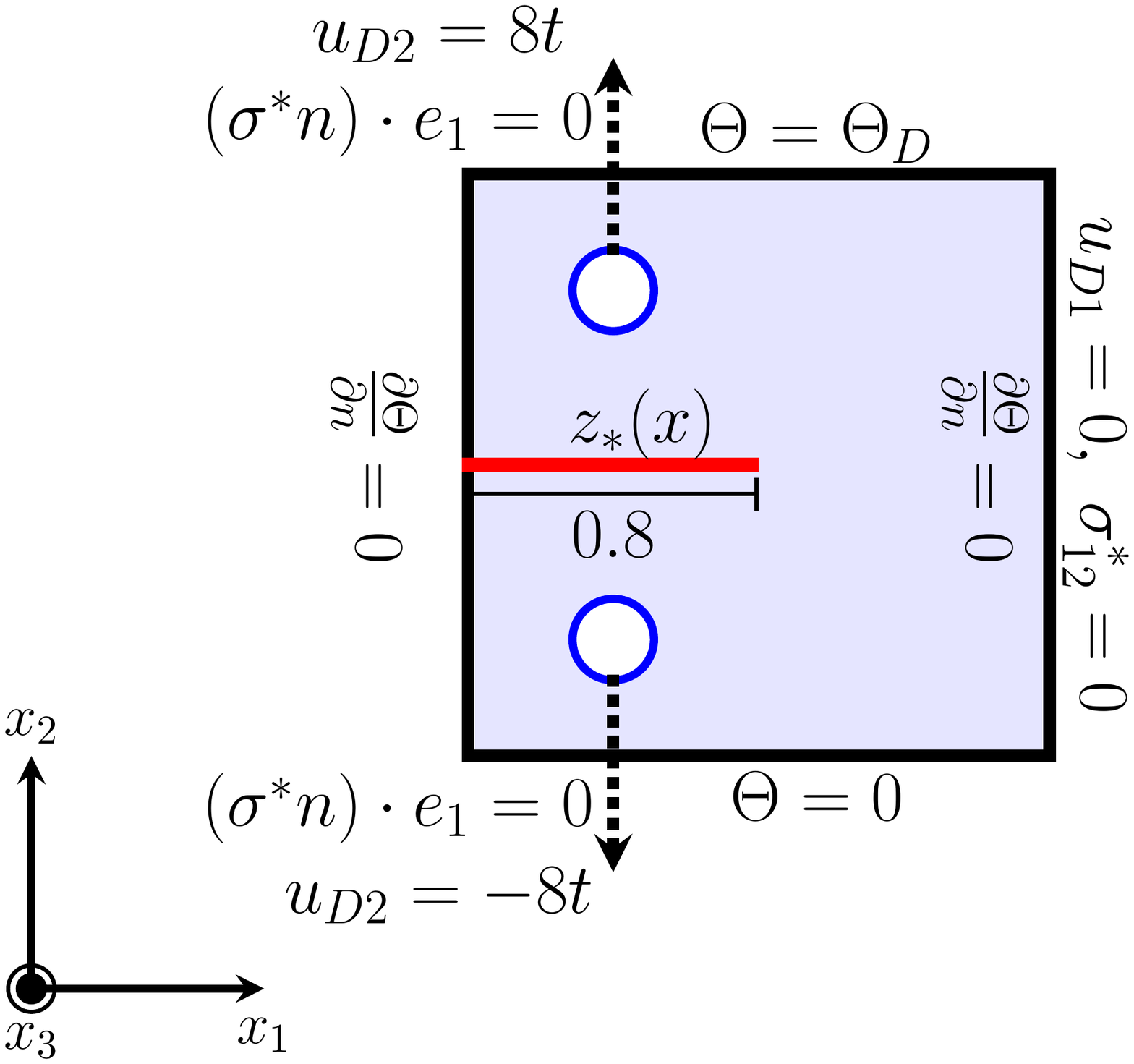}
		&
		{\includegraphics[trim=0cm 5.75cm 1cm 6cm,width=0.45\textwidth]{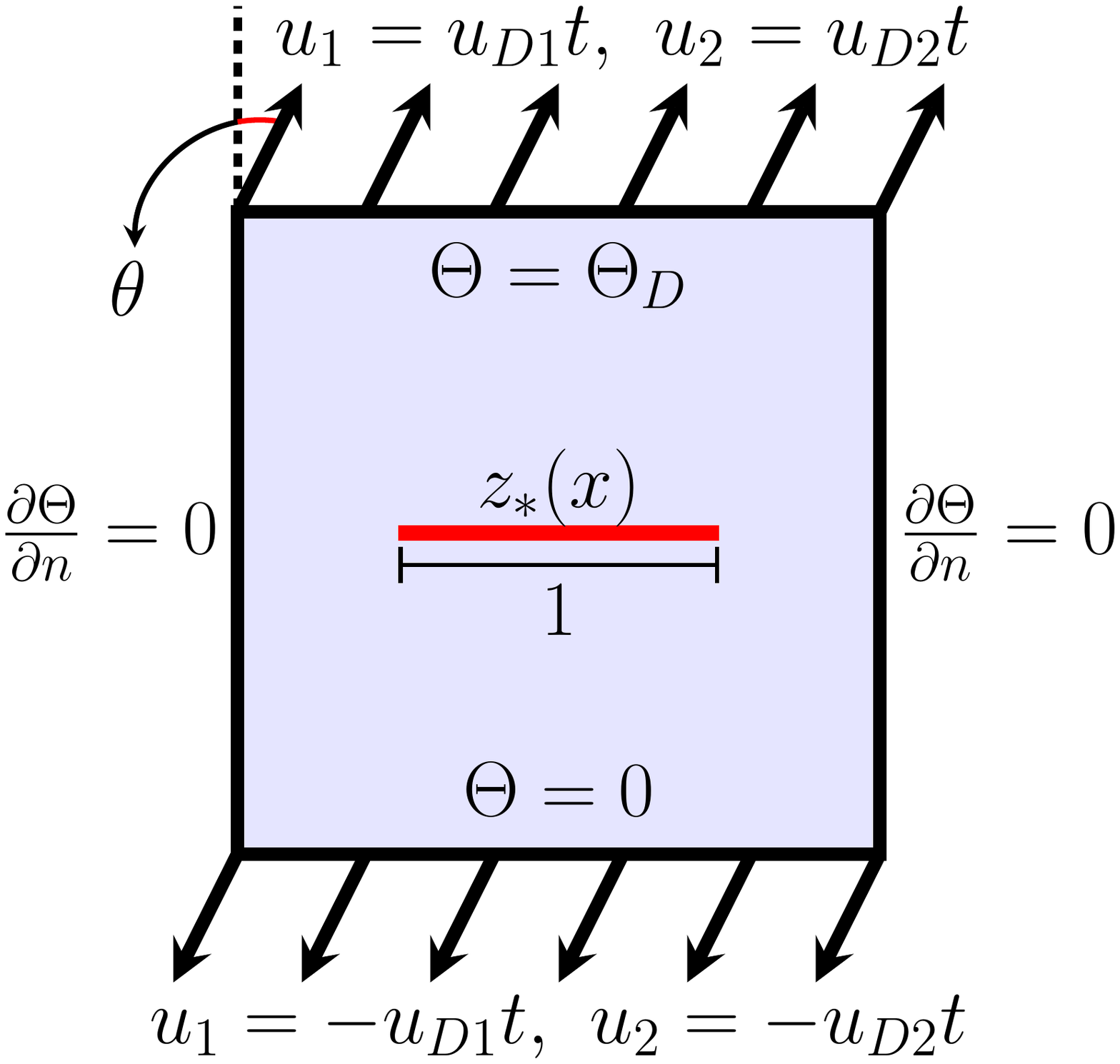}}
	\end{tabular}\\[10pt]
	\caption{Mode I ({left}) and Mode I+II (right) for the study of the crack path under thermal expansion and the loading process. Here, the initial damage $z_{*}(x)$ is illustrated by the red initial crack in the figures.} 
	\label{fig:crackpath0}
\end{figure}
\subsection{Thermoelastic effect on the crack path}\label{SubSec4.4}
In this section, we investigate the effect of the thermoelastic coupling parameter on crack path selection using our proposed models. Under a given temperature gradient, we consider crack propagation of an opening mode (Mode I) and a mixed mode (Mode I+II). In the following numerical examples, we also use the parameters in Table \ref{tab:4}.
\subsubsection{Mode I}\label{SubSec4.4.1}
We use an edge-cracked square domain, which is shown in Figure \ref{fig:crackpath0} (left). We set the domain as follows:
\begin{align}
& C_{\pm} := \left(\begin{array}{l}
-\frac{1}{2}\\
\pm\frac{5}{8}
\end{array}\right) \in \mathbb{R}^{2},\notag\\
& H_{\pm} := \Big\{ x\in \mathbb{R}^{2}; ~\left|x - C_{\pm}\right| \leq \frac{3}{20}\Big\},\notag\\
& \Omega := (-1,1)^{2}\setminus (H_{+} \cup H_{-}),\notag
\end{align}
and we define
\begin{align*}
& \Gamma_{DN1}^{u} := \Gamma \cap \{x_{1} = 1\},~\Gamma_{DN2}^{u} := \partial H_{+} \cup \partial H_{-},~\Gamma_{N}^{u} := \Gamma\setminus(\Gamma_{DN1}^{u} \cup \Gamma_{DN2}^{u}),\\
&\Gamma_{\pm D}^{\Theta} := \Gamma \cap \{x_{2} = \pm 1\},~\Gamma_{N}^{\Theta} := \Gamma\setminus (\Gamma_{+D}^{\Theta} \cup \Gamma_{-D}^{\Theta}).
\end{align*}
The boundary conditions for $u$ and $\Theta$ are given as follows:
\begin{align}
&\left\{
\begin{array}{l}
u_{1}=0\\
\sigma_{12}^{*} = 0
\end{array}
\right. ~ \mbox{on}~ \Gamma_{DN1}^{u}, \quad  \left\{
\begin{array}{l}
(\sigma^{*}n)\cdot e_{1} =0\\
u_{2} = \pm 8t
\end{array}
\right. ~ \mbox{on}~ {\partial H_{\pm}}, \quad
\sigma^{*}[u,\Theta]n = 0~ \mbox{on}~ \Gamma_{N}^{u}, \notag\\[8pt]
&\Theta=\Theta_{D} ~ \mbox{on}~ \Gamma_{+D}^{\Theta}, \quad\Theta = 0 ~ \mbox{on}~ \Gamma_{-D}^{\Theta},  \quad\frac{\partial\Theta}{\partial n} = 0~ \mbox{on}~ \Gamma_{N}^{\Theta}.\notag
\end{align}
The initial condition for $\Theta$ is given as $ \Theta_{*} = 0$.

For $z$, similar to the previous example (Section \ref{SubSec4.3}) , we set $\frac{\partial z}{\partial n} = 0$ on $\Gamma$ and choose the initial value as $z_{*}(x) := \exp{(-(x_{2}/\eta)^{2})}/(1 + \exp{((x_{1}+0.2)/\eta)})$ with $\eta = 1.5\times 10^{-2}$. In this numerical experiment, we apply the thermoelastic coupling parameter $\delta = 0.5$.

\begin{figure}[!h]
	\begin{center}
		\begin{tabular}{cccc}
			{\includegraphics[trim=3cm 0.5cm 2cm 0cm, scale=0.1]{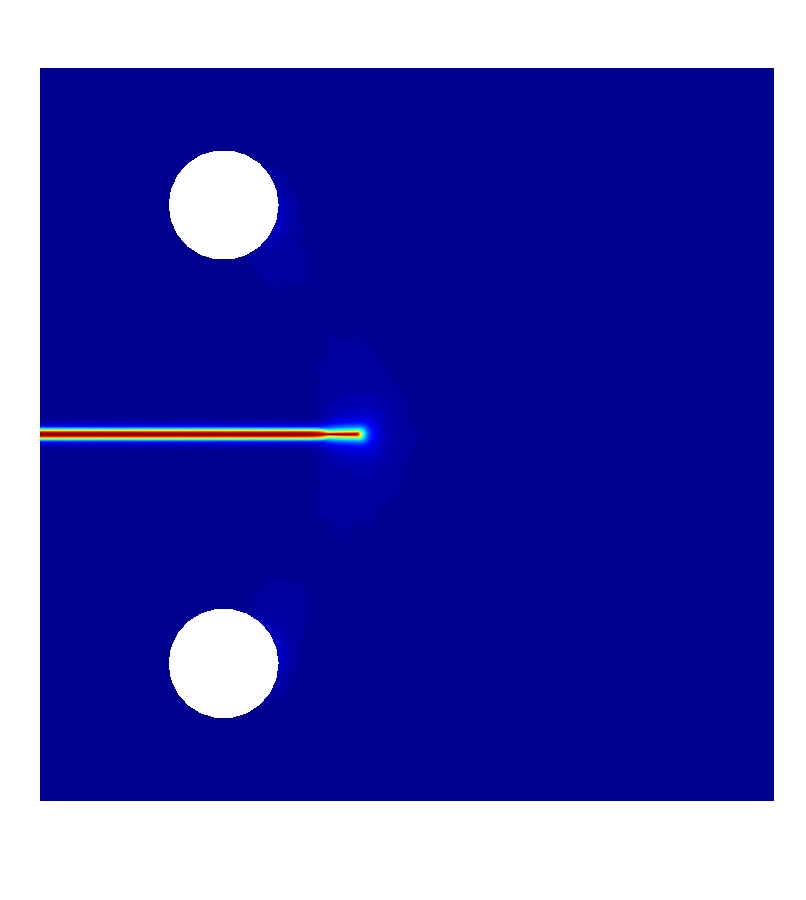}} & \hspace{-10pt}
			{\includegraphics[trim=0cm 0.5cm 2cm 0cm, scale=0.1]{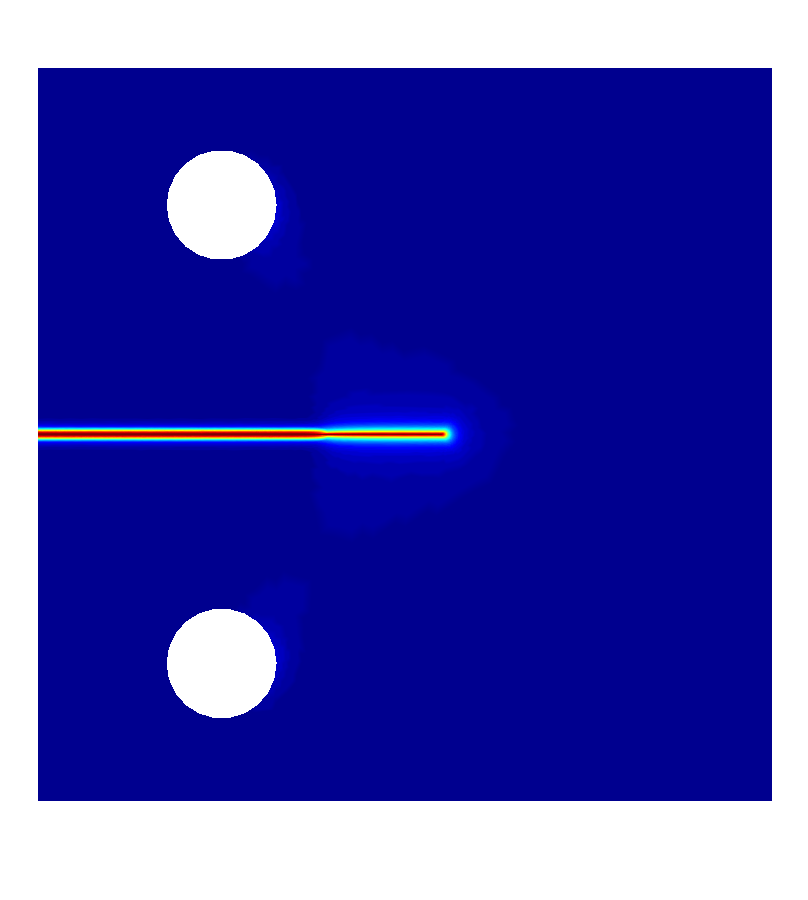}} & \hspace{-10pt}
			{\includegraphics[trim=0cm 0.5cm 2cm 0cm, scale=0.1]{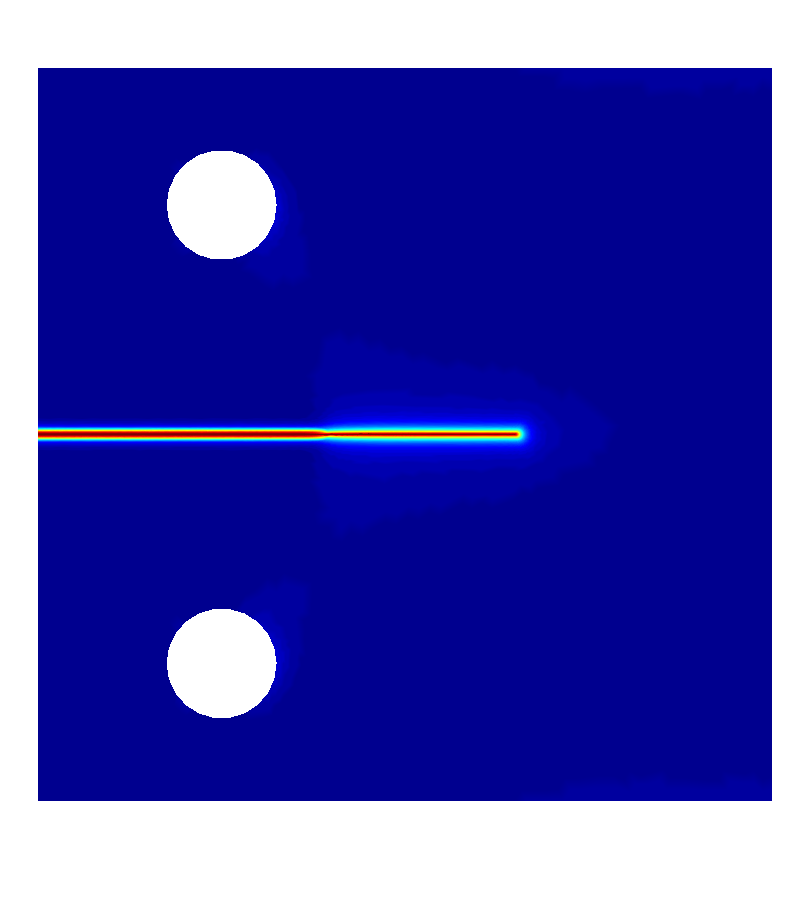}} & \hspace{-10pt}
			{\includegraphics[trim=0cm 0.5cm 2cm 0cm, scale=0.1]{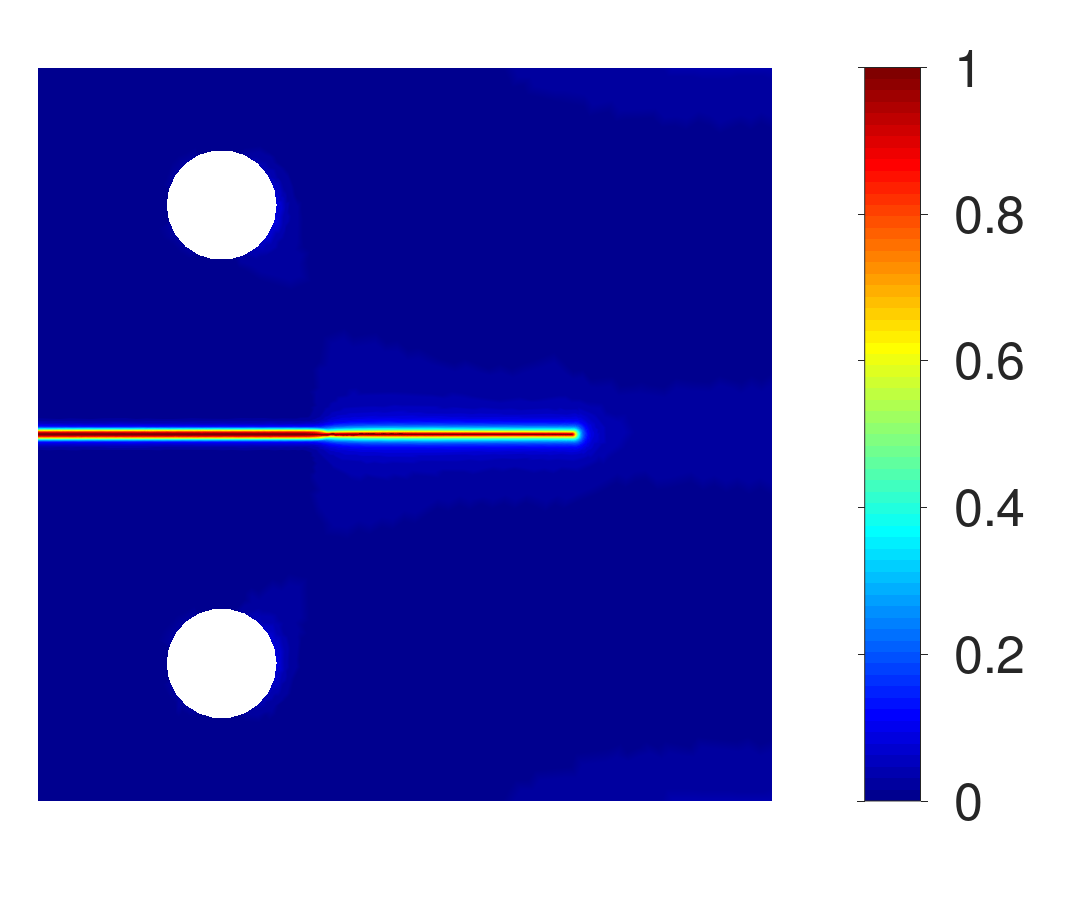}} \\[-10pt]
			{\includegraphics[trim=3cm 0.5cm 2cm 0cm, scale=0.1]{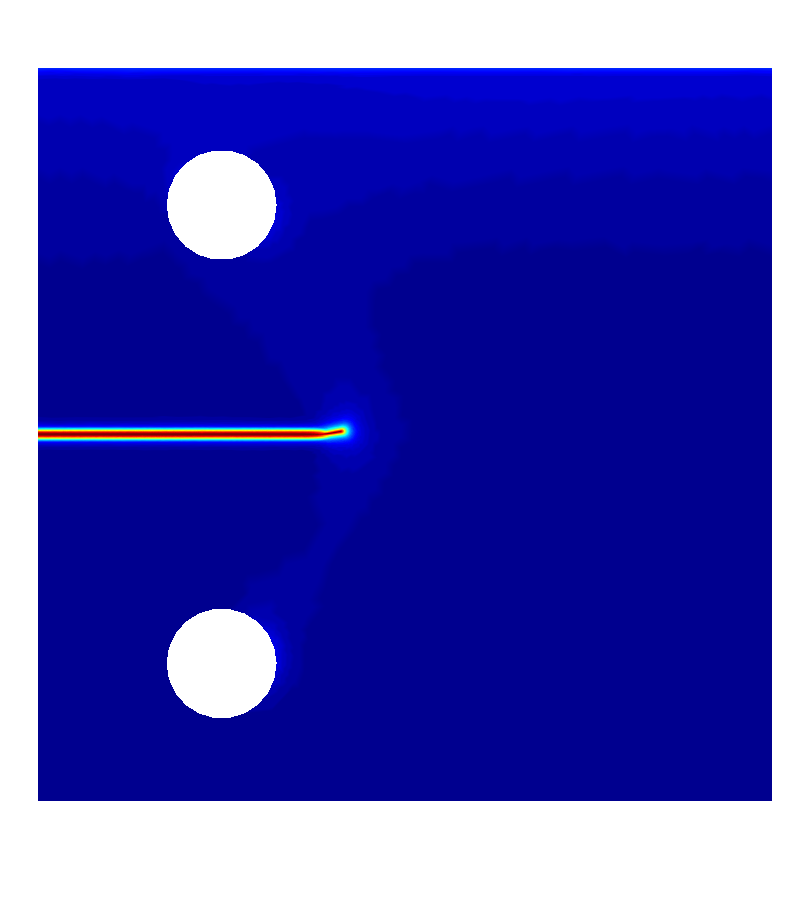}} & \hspace{-10pt}
			{\includegraphics[trim=0cm 0.5cm 2cm 0cm, scale=0.1]{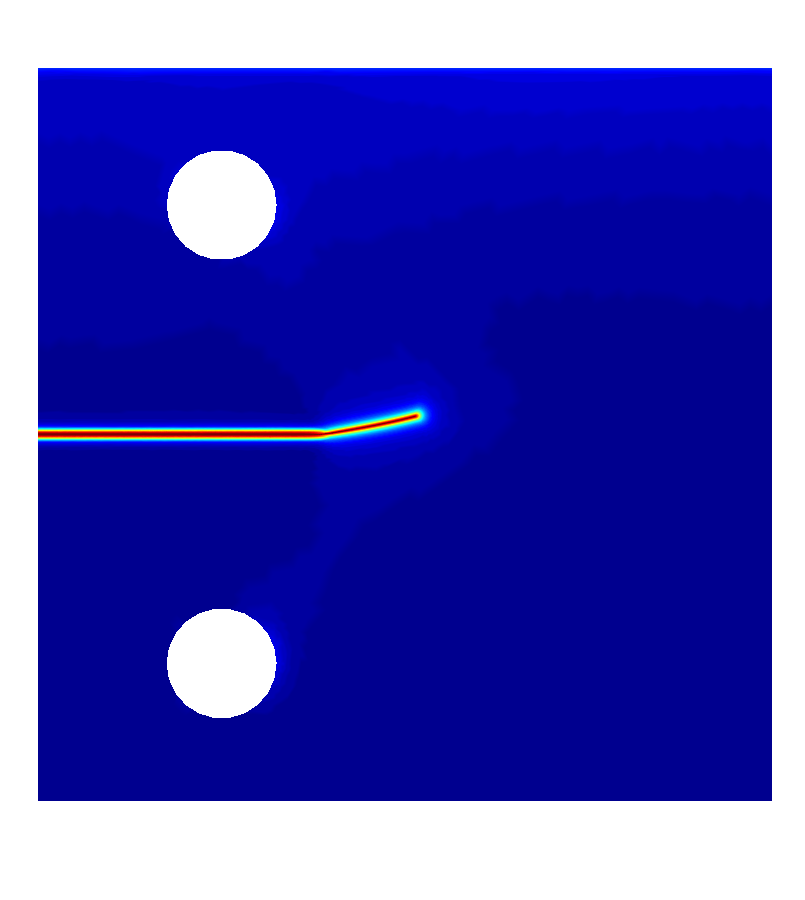}} & \hspace{-10pt}
			{\includegraphics[trim=0cm 0.5cm 2cm 0cm, scale=0.1]{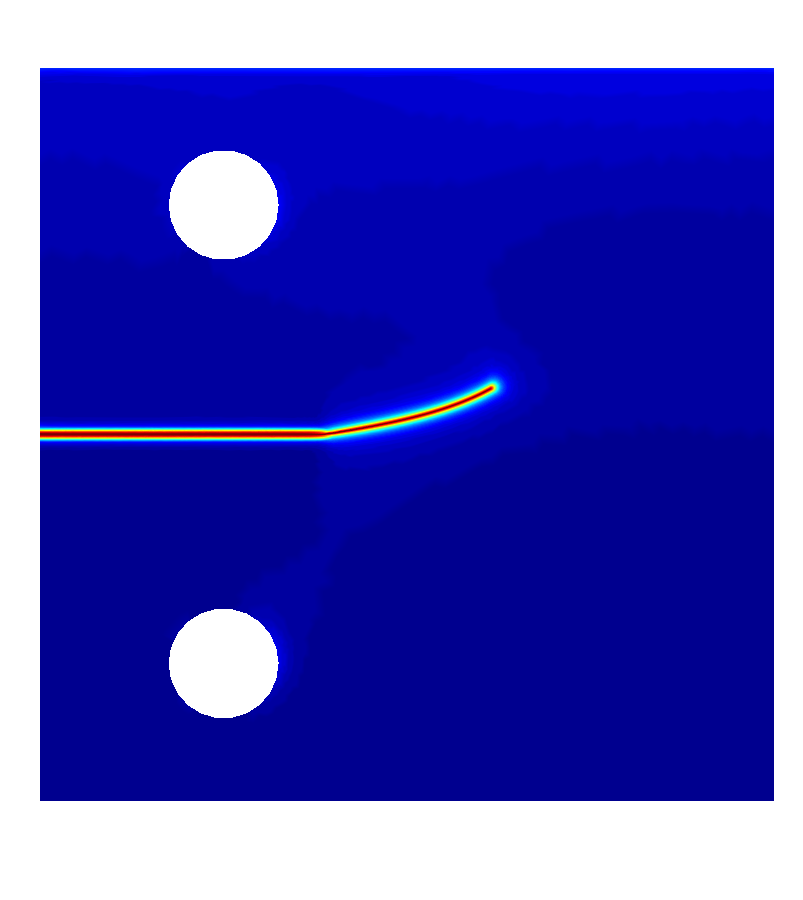}} & \hspace{-10pt}
			{\includegraphics[trim=0cm 0.5cm 2cm 0cm, scale=0.1]{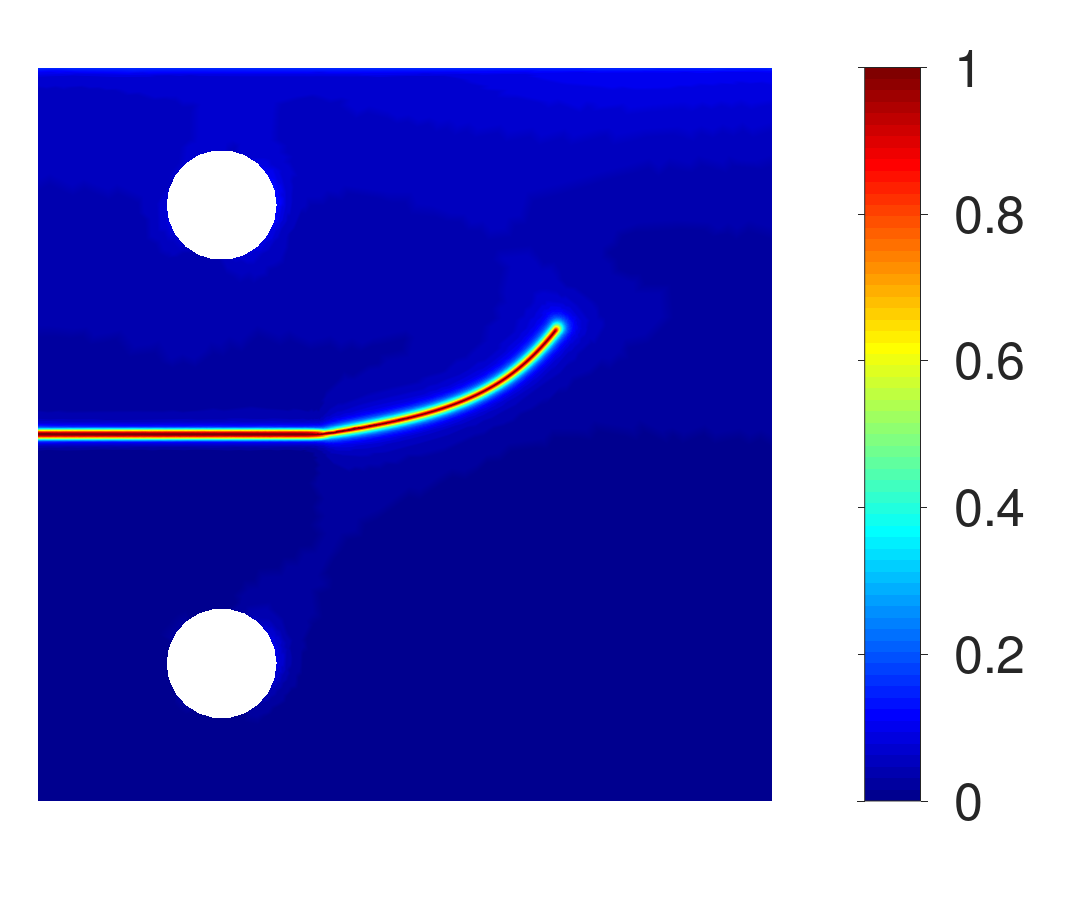}} \\[-10pt]
			{\includegraphics[trim=3cm 0.5cm 2cm 0cm, scale=0.1]{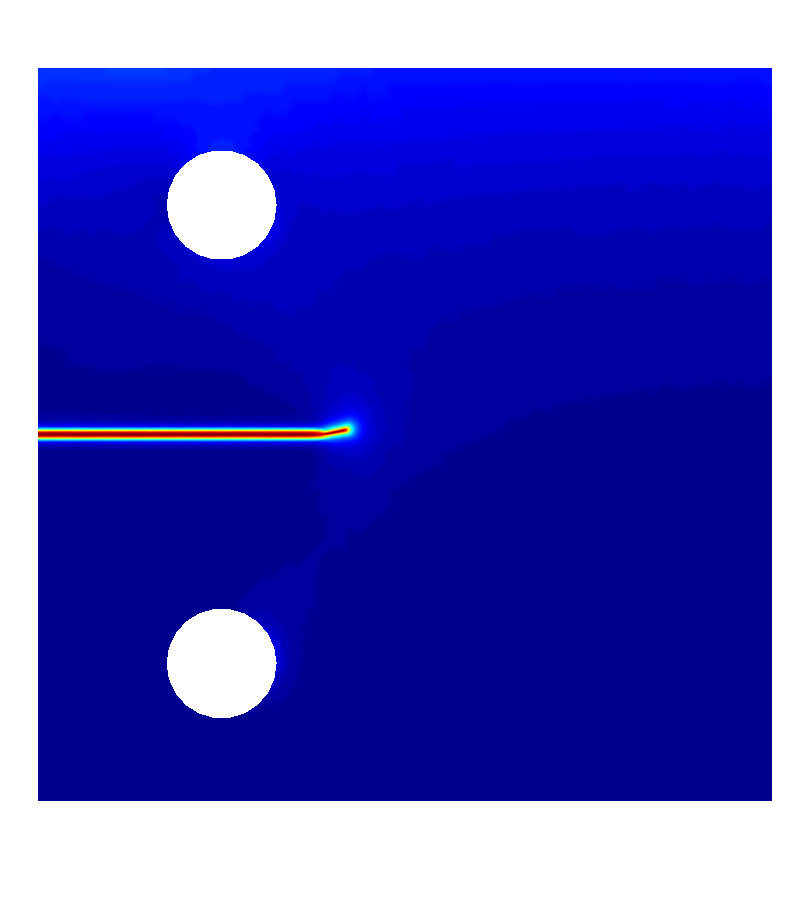}} & \hspace{-10pt}
			{\includegraphics[trim=0cm 0.5cm 2cm 0cm, scale=0.1]{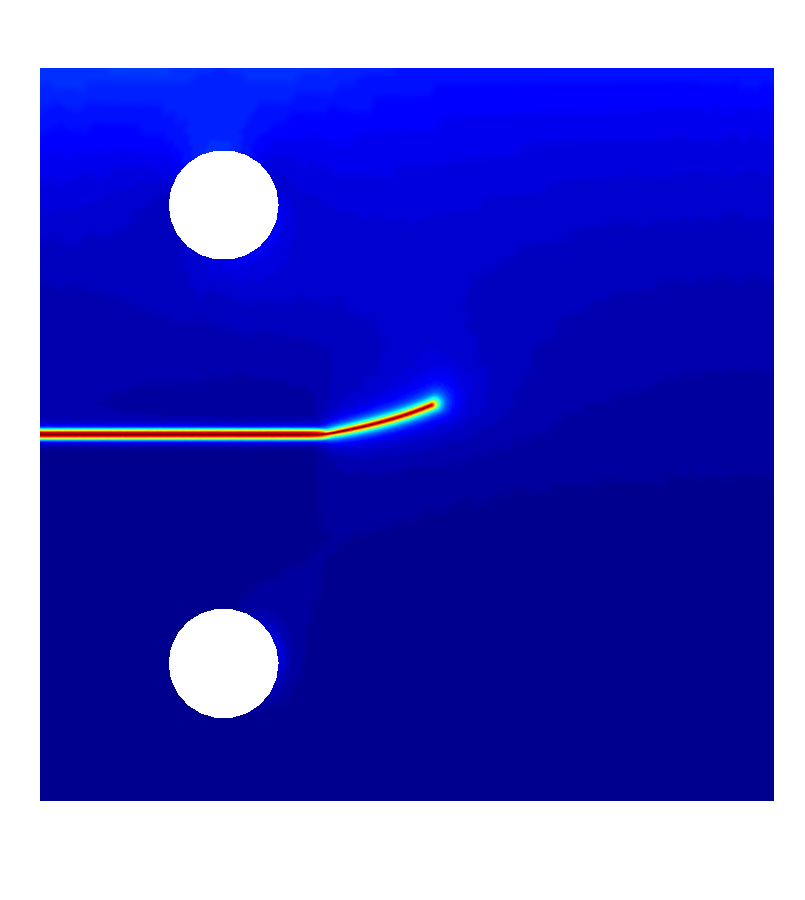}} & \hspace{-10pt}
			{\includegraphics[trim=0cm 0.5cm 2cm 0cm, scale=0.1]{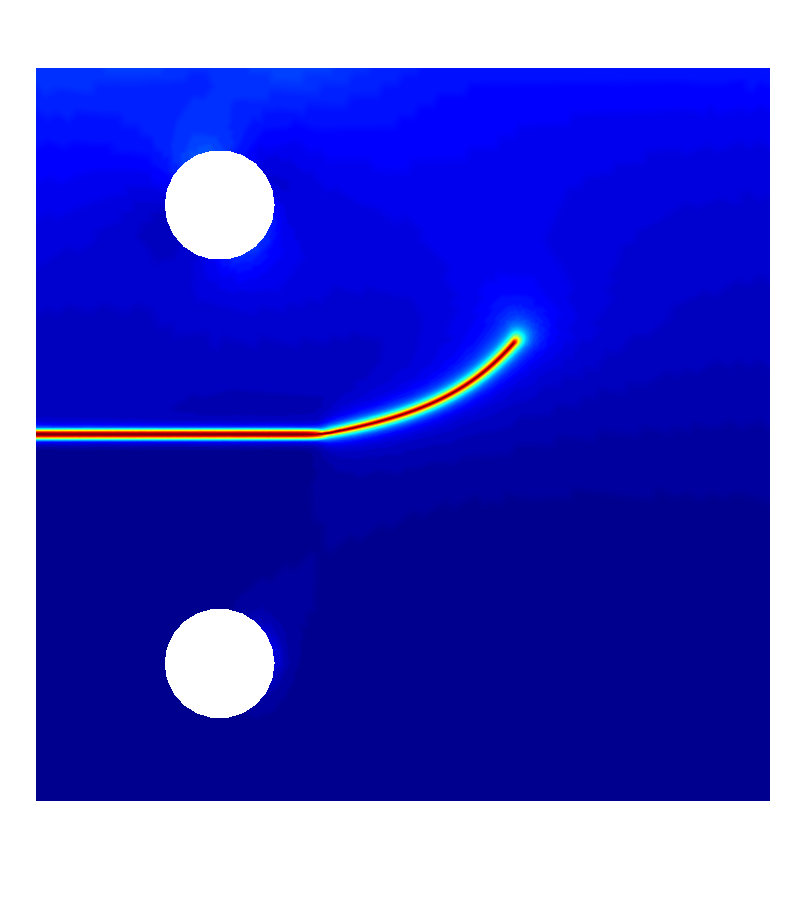}} & \hspace{-10pt}
			{\includegraphics[trim=0cm 0.5cm 2cm 0cm, scale=0.1]{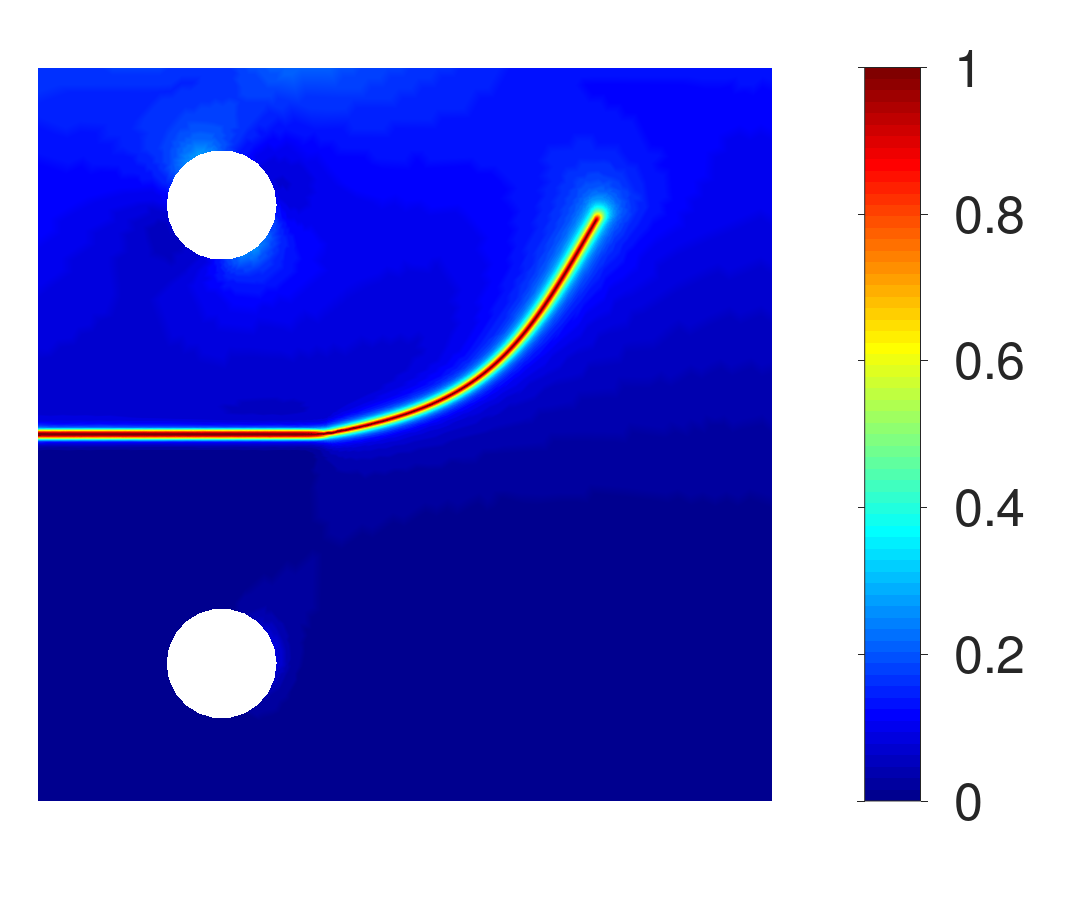}} 
		\end{tabular}\\[-25pt]
		\caption{Snapshots of the crack paths.  F-PFM (upper),  TF-PFM1 (middle),  and TF-PFM2 (lower) at $t = 0.4, ~0.6, ~0.8, ~1$ (left to right). For TF-PFM1 and TF-PFM2, we set $\Theta_{D} = 10$ and $\delta = 0.5$. Here, the color represents the value of $z$.}\label{fig:crackpath1}
	\end{center}
\end{figure}

Figure \ref{fig:crackpath1} shows the different crack paths obtained by the three models when $\Theta_{D} = 10$. Straight cracks occur in the F-PFM path since the thermal effect is ignored there. On the other hand, crack curves occur in the TF-PFM1 and TF-PFM2 paths. Here, the crack path is more curved in the TF-PFM2 path than in the TF-PFM1 path. These results show good qualitative agreement with the results reported in \cite{Jaskowiec2017}.  

\begin{figure}[!h]
	\begin{center}
		\begin{tabular}{cccc}
			{\includegraphics[trim=1cm 0.25cm 2cm 0cm, scale=0.225]{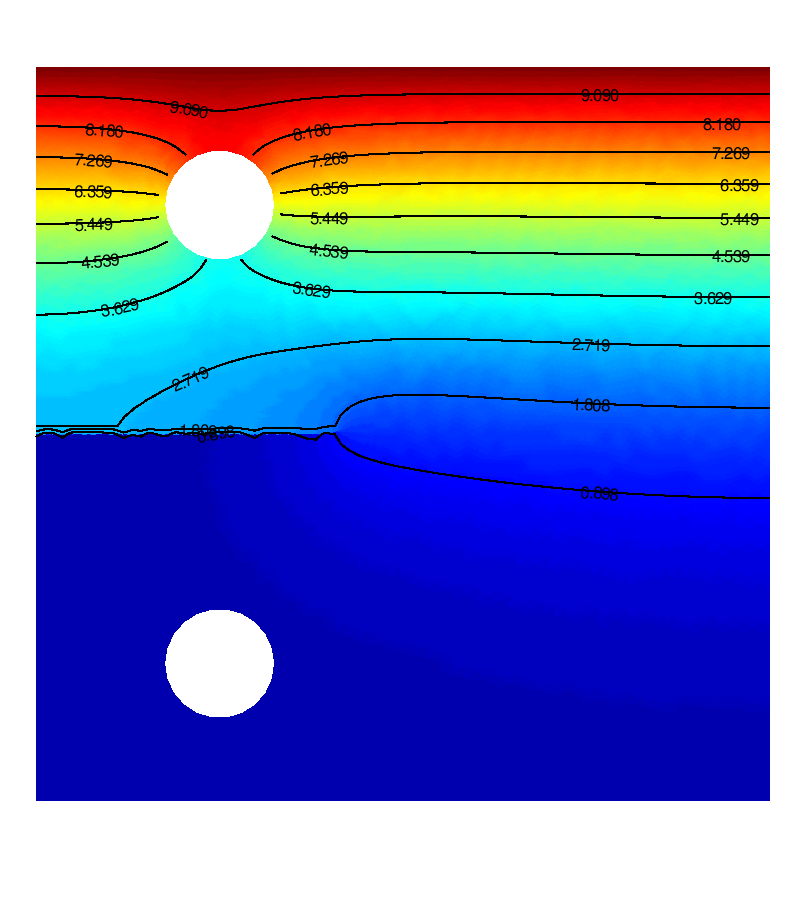}} & \hspace{-10pt}
			{\includegraphics[trim=0cm 0.25cm 2cm 0cm, scale=0.225]{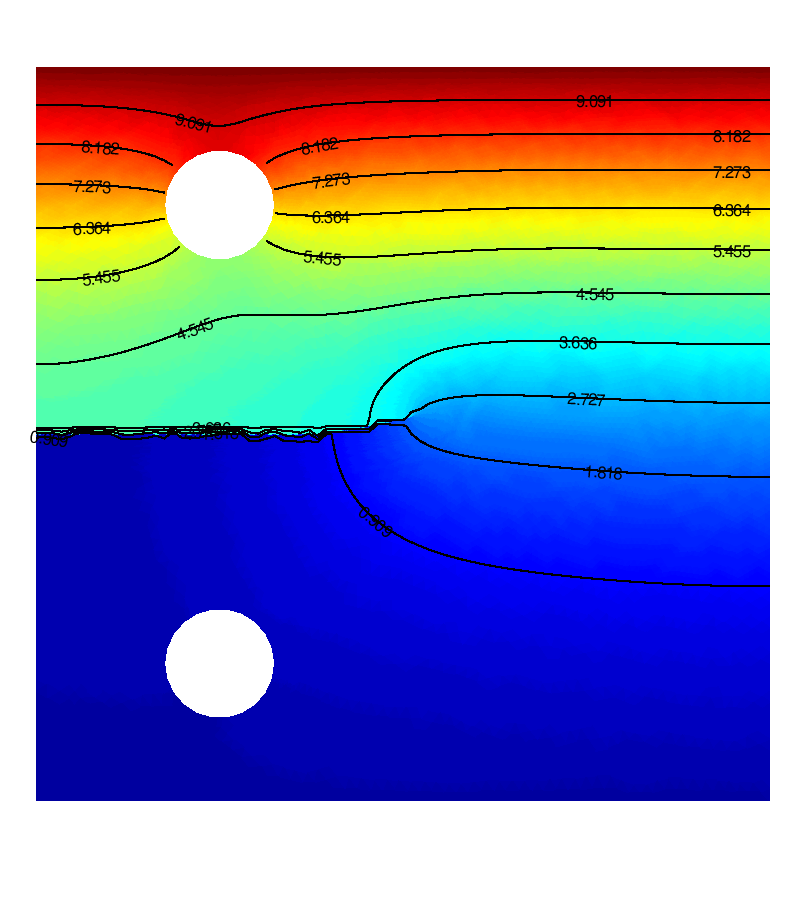}} & \hspace{-10pt}
			{\includegraphics[trim=0cm 0.25cm 2cm 0cm, scale=0.225]{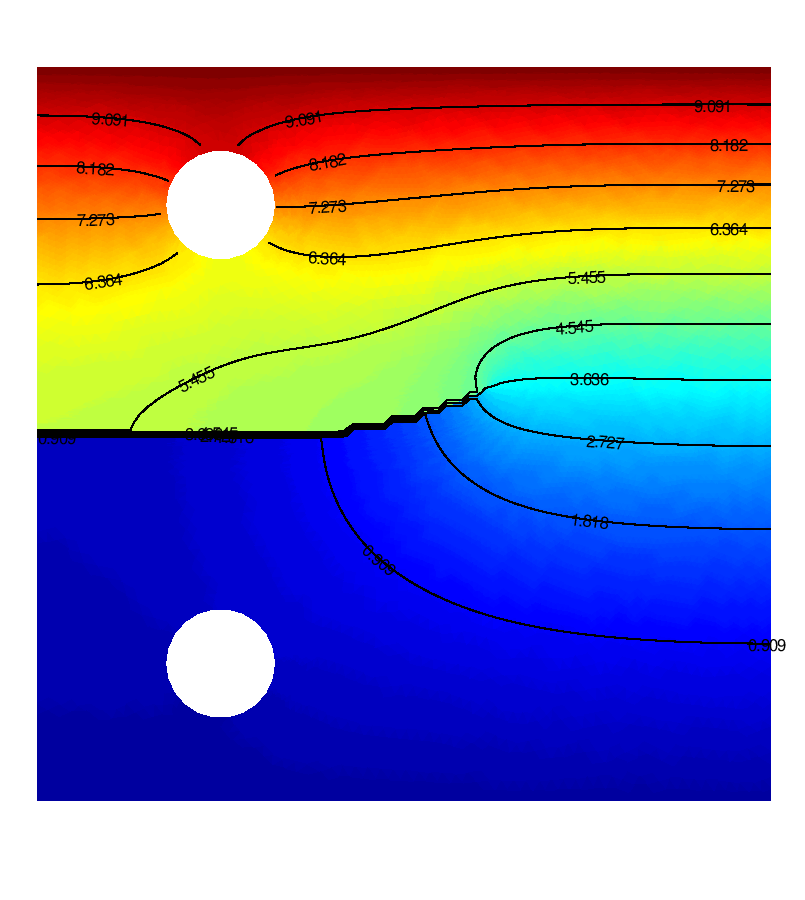}} & \hspace{-10pt}
			{\includegraphics[trim=0cm 0.25cm 2cm 0cm, scale=0.225]{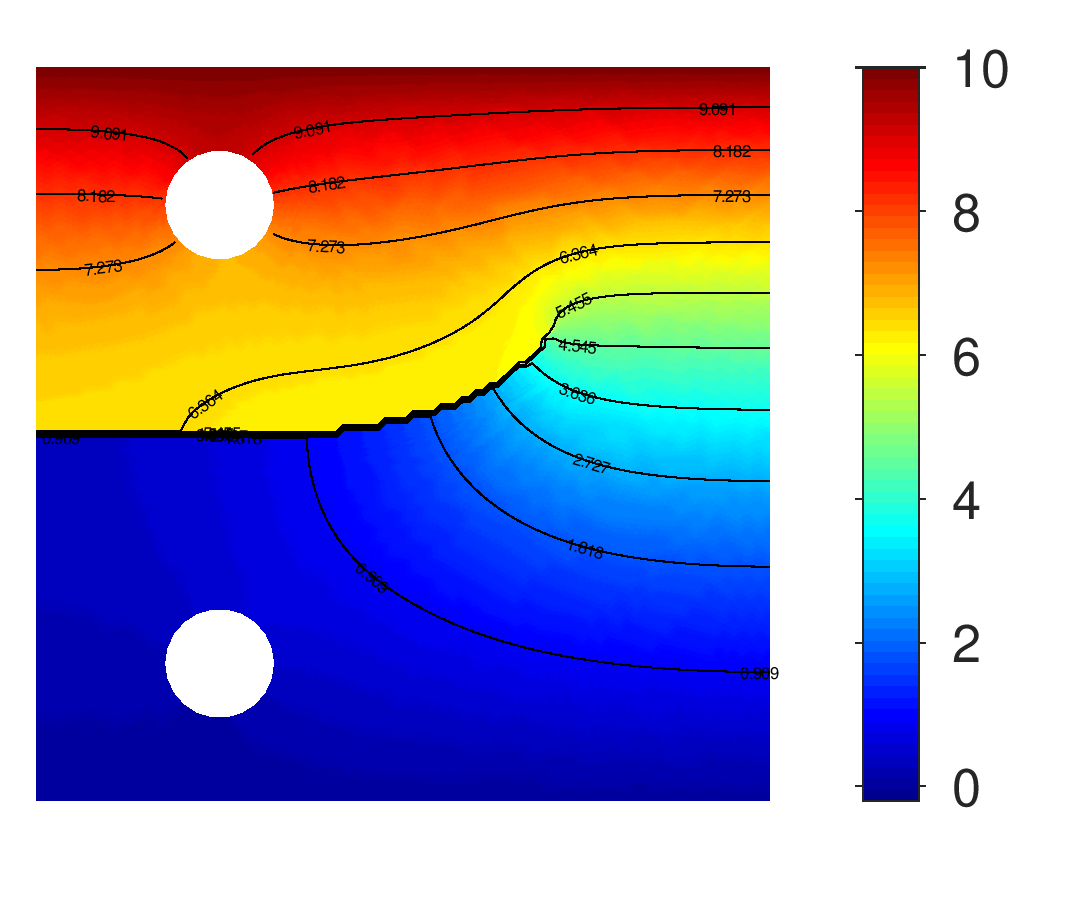}} \\[-12pt]
			{\includegraphics[trim=1cm 0.25cm 2cm 0cm, scale=0.225]{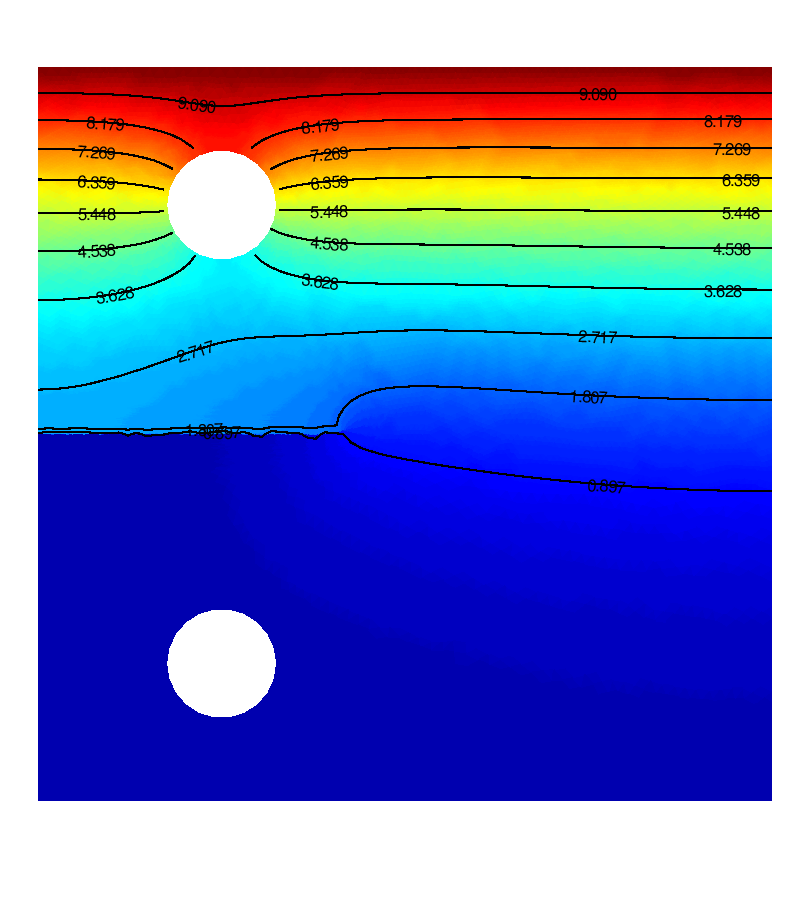}} & \hspace{-10pt}
			{\includegraphics[trim=0cm 0.25cm 2cm 0cm, scale=0.225]{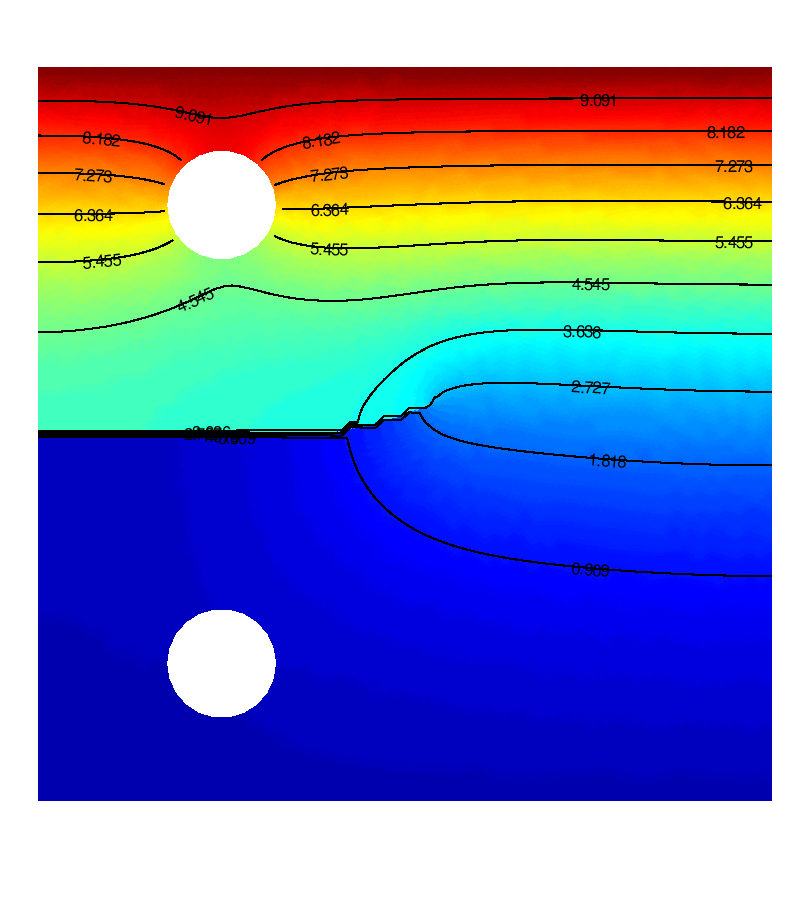}} & \hspace{-10pt}
			{\includegraphics[trim=0cm 0.25cm 2cm 0cm, scale=0.225]{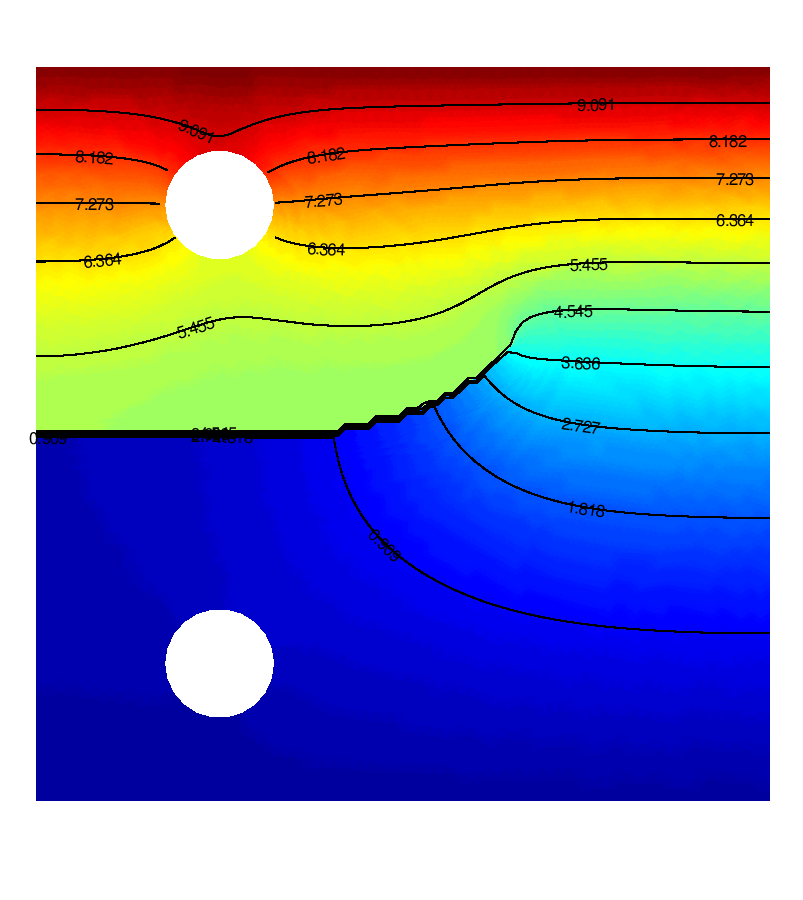}} & \hspace{-10pt}
			{\includegraphics[trim=0cm 0.25cm 2cm 0cm, scale=0.225]{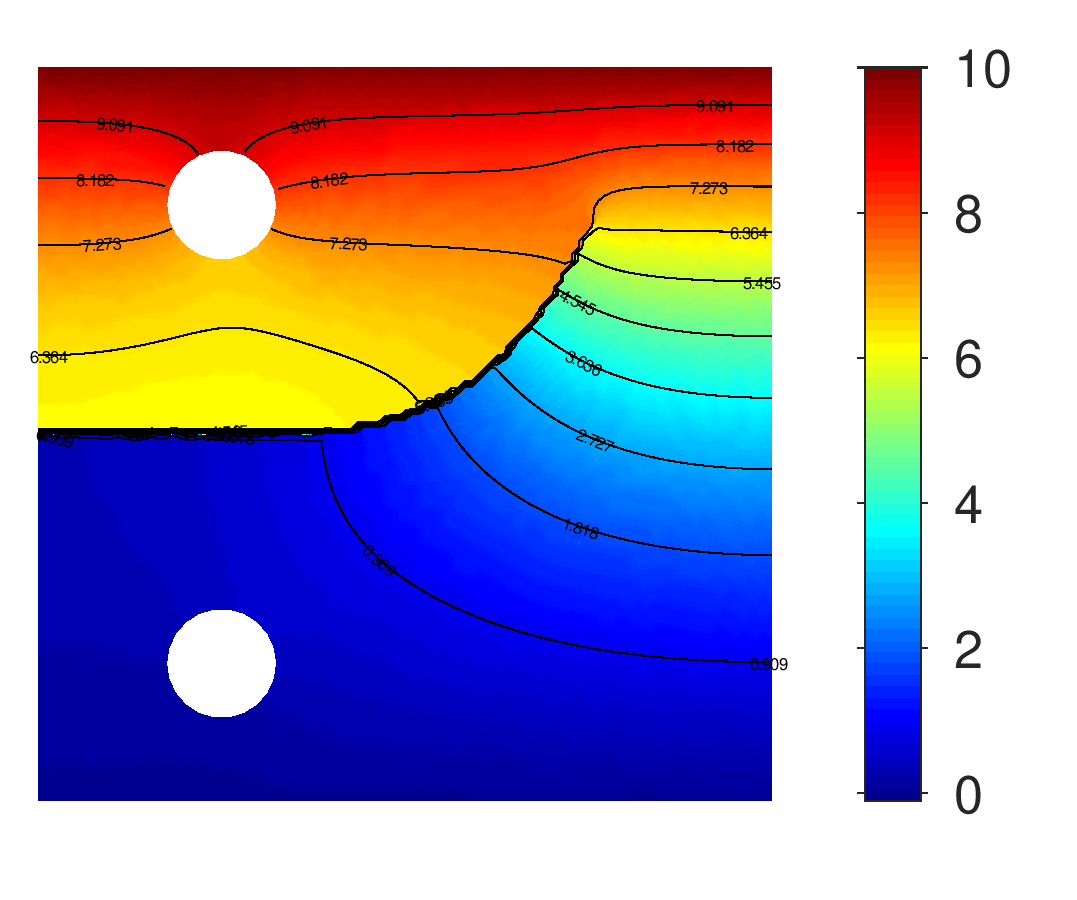}} 
		\end{tabular}\\[-25pt]
		\caption{Snapshots of the temperature gradient during thermal expansion and crack growth under the given temperature $\Theta_{D} = 10$. TF-PFM1 (top)  and TF-PFM2 (bottom) at $t = 0.4, ~0.6, ~0.8, ~1$ (left to right); the color represents the value of $\Theta$. }\label{fig:crackpath2}
	\end{center}
\end{figure}
\begin{figure}[!h]
	\begin{tabular}{ccc}
		\includegraphics[trim=4cm 7cm 3cm 6cm,width=0.45\textwidth]{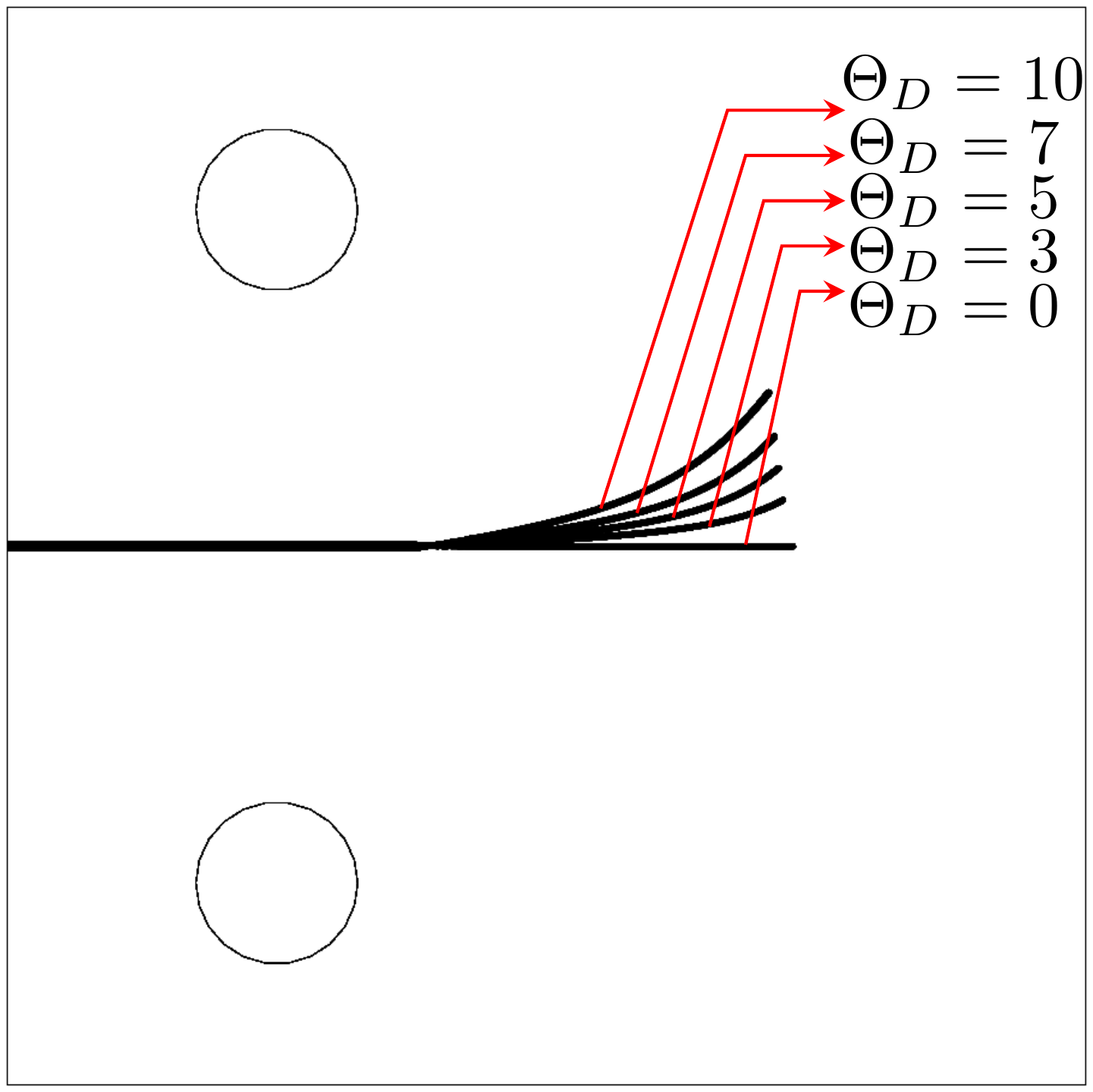}
		&
		{\includegraphics[trim=4cm 7cm 3cm 6cm,width=0.45\textwidth]{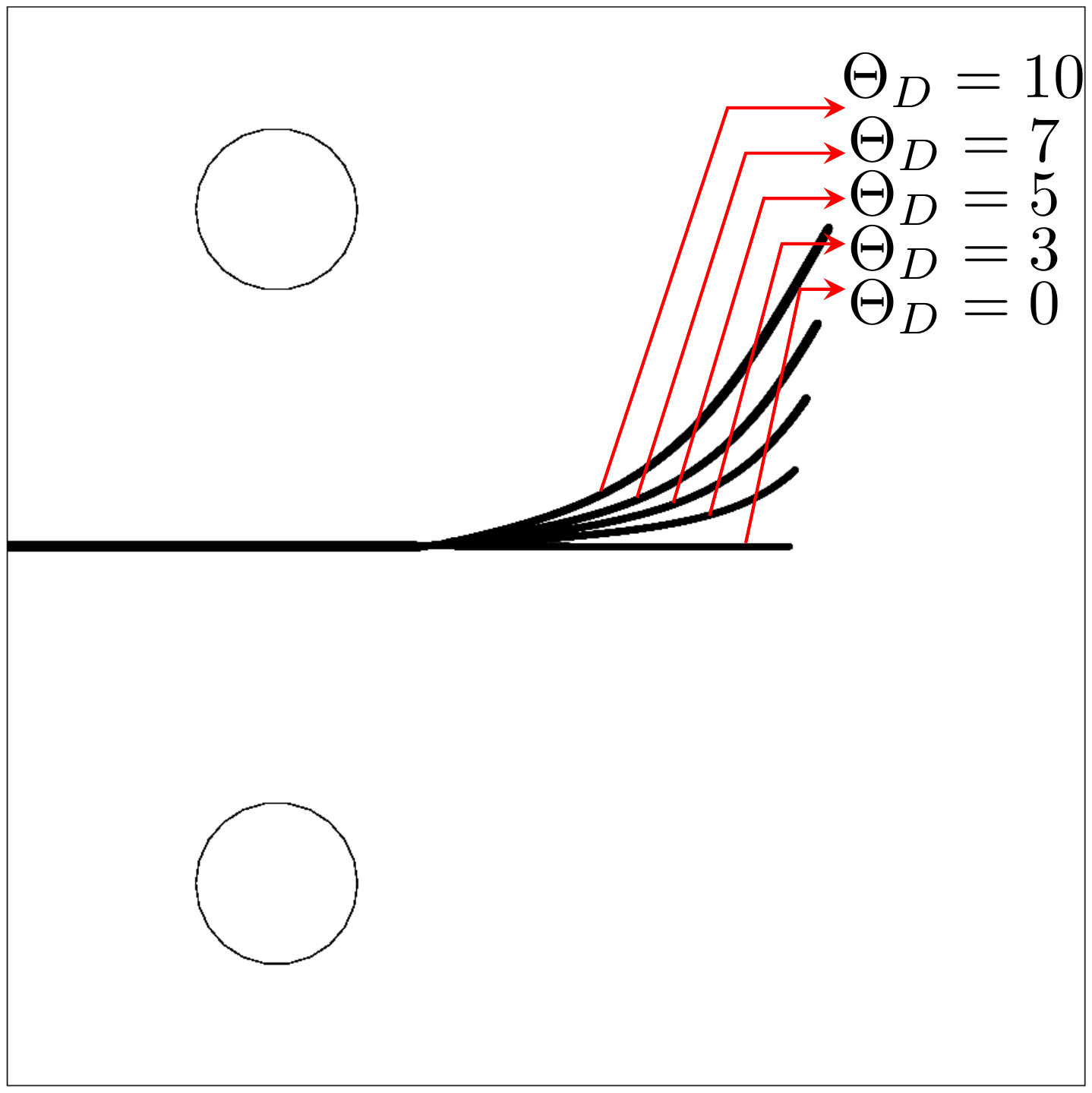}}
	\end{tabular}
		\caption{Comparison of the crack paths using TF-PFM1 (left) and TF-PFM2 (right) with the given temperature variations under Mode I at the final computational time $t=1$.}\label{fig:crackpath3}
\end{figure}

Figure \ref{fig:crackpath3} shows the crack paths for different temperature gradients $\Theta_{D} = 0,~3,~5,~7,~10$ obtained by TF-PFM1 (left) and TF-PFM2 (right). A larger temperature gradient generates a more curved crack path, and TF-PFM2 obtains a more curved crack path than TF-PFM1. Both have significant differences in the magnitude of angle deviation but have the same crack path directions. Therefore, it is clear that thermal expansion changes the crack path.

The temperature distributions during crack growth are shown in Figure \ref{fig:crackpath2}. There exists a temperature discontinuity along the crack path, which is caused by $\kappa(z) =(1-z)^{2}\kappa_{0}$. It approximately represents a thermal insulation condition across the crack. 

\subsubsection{Mode I+II}\label{SubSec4.4.2}
According to the numerical experiment in \cite{Kimura2021}, we consider the following setting for mixed mode crack propagation under a thermal gradient. Let $\Omega := (-1,1)^{2} \in \mathbb{R}^{2}$, as shown in Figure \ref{fig:crackpath0} (right), and $\Gamma := \partial\Omega$. We set
\begin{align*}
&\Gamma_{\pm D}^{u} := \Gamma \cap \{x_{2} = \pm 1\},~\Gamma_{N}^{u} := \Gamma \setminus (\Gamma_{+ D}^{u} \cup \Gamma_{- D}^{u}),\\
& \Gamma_{\pm D}^{\Theta} := \Gamma \cap \{x_{2} = \pm 1\},,~\Gamma_{N}^{\Theta} := \Gamma\setminus (\Gamma_{+D}^{\Theta} \cup \Gamma_{-D}^{\Theta}).
\end{align*}
The boundary conditions for $u$ are given as follows:
\begin{align}
&\left\{
\begin{array}{ll}
u_{1}=\pm 3\sin(\pi/3)t,\\
u_{2} = \pm 3\cos(\pi/3)t
\end{array}
\right. ~ \mbox{on}~ \Gamma_{\pm D}^{u}, \qquad \sigma^{*}[u,\Theta]n = 0~ \mbox{on}~ \Gamma_{N}^{u}. \notag
\end{align}
The boundary conditions for $\Theta$ and $z$ are the same as those in Section \ref{SubSec4.4.1}. The initial crack profile is given as $z_{*}(x) := \exp{(-(x_{2}/\eta)^{2})}/(1 + \exp{((x_{1}-0.5)/\eta)}) - \exp{(-(x_{2}/\eta)^{2})}/(1 + \exp{((x_{1}+0.5)/\eta)})$ with $\eta = 1.5\times 10^{-2}$. We fix the thermoelastic coupling parameter $\delta = 0.15$ and change the temperature gradient to $\Theta_{D} = 0,~ 2, ~ 3,~ 5, ~ 6$.
\begin{figure}[!h]
	\begin{tabular}{cc}
		\includegraphics[trim=4cm 7cm 3cm 5cm,width=0.45\textwidth]{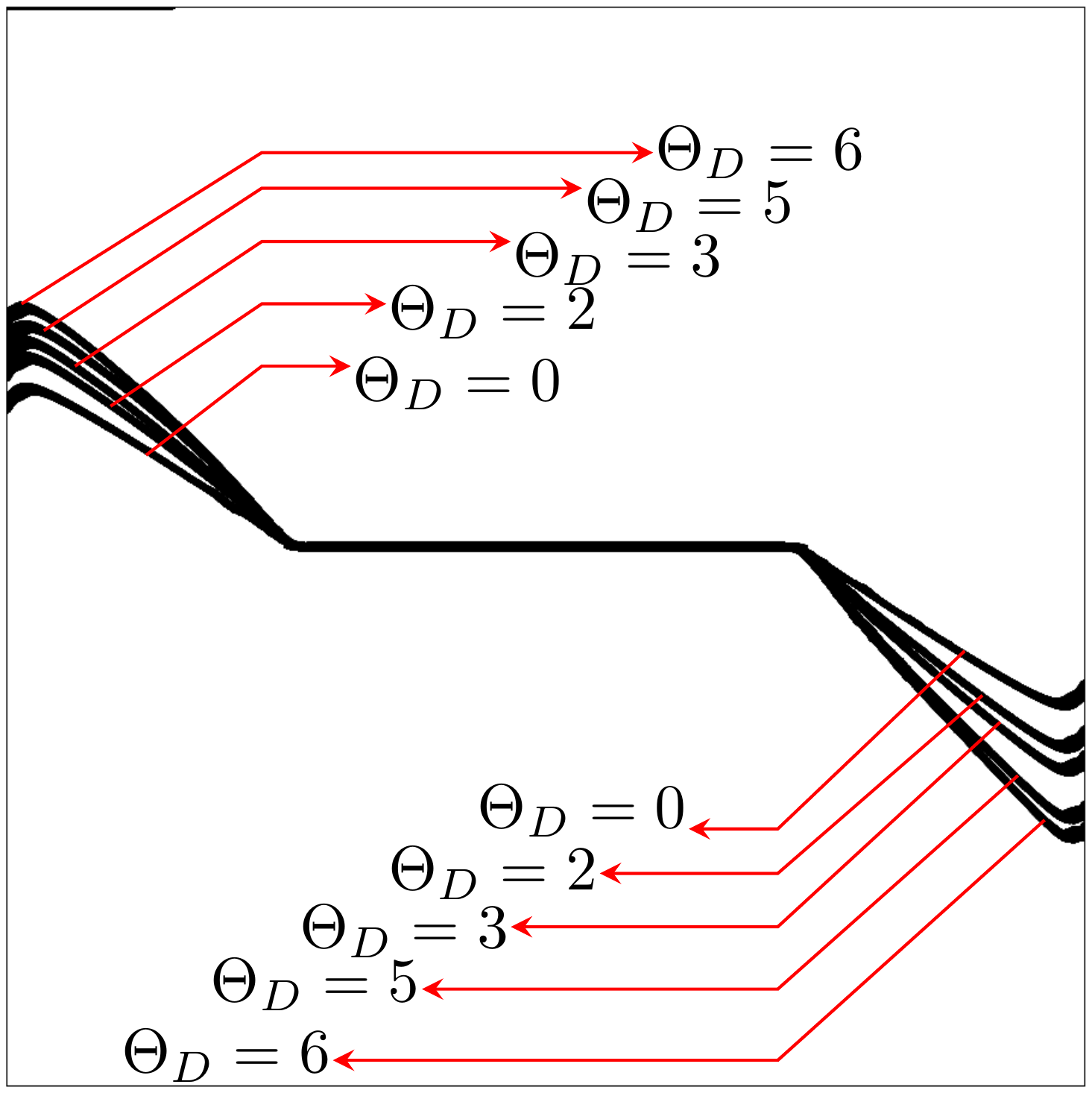}
		&
		{\includegraphics[trim=4cm 7cm 3cm 5cm,width=0.45\textwidth]{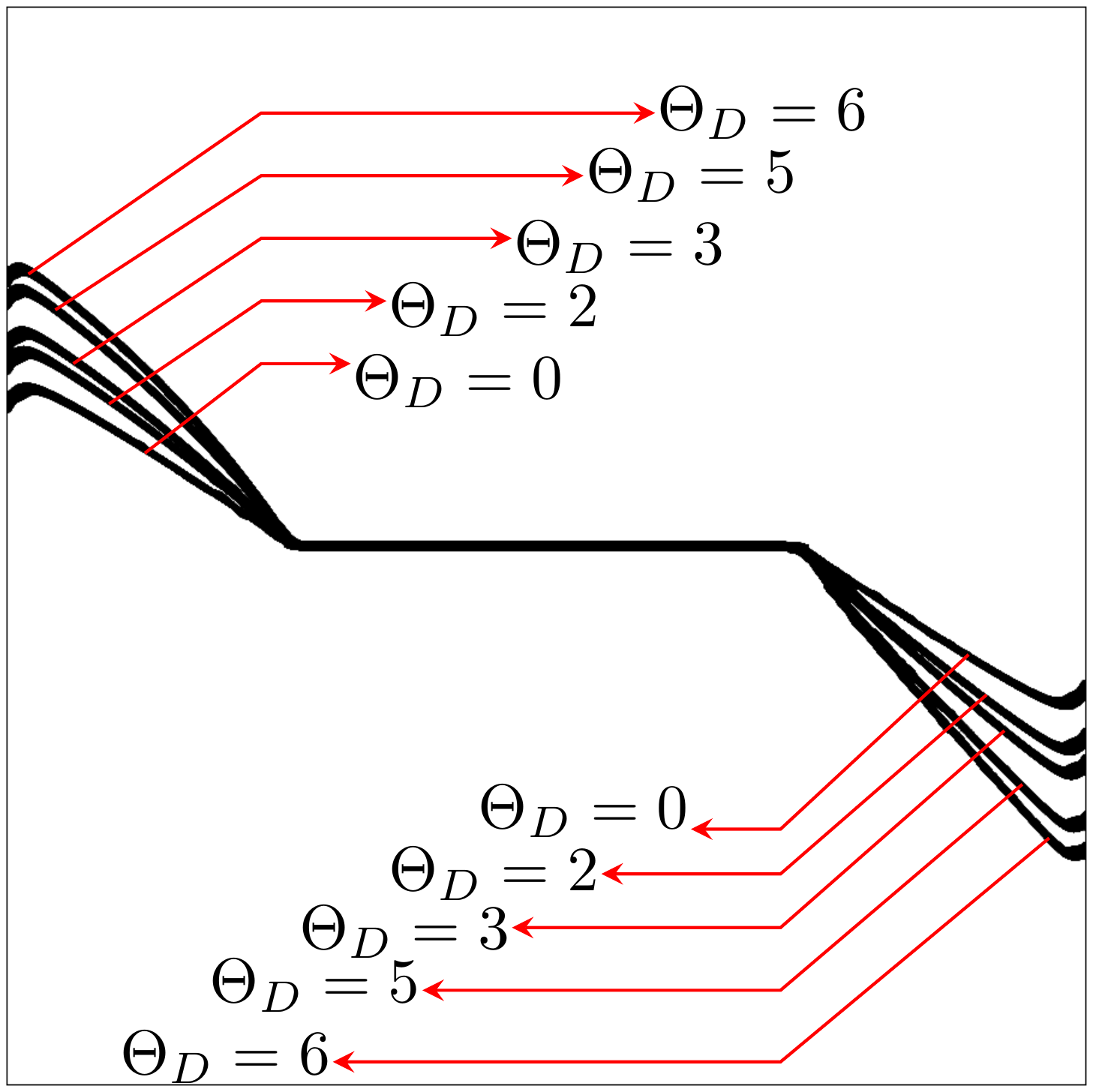}}
	\end{tabular}
	\caption{Comparison of the crack paths using TF-PFM1 (left) and TF-PFM2 (right) with the given temperature variations under Mode I+II at the final computational time.} 
	\label{fig:crackpathMIXED}
\end{figure}

\begin{figure}[!h]
	\centering
	\begin{tabular}{cc}
		\includegraphics[trim=4cm 7cm 3cm 6cm,width=0.45\textwidth]{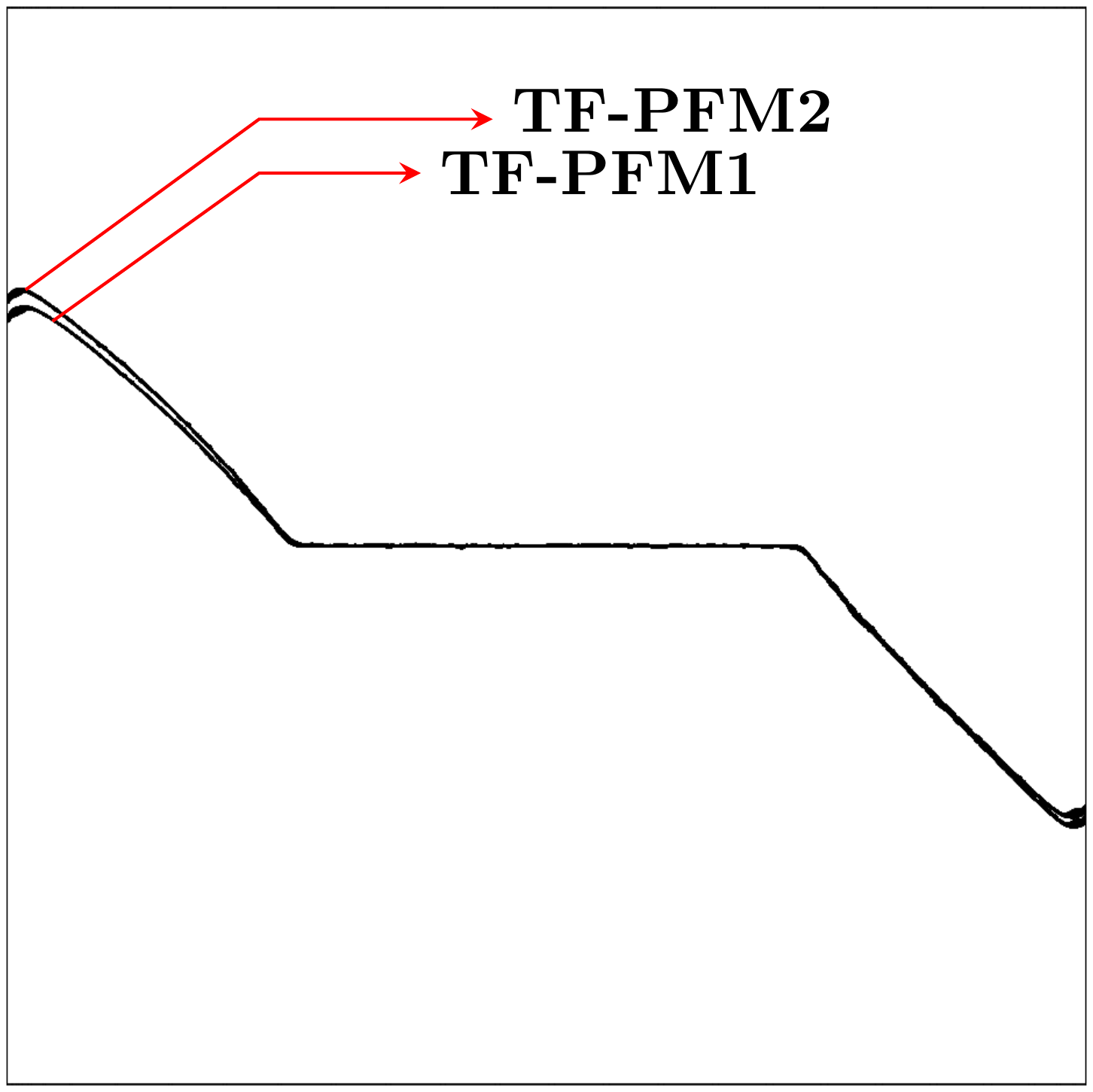}
		&
		\includegraphics[trim=4cm 7cm 3cm 6cm,width=0.45\textwidth]{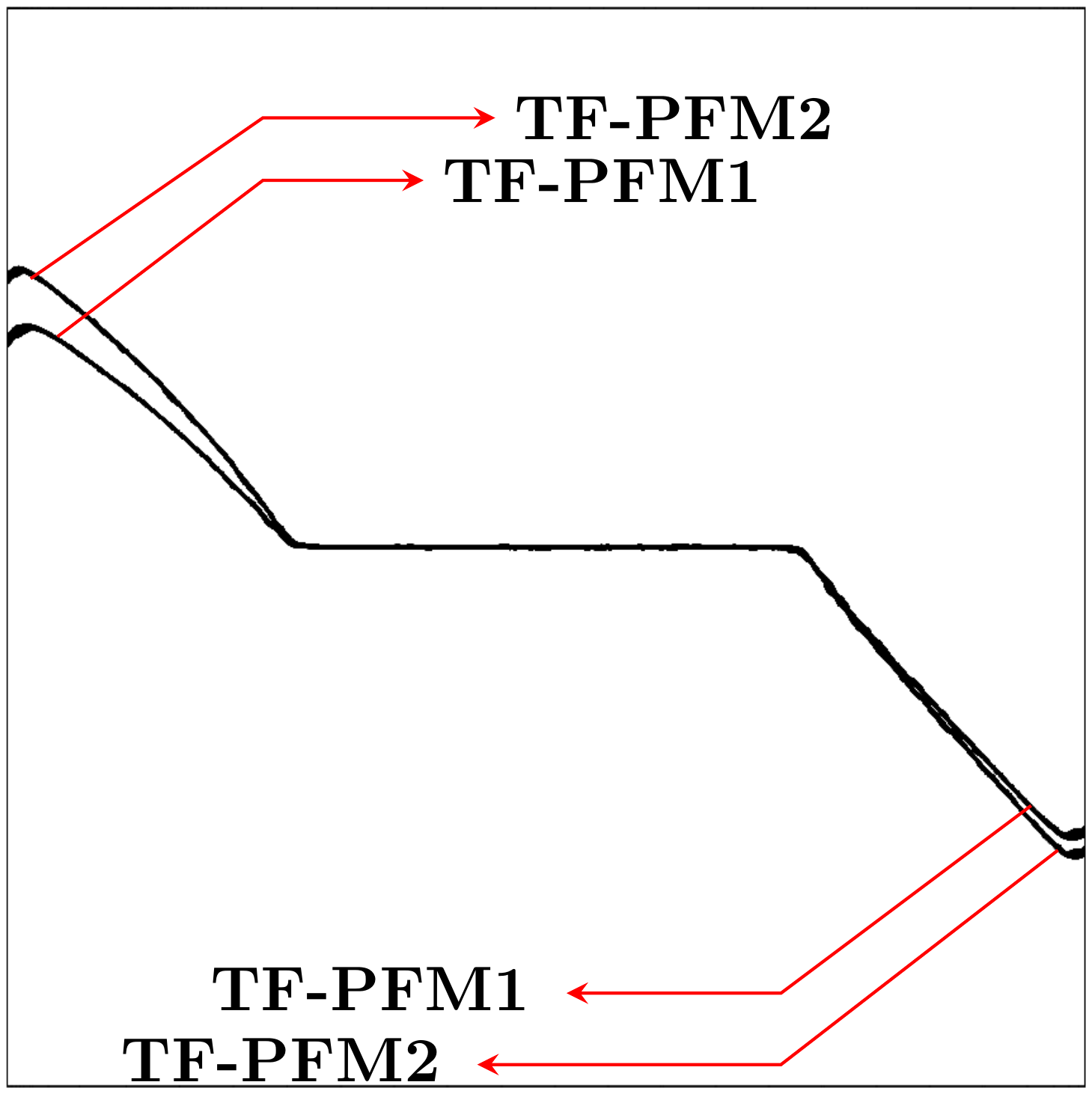}
	\end{tabular}
	\caption{Comparison of the crack paths using TF-PFM1 and TF-PFM2 when $\Theta = 5$  (left) and $\Theta = 6$ (right) at the final computational time.} 
	\label{fig:Compare}
\end{figure}

Figure \ref{fig:crackpathMIXED} shows the crack paths obtained by TF-PFM1 and TF-PFM2. The cracks are kinked, and the kink angle becomes larger when the thermal gradient $\Theta_{D}$ increases. The two models provide similar results, but the kink angle in the TF-PFM2 crack is larger than that in the TF-PFM1 crack, as shown in Figure \ref{fig:Compare}. Therefore,  we conclude that thermal expansion changes the crack path. 

Here, we do not show the temperature distribution during thermal expansion. We observe that the temperature distribution is quite similar to that of Mode I in Section \ref{SubSec4.4.1}, and a temperature discontinuity exists along the crack path during temperature injection.

\section{Summary and Future Works} \label{sec:4}
We proposed two thermal fracturing phase field models, TF-PFM1 and TF-PFM2, by coupling the Biot thermoelasticity model \cite{Biot1956} and the fracturing phase field model (F-PFM) by Takaishi-Kimura \cite{Kimura2021,Kimura2009}.

For the Biot model, we studied  a variational principle (Proposition \ref{Proposition:2}) and energy equality (Theorem \ref{Theorem1}), which were related to different energies $E^{*}_{el}(u,\Theta)$ and $E_{el}(u) + E_{th}(\Theta)$, respectively (see Tables \ref{tab:1} and \ref{tab:2}).

On the other hand, F-PFM has a gradient flow structure with respect to the total energy $\mathcal{E}_{el}(u,z)+ E_{s}(z)$ and admits energy equality (Theorem \ref{Theorem2:a}).

As the first model, TF-PFM1 was derived based on the variational principle of the Biot model and the gradient flow structure of F-PFM, while TF-PFM2 is based on the energy equalities of the Biot model and F-PFM. The difference between them is the driving force term for the crack: ${W}^{*}(u,\Theta)$ in TF-PFM1 \eqref{TMF1b} and $W(u)$ in TF-PFM2 \eqref{TMF2b}.

Consequently, we established partial energy equality for TF-PFM1 (Theorem \ref{Theorem2:b}) and energy equality for TF-PFM2 (Theorem \ref{Theorem3}). From the viewpoint of energy consistency, both models are satisfactory, but TF-PFM2 is more energetically consistent than TF-PFM1. 

Based on the obtained numerical experiments, the following conclusions can be drawn.
\begin{enumerate}
	\item The thermoelastic coupling parameter $\delta$ in TF-PFM1 and TF-PFM2 enhances crack propagation (Figure \ref{fig1}).
	\item TF-PFM1 accelerates the crack speed more than TF-PFM2 (Figure \ref{Compare01}). On the other hand, the effect of the temperature gradient on the crack path in TF-PFM2 is larger than that in TF-PFM1 (Figure \ref{fig:crackpath3} - \ref{fig:Compare}).  
\end{enumerate}
The analytical and numerical comparisons between the two models are briefly summarized in Table \ref{tab:5}
\newcolumntype{A}{>{\raggedright\arraybackslash}p{18mm}}
\newcolumntype{D}{>{\raggedright\arraybackslash}p{15mm}}
\newcolumntype{L}{>{\raggedright\arraybackslash}p{46mm}}
\newcolumntype{C}{>{\raggedright\arraybackslash}p{15.75mm}}
\newcolumntype{E}{>{\raggedright\arraybackslash}p{13mm}}
\linespread{1.}{
\begin{table}[!h]
	\caption{Numerical comparison of TF-PFM1 and TF-PFM2}\label{tab:5}
	\begin{longtable}{C|LADE}
		\toprule
		 Models &
		Driving force &
		Energy consistency &
		Straight crack speed &
		Crack path \\
		\midrule
		\endhead
		TF-PFM1 &
		$W^{*}(u,\Theta)=\sigma^{*}[u,\Theta]:e^{*}[u,\Theta]$ &
		Partially satisfied &
		Faster &
		Less curved \\[20pt]
		
		TF-PFM2 &
		$W(u)= \sigma[u]:e[u]$ &
		Fully satisfied &
		Slower &
		More curved \\
		\midrule
		Remarks &
		$W^{*}(u,\Theta) >  W(u)$ (Figure \ref{fig:sigmaustar}) &
		Theorems \ref{Theorem2:b}, \& \ref{Theorem3} &
		Figure \ref{fig8} &
		Figures \ref{fig:crackpath3} \& \ref{fig:Compare} \\
		\bottomrule
	\end{longtable}
\end{table}
}

In this study, we did not consider the unilateral contact condition along the crack for the sake of simplicity. To further improve TF-PFM, the ideal unilateral condition for fracturing PFM \cite{Amor2009,Kimura2021} should be introduced in our PFM. 

\appendix
\gdef\thefigure{\arabic{figure}}    
\section{Weak forms for \eqref{SD}} \label{AppendixA}
The implicitly time-discretized problem \eqref{SD} is solved with the following boundary conditions. We set the initial temperature $\Theta^{0} := \Theta_{*}$ and set $\Theta^{-1} = \Theta_{*}$, which is a temperature of $t=-\Delta t$. For a given $\Theta^{k-1}$, the boundary value problem of $u^{k}$ is given as follows:
\begin{align}
&\left\{
\begin{array}{ll}
-{{\text{div}}{\sigma}^{*}[{u}^{k},\Theta^{k-1}]= 0 }& \mbox{in}~\Omega,\\
u^{k} = u_{D}(\cdot,t_{k}) & \text{on}~\Gamma_{D}^{u},\\
\sigma^{*}[u^{k},\Theta^{k-1}]n = 0 & \text{on}~\Gamma_{N}^{u},
\end{array}
\right. \qquad (k = 0,1,2,\cdots).  \label{Appendix1}
\end{align}
We define a weak form for \eqref{Appendix1} as
\begin{align}
\left\{
\begin{array}{l}
\displaystyle u^{k} \in V^{u}(u_{D}(\cdot,t_{k})),\\
\displaystyle\int_{\Omega} \sigma^{*}[u^{k},\Theta^{k-1}]: e[v]~dx = 0, \qquad(\mbox{for all}~ v \in V^{u}(0)),
\end{array}
\right.\label{WF-Appendix}
\end{align}
where $V^{u}(\cdot)$ is defined by \eqref{Original-Thermo-Elastic1}. The second equation of \eqref{WF-Appendix} is equivalent to
\begin{align}
& \int_{\Omega} \sigma[u^{k}]:e[v]~dx = \int_{\Omega} \Theta^{k-1}\mbox{div}v~dx.
\end{align}
Similarly,  for given $u^{k-1}$ and $u^{k}$, the boundary value problem of $\Theta^{k}$ is given as follows:
\begin{align}
&{\left\{
\begin{array}{ll}
\displaystyle\frac{\Theta^{k} - \Theta^{k-1}}{\Delta t} - \Delta\Theta^{k} + \delta\mbox{div}\left(\frac{u^{k} - u^{k-1}}{\Delta t}\right)= 0& \mbox{in}~\Omega,\\
\displaystyle\Theta^{k} = 0 & \text{on}~\Gamma_{D}^{\Theta},\\
\displaystyle\frac{\partial\Theta^{k}}{\partial n} = 0 & \text{on}~\Gamma_{N}^{\Theta},
\end{array}
\right.~ (k = 1,2,\cdots)}.  \label{Appendix2}
\end{align}
We define a weak form for \ref{Appendix2} as
\begin{align}
\left\{
\begin{array}{ll}
\displaystyle\Theta^{k} \in V^{\Theta},\\
\displaystyle\int_{\Omega}\left(\frac{\Theta^{k} - \Theta^{k-1}}{\Delta t}\right)\psi~dx + \int_{\Omega} \nabla\Theta^{k}\cdot\nabla\psi~dx \\
\displaystyle\hspace{50pt}+ \delta\int_{\Omega} \mbox{div}\left(\frac{u^{k} - u^{k-1}}{\Delta t}\right)\psi~dx = 0 ~(\text{for all}~ \psi \in V^{\Theta}(0)), 
\end{array}
\right. \label{WF-Appendix2} 
\end{align}
where $V^{\Theta} := \{\psi \in H^{1}(\Omega); ~\psi|_{\Gamma_{D}^{\Theta}}= 0\}$.
\begin{Propositionn}
We suppose that the $(d-1)$-dimensional volume of $\Gamma_{D}^{u}$ is positive. If $\Theta^{0} = \Theta_{*} \in L^{2}(\Omega)$ and $u_{D}(\cdot,t_{k}) \in H^{\frac{1}{2}}(\Gamma_{D}^{u})$ $(k=0,1,2,\cdots)$, then weak solutions $u^{k}$ $(k=0,1,2,\cdots)$ for \eqref{WF-Appendix} and $\Theta^{k}$ $(k=1,2,\cdots)$ for \eqref{WF-Appendix2} uniquely exist.
\normalfont\begin{Proof}
\normalfont At each time step, the unique solvabilities of \eqref{WF-Appendix} and \eqref{WF-Appendix2} follow from the Lax-Milgram theorem \cite{Girault1979,Duvaut1976}. More precisely, first we solve $u^{0}$ by \eqref{WF-Appendix}. Then, for $k=1,2,\cdots$, we can obtain $u^{k}$ by \eqref{WF-Appendix} and $\Theta^{k}$ by \eqref{WF-Appendix2}, sequentially. \qed
\end{Proof}
\end{Propositionn}
 
\section{Divergence of $u$ around the crack tip} \label{AppendixB}
We want to observe the contracting and expanding areas around the crack tip area. Here, we show an analytical solution for $\mbox{div}u$ around the crack tip. We consider Mode I as the type of loading; then, we analytically obtain the following crack tip displacement field:
\begin{align}
&u_{1} = \frac{K_{I}}{2\mu}\sqrt{\frac{r}{2\pi}}\cos\left(\frac{\theta}{2}\right)\left[\xi -1 + 2\sin^{2}\left(\frac{\theta}{2}\right)\right], \label{u1analytic}\\
&u_{2} = \frac{K_{I}}{2\mu}\sqrt{\frac{r}{2\pi}}\sin\left(\frac{\theta}{2}\right)\left[\xi +1 - 2\cos^{2}\left(\frac{\theta}{2}\right)\right]\label{u2analytic},
\end{align}
where $K_{I}$, $\mu$, $\xi = 3-4\nu_{P}$, and ($r$, $\theta$) are the Mode I stress intensity factor, Lam\'{e}'s constant, plane strain and polar coordinates for the crack tip, respectively.
  
\begin{figure}[!h]
	\centering 
	\includegraphics[trim=7cm 0cm 2cm 0cm, scale=0.15]{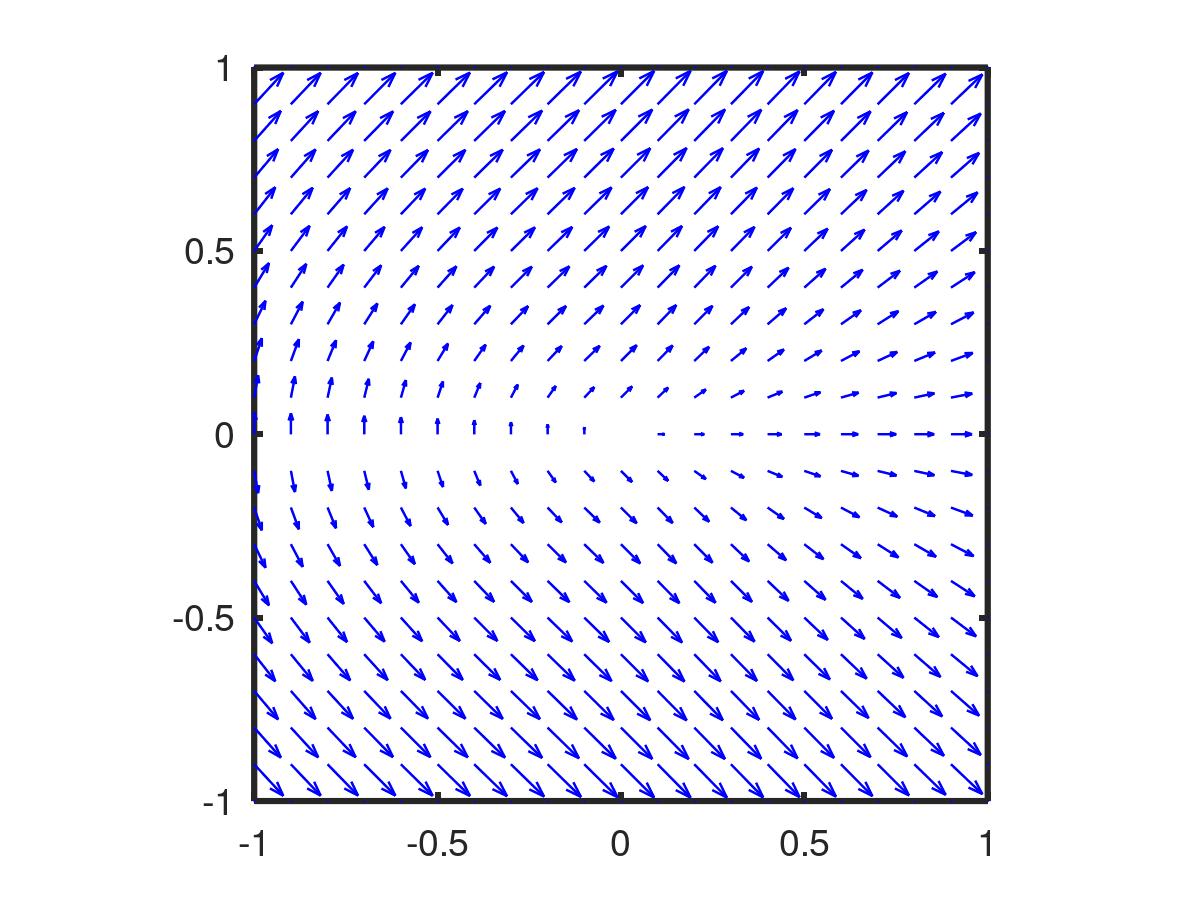}\quad
	\includegraphics[trim=7cm 0cm 2cm 0cm, scale=0.15]{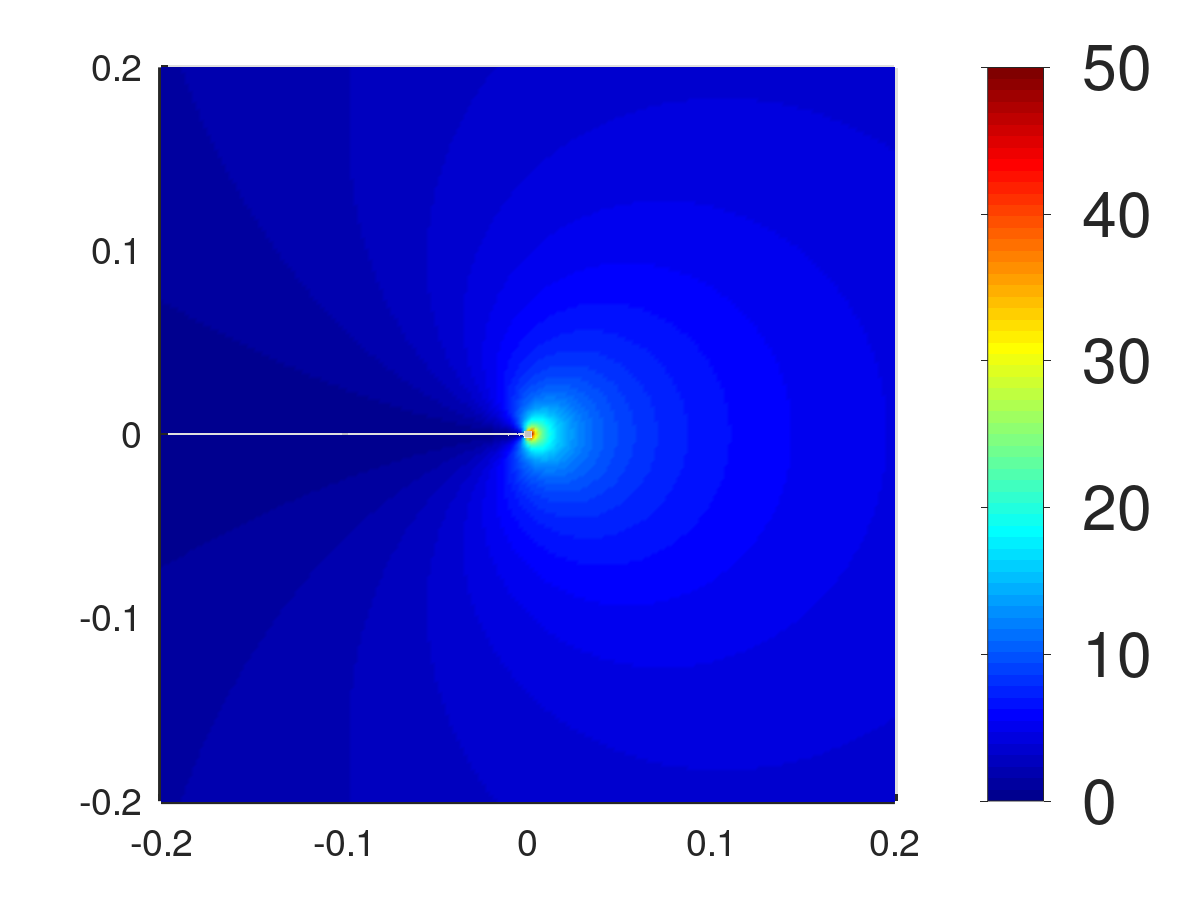}\par\medskip
	\includegraphics[trim=1cm 0cm 2cm 2cm, scale=0.15]{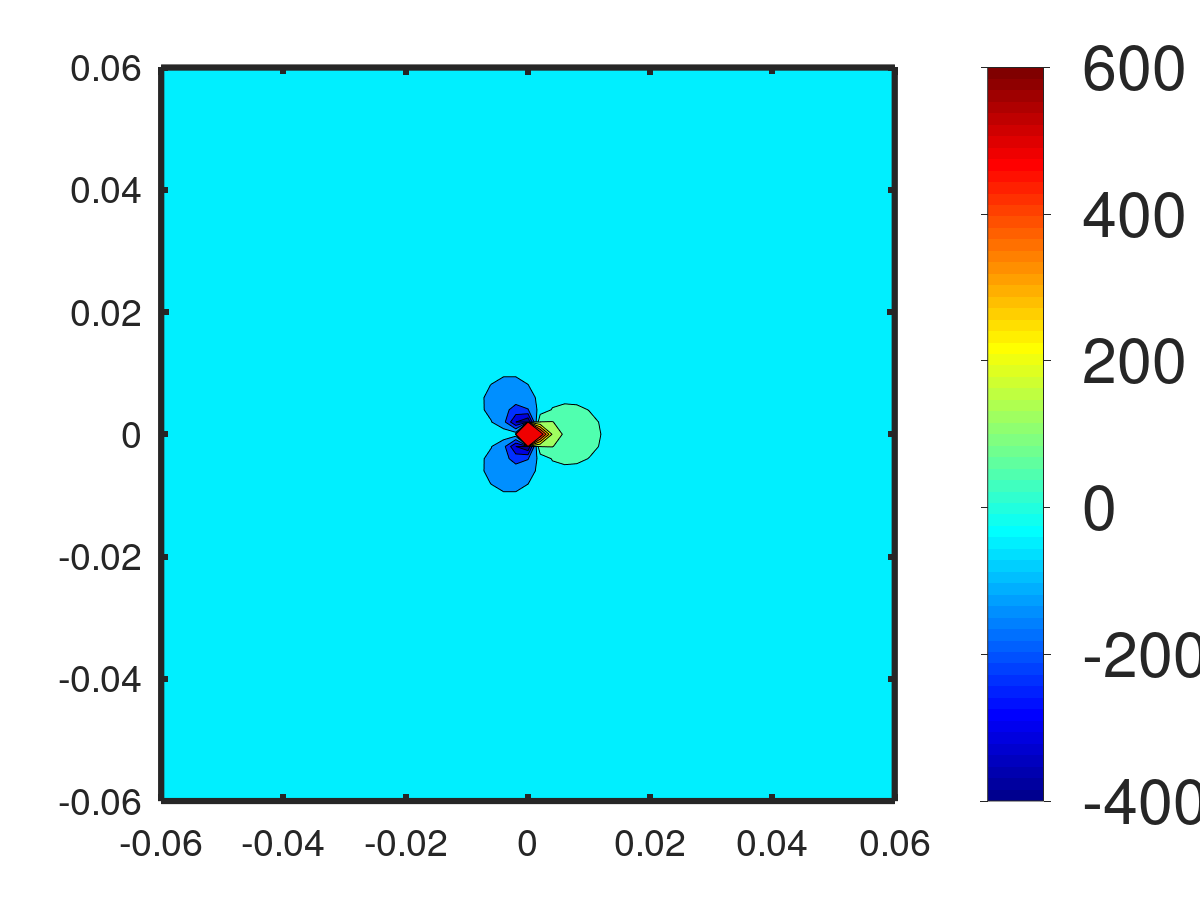}
	\caption{Profile of displacement $[u_{1},u_{2}]$ (left), $\mbox{div}u$ (right), and $\frac{\partial}{\partial x_{1}}(\mbox{div}{u})$ (bottom) around the crack under Mode I.}\label{fig:divergence}
\end{figure}

Now, we can calculate $\mbox{div}u$ as follows:
\begin{eqnarray}
\mbox{div}u & =& \left(\frac{\partial r}{\partial x_{1}}\frac{\partial}{\partial r} + \frac{\partial \theta}{\partial x_{1}}\frac{\partial}{\partial \theta}\right)u_{1} + \left(\frac{\partial r}{\partial x_{2}}\frac{\partial}{\partial r} + \frac{\partial \theta}{\partial x_{2}}\frac{\partial}{\partial \theta}\right)u_{2}\notag\\
&=& \frac{K_{I}}{2\mu\sqrt{2\pi r}} \Big(\xi\cos\Big(\frac{\theta}{2}\Big)- \cos\Big(\frac{\theta}{2}\Big)\Big)= \frac{K_{I}(\xi-1)}{2\mu\sqrt{2\pi r}} \cos\left(\frac{\theta}{2}\right).\label{divuan3}
\end{eqnarray}
Assume a crack is growing as
\begin{align}
& \Sigma(t) = \left\{(x_{1},0)^{T} \big|~ -\infty < x_{1} \leq v_{0}t\right\}. \notag
\end{align}
Then, we obtain the following displacement at time $t$
\begin{align}
&\tilde{u}(x,t) \approx u(x-v_{0}te_{1}),~ \text{where}~ e_{1} := (1,0)^{T}, \notag
\end{align}
and we also obtain $\mbox{div}u$ at time $t$
\begin{align}
\mbox{div}\tilde{u}(x,t) &= \mbox{div}u(x-v_{0}te_{1})\notag\\
\frac{\partial}{\partial t}\mbox{div}\tilde{u}(x,t)\big|_{t=0} & = -v_{0}\frac{\partial}{\partial x_{1}}\mbox{div}u = \frac{v_{0}K_{I}(\xi-1)}{4\mu\sqrt{2\pi r^{3}}}\cos\left(\frac{3\theta}{2}\right).\label{dxdivu}
\end{align}
Now, we set $E_{Y} = 1$, $\nu_{P} = 0.3$, $K_{I} = 5$, and $v_{0} = 0.05$, and then we obtain the displacement, $\mbox{div}u$ and the $\frac{\partial}{\partial x_{1}}\mbox{div}u$ profiles through \eqref{u1analytic}, \eqref{u2analytic}, \eqref{divuan3}, and \eqref{dxdivu}, respectively. From Figure \ref{fig:divergence}, a compressing area exists at the crack tip.

\section*{Acknowledgments}
This  work  was  supported  by  the MEXT  (the Ministry  of  Education,  Culture,  Sports,  Science,  and Technology) scholarship in Japan. This work was also partially supported by JSPS KAKENHI, grant numbers JP20H01812 and JP20KK0058.



\begin{thebibliography}{}
	\bibitem{Weilong2019}
	Ai, W., and Augarde, C. E., 2019. Thermoelastic fracture modelling in 2D by an adaptive cracking particle method without enrichment functions. International Journal of Mechanical Sciences. 160, 343-357. \href{https://doi.org/10.1016/j.ijmecsci.2019.06.033}{https://doi.org/10.1016/j.ijmecsci.2019.06.033}.
	
	\bibitem{Alfat2018}
	Alfat, S.,  Kimura, M., Firihu, M. Z., and Rahmat, 2018. Numerical investigation of shape domain effect to its elasticity and surface energy using adaptive finite element method, AIP Conference Proceedings. 1964(1), 020011. \href{https://doi.org/10.1063/1.5038293}{https://doi.org/10.1063/1.5038293}.

	\bibitem{Amor2009}
	Amor, H., Marigo, J.-J., and Maurini, C., 2009.
	Regularized formulation of the variational brittle fracture with unilateral contact: Numerical experiments. J. Mech. Phys. Solids. 57, 1209-1229. \href{https://doi.org/10.1016/j.jmps.2009.04.011}{https://doi.org/10.1016/j.jmps.2009.04.011}.
	
	\bibitem{Anderson2017}
	Anderson, T. L., 2017. Fracture Mechanics: Fundamentals and Applications, fourth ed. CRC press. 
	
	\bibitem{Biot1956}
	Biot, M. A., 1956. Thermoelasticity and irreversible thermodynamics. Journal of applied physics. 27(3),  240-253. \href{https://doi.org/10.1063/1.1722351}{https://doi.org/10.1063/1.1722351}.
	
	\bibitem{Bourdin2000}
	Bourdin, B., Francfort, G. A., and Marigo, J.-J., 2000. Numerical experiments in revisited brittle fracture. Journal of the Mechanics and Physics of Solids. 48(4), 797-826. \href{https://doi.org/10.1016/S0022-5096(99)00028-9}{https://doi.org/10.1016/S0022-5096(99)00028-9}. 
	
	\bibitem{Bourdin2007}
	Bourdin, B., 2007. Numerical implementation of the variational formulation of brittle fracture. Interfaces Free Boundaries. 9, 411-430. \href{doi: 10.4171/IFB/171}{doi: 10.4171/IFB/171}.
	
	\bibitem{Bourdin2008}
	Bourdin, B.,  Francfort, G. A.,  and Marigo, J.-J., 2008.
	The Variational Approach to Fracture. Journal of Elasticity. 91, 5–148. \href{https://doi.org/10.1007/s10659-007-9107-3}{https://doi.org/10.1007/s10659-007-9107-3}.
	
	\bibitem{Bourdin2014}
	Bourdin, B., Marigo, J.-J., Maurini, C., and Sicsic, P., 2014. Morphogenesis and propagation of complex cracks induced by thermal shocks. Physical review letters. 112(1), 014301. \href{https://doi.org/10.1103/PhysRevLett.112.014301}{https://doi.org/10.1103/PhysRevLett.112.014301}.
	
	\bibitem{Ciarlet2002}
	Ciarlet, P.G., 2002. The Finite Element Method for Elliptic Problems. Society for Industrial and Applied Mathematics.
	
	\bibitem{Duflot2008}
	Duflot, M., 2008. The extended finite element method in thermoelastic fracture mechanics. International Journal for Numerical Methods in Engineering. 74(5), 827-847. \href{https://doi.org/10.1002/nme.2197}{https://doi.org/10.1002/nme.2197}. 
	
	\bibitem{Duvaut1976}
	Duvaut, G., and Lions, J. L., 1976. Inequalities in Mechanics and Physics. Springer, Berlin, Heidelberg. 
	
	\bibitem{Dwivedi2018}
	Dwivedi, S. K., and Vishwakarma, M., 2018. Hydrogen embrittlement in different materials: a review. International Journal of Hydrogen Energy. 43(46), 21603-21616. \href{https://doi.org/10.1016/j.ijhydene.2018.09.201}{https://doi.org/10.1016/j.ijhydene.2018.09.201}. 
	
	
	\bibitem{Entezari2018}
	Entezari, A., Filippi, M., Carrera, E., and Kouchakzadeh, M. A., 2018. 3D dynamic coupled thermoelastic solution for constant thickness disks using refined 1D finite element models. Applied Mathematical Modelling. 60, 273-285. \href{https://doi.org/10.1016/j.apm.2018.03.015}{https://doi.org/10.1016/j.apm.2018.03.015}.
	
	\bibitem{Francfort1998}
	Francfort, G. A., and Marigo, J.-J., 1998. Revisiting brittle fracture as an energy minimization problem. Journal of the Mechanics and Physics of Solids. 46(8), 1319-1342. \href{https://doi.org/10.1016/S0022-5096(98)00034-9}{https://doi.org/10.1016/S0022-5096(98)00034-9}. 
	
	\bibitem{Freiman1984}
	Freiman, S. W., 1984. Effects of chemical environments on slow crack growth in glasses and ceramics. Journal of Geophysical Research: Solid Earth. 89(B6), 4072-4076. \href{https://doi.org/10.1029/JB089iB06p04072}{https://doi.org/10.1029/JB089iB06p04072}. 
	
	\bibitem{Hecht2012}
	Hecht, F., 2012. New development in FreeFem++. Journal of numerical mathematics. 20(3-4), 251-266. \href{https://doi.org/10.1515/jnum-2012-0013}{https://doi.org/10.1515/jnum-2012-0013}. 
	
	\bibitem{Gao2019}
	Gao, Y., and Oterkus, S., 2019. Ordinary state-based peridynamics modelling for fully coupled thermoelastic problems. Continuum Mechanics and Thermodynamics. 31(4), 907-937. \href{https://doi.org/10.1007/s00161-018-0691-1}{https://doi.org/10.1007/s00161-018-0691-1}. 
	
	\bibitem{Girault1979}
	Girault, V., and Raviart, P. A., 1979. Finite element approximation of the Navier-Stokes equations. Lecture Notes in Mathematics, Berlin Springer Verlag, 749.
	
	\bibitem{GreenLindsay1972}
	Green, A. E., and Lindsay, K. A., 1972.  Thermoelasticity. Journal of Elasticity. 2, 1-7. \href{https://doi.org/10.1007/BF00045689}{https://doi.org/10.1007/BF00045689}.
	
	\bibitem{GreenNaghdi1991}
	Green, A. E., and Naghdi, P. M., 1991. A re-examination of the basic postulates of thermomechanics. Proceedings of the Royal Society of London, Series A: Mathematical and Physical Sciences. 432(1885), 171-194. \href{https://doi.org/10.1098/rspa.1991.0012}{https://doi.org/10.1098/rspa.1991.0012}.
	
	\bibitem{Jaskowiec2017}
	Jaskowiec, J., 2017. A model for heat transfer in cohesive cracks. Computers \& Structures, 180, 89-103. \href{https://doi.org/10.1016/j.compstruc.2016.01.009}{https://doi.org/10.1016/j.compstruc.2016.01.009}. 
	
	\bibitem{Karma2001}
	Karma, A., Kessler, D. A., and Levine, H., 2001. Phase-field model of mode III dynamic fracture. Physical Review Letters. 87(4), 045501. \href{https://doi.org/10.1103/PhysRevLett.87.045501}{https://doi.org/10.1103/PhysRevLett.87.045501}. 
	
	\bibitem{Kimura2021}
	Kimura, M., Takaishi, T., Alfat, S., Nakano, T., and Tanaka, Y.,
	Irreversible phase field models for crack growth in industrial
	applications: thermal stress, viscoelasticity, hydrogen embrittlement.
	SN Applied Sciences. 3(781). \href{https://doi.org/10.1007/s42452-021-04593-6}{https://doi.org/10.1007/s42452-021-04593-6}.
	
	\bibitem{Kouchakzadeh2015}
	Kouchakzadeh, M. A., and Entezari, A., 2015. Analytical solution of classic coupled thermoelasticity problem in a rotating disk. Journal of Thermal Stresses. 38(11), 1267-1289. \href{https://doi.org/10.1080/01495739.2015.1073529}{https://doi.org/10.1080/01495739.2015.1073529}. 
	
	\bibitem{Kou2019}
	Kou, M., Liu, X., Tang, S., and Wang, Y., 2019. 3-D X-ray computed tomography on failure characteristics of rock-like materials under coupled hydro-mechanical loading. Theoretical and Applied Fracture Mechanics. 104, 102396. \href{https://doi.org/10.1016/j.tafmec.2019.102396}{https://doi.org/10.1016/j.tafmec.2019.102396}. 
	
	\bibitem{LordShulman1967}
	Lord, H. W., and Shulman, Y., 1967. A generalized dynamical theory of thermoelasticity. Journal of the Mechanics and Physics of Solids, 15(5), 299-309. \href{https://doi.org/10.1016/0022-5096(67)90024-5}{https://doi.org/10.1016/0022-5096(67)90024-5}.
	
	\bibitem{Louthan1972}
	Louthan Jr., M. R., Caskey Jr., G. R., Donovan, J. A., and Rawl Jr., D. E., 1972. Hydrogen embrittlement of metals. Materials Science and Engineering. 10, 357-368. \href{https://doi.org/10.1016/0025-5416(72)90109-7}{https://doi.org/10.1016/0025-5416(72)90109-7}.
	
	\bibitem{Mackin2002}
	Mackin, T. J., et al., 2002. Thermal cracking in disc brakes. Engineering Failure Analysis. 9(1), 63-76. \href{https://doi.org/10.1016/S1350-6307(00)00037-6}{https://doi.org/10.1016/S1350-6307(00)00037-6}. 
	
	\bibitem{Miehe2015}
	Miehe, C., Schaenzel, L. M., and Ulmer, H., 2015. Phase field modeling of fracture in multi-physics problems. Part I. Balance of crack surface and failure criteria for brittle crack propagation in thermo-elastic solids. Computer Methods in Applied Mechanics and Engineering. 294, 449-485. \href{https://doi.org/10.1016/j.cma.2014.11.016}{https://doi.org/10.1016/j.cma.2014.11.016}. 
	
	\bibitem{Nara2011}
	Nara, Y., Morimoto, K., Yoneda, T., Hiroyoshi, N., and Kaneko, K., 2011. Effects of humidity and temperature on subcritical crack growth in sandstone. International Journal of Solids and Structures. 48(7-8), 1130-1140. \href{https://doi.org/10.1016/j.ijsolstr.2010.12.019}{https://doi.org/10.1016/j.ijsolstr.2010.12.019}. 
	
	\bibitem{Nguyen2017}
	Nguyen, M. N., Bui, T. Q., Nguyen, N. T., and Truong, T. T., 2017. Simulation of dynamic and static thermoelastic fracture problems by extended nodal gradient finite elements. International Journal of Mechanical Sciences. 134, 370-386. \href{https://doi.org/10.1016/j.ijmecsci.2017.10.022}{https://doi.org/10.1016/j.ijmecsci.2017.10.022}. 
	
	\bibitem{Nguyen2016}
	Nguyen, T.T., Yvonnet, J., Bornert, M., et al., 2016. On the choice of parameters in the phase field method for simulating crack initiation with experimental validation. Int J Fract 197. 213–226. \href{https://doi.org/10.1007/s10704-016-0082-1}{https://doi.org/10.1007/s10704-016-0082-1}.
	
	\bibitem{Kimura2009}
	Takaishi, T., and Kimura, M., 2009. Phase field model for mode III crack growth in two dimensional elasticity. Kybernetika.  45(4),  605-614.
	
	\bibitem{Vivekanandan2020}
	Vivekanandan, A.,  and Ramesh, K., 2020. Study of crack interaction effects under thermal loading by digital photoelasticity and finite elements. Experimental Mechanics, 60(3), 295-316. \href{https://doi.org/10.1007/s11340-019-00561-9}{https://doi.org/10.1007/s11340-019-00561-9}. 
	
	\bibitem{Zheng2015}
	Zheng, B. J., Gao, X. W., Yang, K., Zhang, C. Z., 2015. A novel meshless local Petrov–Galerkin method for dynamic coupled thermoelasticity analysis under thermal and mechanical shock loading. Engineering Analysis with Boundary Elements. 60, 154-161. \href{https://doi.org/10.1016/j.enganabound.2014.12.001}{https://doi.org/10.1016/j.enganabound.2014.12.001}. 
	
	\bibitem{Zhou2011}
	Zhou, F. X., Li, S. R., and Lai, Y. M., 2011. Three-dimensional analysis for transient coupled thermoelastic response of a functionally graded rectangular plate. Journal of Sound and Vibration, 330(16), 3990-4001. \href{https://doi.org/10.1016/j.jsv.2011.03.015}{https://doi.org/10.1016/j.jsv.2011.03.015}. 
\end{thebibliography}
\end{document}